\documentclass[twoside,final,letterpaper,11pt]{amsart}

\usepackage{amsmath}               % AmSLaTeX
\usepackage{amsthm}
\usepackage{amscd}
\usepackage{times}

 \newfont{\Bb}{msbm10 scaled\magstephalf}

 \newsymbol\pitchfork 1374

 \title[Floer and quantum homologies]{An explicit isomorphism between Floer\\ homology\\
   and quantum homology}

 \author[G. C. Lu]{Guangcun Lu}
 \address{Department of Mathematics, Beijing Normal University\\
              Beijing 100875, P.R.China}
 \email{gclu@bnu.edu.cn}

 \thanks{The author was supported by the NNSF 19971045 and 10371007 of China.}

 \newtheorem{thm}{Theorem}[section]

 \theoremstyle{definition}

 \theoremstyle{remark}

 \date{November 18, 2000 /  Revised June 6, 2003}
\begin{document}

 \begin{abstract}
  We use Liu-Tian's virtual moduli cycle methods
  to construct detailedly the explicit isomorphism between  Floer
  homology and quantum homology for any closed symplectic manifold
  that was first outlined by Piunikhin, Salamon and Schwarz for the
  case of the semi-positive symplectic manifolds.
 \end{abstract}
 \maketitle

\section{Introduction}\label{sec:1}

\subsection{Background and motivation}\label{1.1}

It is one of the exciting mathematics achievements in the last few
years that the Floer and quantum homologies were established for
all closed symplectic manifolds (see \cite{FuO, LiT, LiuT1, R,
Sie} and \cite{HS2}). (Less general versions had been obtained
before by Floer \cite{F}, Hofer-Salamon \cite{HS1}, Ruan-Tian
\cite{RT1} and McDuff-Salamon \cite{McS}.) The purpose of this
paper is to construct detailedly an explicit ring isomorphism
between them. Such an isomorphism was first outlined by Piunikhin,
Salamon and Schwarz in the semi-positive case \cite{PSSc}. Our
argument is based on Liu-Tian's virtual cycle methods in
\cite{LiuT1}-\cite{LiuT3}. The isomorphism is necessary and
convenient for studies of some symplectic topology problems, e.x.,
the topology and geometry of the group ${\rm Ham}(M,\omega)$ of
Hamiltonian automorphisms of a symplectic manifold $(M,\omega)$.
Let $\widetilde G$ be the group of pairs $(g,\tilde g)$ consisting
of a smooth loop  $g: S^1\to {\rm Ham}(M,\omega)$ such that
$g(0)=Id$ and a lift $\tilde g:\tilde{\mathcal L}(M)\to
\tilde{\mathcal L}(M)$  of the action of $g$ to a covering of the
space ${\mathcal L}(M)$ of contractible loops in $M$ (see $\S 1.2$
below). In a beautiful paper \cite{Se} by Seidel, for every pair
$(g,\tilde g)\in\widetilde G$ there is assigned an automorphism
$HF_\ast(g,\tilde g)$ of the Floer homology $HF_\ast(M,\omega)$;
he constructed a homomorphism $q$ from $\widetilde G$ to the group
$QH_\ast(M,\omega)^\times$ of homogeneous even-dimensional
invertible elements of the quantum homology ring
$QH_\ast(M,\omega)$ and proved  his main result:
$$HF_\ast(g,\tilde g)(b)=\Psi^+(q(g,\tilde g))\ast_{PP}b$$
for any $(g,\tilde g)\in\widetilde G$ and $b\in
HF_\ast(M,\omega)$. Here  $\ast_{PP}$ and $\Psi^+$ are
 the `pair-of-pants' product in $HF_\ast(M,\omega)$
and the canonical isomorphism from $QH_\ast(M,\omega)$ to
$HF_\ast(M,\omega)$ constructed in \cite{PSSc} respectively. A key
step in the proof of his main result is Theorem 8.2 on the page
1080 of \cite{Se}, whose proof was based on the arguments of
\cite{PSSc}.

 Schwarz \cite{Sch3} defined and analyzed a bi-invariant metric on
${\rm Ham}(M,\omega)$ with the construction of such an explicit
isomorphism  on a closed symplectic manifold $(M,\omega)$ with
$c_1|_{\pi_2(M)}=\omega|_{\pi_2(M)}=0$.  Recently, Oh \cite{Oh}
obtained the corresponding results on arbitrary closed symplectic
manifolds. As pointed out in $\S 5.3$ of \cite{Oh} it would seem
more natural to use the Piunikhin-Salamon-Schwarz map in the
definition of his mini-max value function $\rho$. Entov \cite{En}
studied the relations between the {\rm K}-area for Hamiltonian
fibrations with a strongly semi-positive typical fiber
$(M,\omega)$ over a surface with boundary and the Hofer geometry
on the group ${\rm Ham}(M,\omega)$.  Such a ring isomorphism was
used  to obtain a key estimate in his work.

Since these applications used the Piunikhin-Salamon-Schwarz
isomorphism in \cite{PSSc} our detailed generalization to
arbitrary closed symplectic manifolds may be used to generalize
their results to the desired forms  more directly and
conveniently. Moreover the method to construct the ring
isomorphism has actually more uses than the isomorphism itself
because  not only the isomorphism itself but also the map of
Piunikhin-Salamon-Schwarz'type in the chain level were used in
some applications.  The construction of another ring isomorphism
was given by Liu-Tian \cite{LiuT3} (a less general version was
announced before by Ruan-Tian \cite{RT2}). Without doubt different
construction methods of the ring isomorphisms between Floer
homology and quantum homology have respective advantages in the
studies of different symplectic topology problems.

\subsection{Outline and the main result}\label{1.2}

 For a smooth non-degenerate
time-dependent function $H: M\times\mbox{\Bb R}/\mbox{\Bb
Z}\to\mbox{\Bb R}$ one may associate a family of the Hamiltonian
vector fields $X_{H_t}$  by $\omega(X_{H_t}, \cdot)=-dH_t$ for
$t\in\mbox{\Bb R}$ and $H_t(\cdot)=H(t, \cdot)$. Let ${\mathcal
P}(H)$ be the set of all contractible $1$-periodic solutions of
the Hamiltonian differential equation: $\dot x(t)=X_{H_t}(x(t))$.
 Denote by ${\mathcal  J}(M,\omega)$ the space of all almost complex
structures  compatible with $\omega$. It determines a unique the
first Chern class $c_1=c_1(TM, J)\in H^2(M, \mbox{\Bb Z})$ via any
$J\in{\mathcal  J}(M,\omega)$(\cite{Gr}). Let
$\phi_{c_1},\phi_\omega: H^S_2(M)\to\mbox{\Bb R}$ be the
homomorphisms by evaluations of $c_1$ and $\omega$ respectively.
Here $H^S_2(M)$ denotes the image of $\pi_2(M)$ in
$H_2(M;\mbox{\Bb Z})$ under the Hurewicz homomorphism modulo
torsion. As usual let ${\mathcal  L}(M)$ be the set of all
contractible loops $x\in C^\infty(S^1, M)$. Consider a pair $(x,
v)$ consisting of $x\in{\mathcal L}(M)$ and a disk $v$ bounding
$x$. Such two pairs  $(x, v)$ and $(y, w)$ are called equivalent
if $x=y$ and $\phi_{c_1},\phi_\omega$ both take zero value on
$v\sharp(-w)$. Denote by $[x, v]$ the equivalence class of a pair
$(x, v)$ and by $\widetilde{{\mathcal  L}(M)}$ the set of all such
equivalence classes. Then the latter is a cover space of
${\mathcal  L}(M)$ with the covering transformation group
$\Gamma=H^S_2(M)/(\rm{ker}\phi_{c_1}\cap\rm{ker}\phi_\omega)$. Its
action is given by  $A\cdot[x, v]=[x, A\sharp v]$ for any
$A\in\Gamma$, where $A\sharp v$ is understood as the connected sum
of any representative  of $A$ in  $\pi_2(M)$ with $v$. In this
paper we shall denote $[x, v]$ by $\tilde x$ and $[x, v\sharp A]$
by $\tilde x\sharp A$ for $A\in\Gamma$ if it is not necessary to
point out the bounding disk $v$ and no confusion occurs.
 Let $\tilde{\mathcal  P}(H)$ be the lifting of ${\mathcal
P}(H)$ in the space $\widetilde{{\mathcal  L}(M)}$. It is exactly
the critical set of the functional
$$\left.\begin{array}{ll}
{\mathcal  F}_H: \widetilde{{\mathcal  L}(M)}\to\mbox{\Bb R},\;
[x, v]\mapsto -\int_D v^*\omega+ \int^1_0 H(t, x(t))dt
\end{array}\right.$$
 on $\widetilde{{\mathcal  L}(M)}$.
 Given  $\tilde x^\pm=[x^{\pm}, v^{\pm}]\in{\tilde{\mathcal  P}}(H)$
 denote by
$${\mathcal  M}(\tilde x^-, \tilde x^+;  H, J)\leqno(1.1)$$
the space of all connecting
 trajectories $u:\mbox{\Bb R}\times S^1\to M$ satisfying the equation
 ${\bar\partial}_{J, H}u=\partial_s u + J(u)\partial_t u+
\nabla{H_t}(u)=0$ with boundary conditions
$\lim_{s\to\pm\infty}u(s,t)\\=x^\pm$ and $[x^+, v^-\sharp u]=[x^+,
v^+]$. After making a small generic perturbation of $H_t$ outside
some small neighborhood of the graph of the elements of ${\mathcal
P}(H)$ one may assume that for any two $\tilde x^\pm\in
{\tilde{\mathcal  P}}(H)$ the space ${\mathcal M}(\tilde
x^-,\tilde x^+; H, J)$ is either an empty set or a manifold of
dimension $\mu(\tilde x^-)- \mu(\tilde x^+)$ (\cite{F}).  Here
$\mu(\tilde x)$ is the Conley-Zehnder index of $\tilde x$
(\cite{SZ}).
 Denote by ${\tilde{\mathcal  P}}_k(H):=\{\tilde x\in {\tilde{\mathcal  P}}(H)\,|\,\mu(\tilde x)=k\}$.
Consider the chain complex whose $k$-th chain group $C_k(H,
J;\mbox{\Bb Q})$ consists of all formal sums $\sum\xi_{\tilde
x}\cdot\tilde x$ with $\xi_{\tilde x}\in\mbox{\Bb Q}$ and $\tilde
x\in{\tilde{\mathcal  P}}_k(H)$ such that the set $\{\tilde x\in
{\tilde{\mathcal  P}}_k(H)\,|\, \xi_{\tilde x}\ne 0,\, {\mathcal
F}_H(\tilde x)>c\}$ is finite for any $c\in\mbox{\Bb R}$. Then
$C_\ast (H, J;\mbox{\Bb Q})={\oplus}_k C_k(H, J;\mbox{\Bb Q})$ is
a graded $\mbox{\Bb Q}$-space of infinite dimension. However its
dimension as a module  over the Novikov ring
$\Lambda_{\omega}=\Lambda_{\omega}(\mbox{\Bb Q})$ is finite. Here
$\Lambda_{\omega}(\mbox{\Bb Q})$ is the collection of all formal
sums $\lambda=\sum\lambda_A\cdot e^A$ with $\lambda_A\in\mbox{\Bb
Q}$  such that the set $\{A\in\Gamma\,|\,\lambda_A\ne 0,
\omega(A)< c\}$ is finite for any $c\in\mbox{\Bb R}$. Its action
on $C_\ast (H, J;\mbox{\Bb Q})$ is defined  by
 $(\lambda\ast\xi)_{\tilde x}=\sum_{A\in\Gamma}\lambda_A
 \xi_{(-A)\cdot\tilde x}$
and the rank of $C_\ast(H, J;\mbox{\Bb Q})$ over
$\Lambda_{\omega}$ is equal to $\sharp{\mathcal  P}(H)$. When $(M,
\omega)$ is a monotone symplectic manifold and $\mu(\tilde x^-)-
\mu(\tilde x^+)=1$  Floer proved that the manifold in (1.1) is
compact and thus first established his Floer homology theory
\cite{F}. Later on his arguments were generalized to the
semi-positive case by Hofer-Salamon \cite{HS1} and Ono \cite{O},
and the case of the product of semi-positive symplectic manifolds
by author \cite{Lu1}. Note that the space in (1.1) is not compact
in general case. It is this noncompactness that impedes the
establishment of Floer homology theory on all closed symplectic
manifolds. Let us outline Liu-Tian's method to overcome this
difficulty since we shall choose their method to realize our
program. Replacing the space in (1.1) they considered
$\overline{\mathcal  M}(\tilde x^-, \tilde x^+;  J, H)$ the space
of the equivalence classes of all $(J, H)$-stable trajectories
from $\tilde x^-$ to $\tilde x^+$(cf. Def.2.1), and used it to
 construct a suitable relative virtual moduli cycle
$C(\overline{\mathcal  M}^\nu(\tilde x^-, \tilde x^+))$ of
dimension $\mu(\tilde x^-)-\mu(\tilde x^+)-1$(see \S2.2.) If
$\mu(\tilde x^-)-\mu(\tilde x^+)=1$ the virtual moduli cycle may
determine a rational number $\sharp\bigl(C(\overline{\mathcal
M}^\nu(\tilde x^-, \tilde x^+))\bigr)$ ( see \cite{LiuT1}). Then
for each $\xi=\sum_{\tilde x}\xi_{\tilde x}\tilde x$ in $C_k(H,
J;\mbox{\Bb Q})$ they defined
$$
{\partial}^F_k\xi=\sum_{\mu(\tilde y)=k-1} \biggl[\sum_{\mu(\tilde
x)=k}\sharp\bigl(C(\overline{\mathcal M}^\nu(\tilde x,\tilde y))
\bigr)\cdot\xi_{\tilde x}\biggr]\tilde y.\leqno(1.2)
$$
and proved
it to be indeed a boundary operator. Let $HF_\ast(M, \omega; H,
J,\nu;\mbox{\Bb Q})$ be the homology of the above chain complex.
Using Floer's deformation ideas they proved that this homology is
invariant under deformations and also isomorphic to
$H_\ast(M,\mbox{\Bb Q})\otimes\Lambda_\omega$, i. e., the quantum
homology of $M$.

For the construction of our isomorphism let us firstly fix a Morse
function $h_0$ on $M$ and a small open neighborhood ${\mathcal
O}(h_0)$ of it in $C^\infty(M)$ such that

 (i) for some $\epsilon>0$ any two different points $a$ and $b$ of ${\rm
Crit}(h_0)=\{a_1,\cdots, a_m\}$ have a distance $d(a,
b)>4\epsilon$ with respect to some distance $d$ on $M$;

(ii) any two different critical points $a^h$ and $b^h$ of
$h\in{\mathcal O}(h_0)$ have a distance $d(a^h, b^h)>3\epsilon$;

(iii) each $h\in{\mathcal  O}(h_0)$ has a unique critical point
$a^h_i$ in each ball $B_d(a_i, \epsilon)=\{c\in M\,|\, d(a_i,
c)<\epsilon\}$ and no other critical points (hence $\sharp{\rm
Crit}(h)=m$);

(iv) for any $h\in{\mathcal O}(h_0)$ the Morse index
$\mu(a^h_i)=\mu(a_i)$, $i=1,\cdots, m$;

(v) the function\vspace{-1mm}
$${\mathcal  O}(h_0)\to B_d(a_1, \epsilon)\times\cdots\times
B_d(a_m, \epsilon),\; h\mapsto (a^h_1, \cdots, a^h_m)$$
 is a smooth surjective map.\vspace{1mm}

 Take $h\in{\mathcal O}(h_0)$ and a Riemannian metric $g$ on $M$ such that $(h, g)$ is
a Morse-Smale pair. As in \cite{PSSc} we may use  the solutions
$\gamma:\mbox{\Bb R}\to M$ of
$$\dot\gamma(s)=-\nabla^gh(\gamma(s))\leqno(1.3)$$
to construct a chain complex expression of the quantum homology
$H_\ast(M, \mbox{\Bb Q})\otimes\Lambda_\omega$ as follows:
 For every integer $k$ let us denote by
$$QC_k(M,\omega; h, g;\mbox{\Bb Q})\leqno(1.4)$$
the set of all formal sums $\zeta=\sum_{\mu(\langle a,
A\rangle)=k}\zeta_{(a, A)} a\oplus A$ such that $\{(a, A)\in ({\rm
Crit}(h)\times\Gamma)_k\,|\, \zeta(a, A)\ne 0,\;
h(a)-\phi_\omega(A)> c\,\}$ is a finite set for all $c\in\mbox{\Bb
R}$. Here $({\rm Crit}(h)\times\Gamma)_k:=\{(a, A)\in {\rm
Crit}(h)\times\Gamma\;|\; \mu(\langle a,
A\rangle):=\mu(a)-2c_1(A)=k\}$. The action of $\Lambda_{\omega}$
on $QC_\ast(M,\omega; h, g;\mbox{\Bb Q})$ is given by
$$\lambda\star\zeta=\sum_{(b, B)\in {\rm Crit}(h)\times\Gamma}
\biggl(\sum_{2c_1(A)=\mu(\langle b,
B\rangle)-k}\lambda_{\alpha\oplus A} \zeta_{\langle{(-\alpha)\cdot
b},B-A\rangle}\biggr) \langle b, B\rangle.\leqno(1.5)$$
 The boundary operator $\partial^Q_k:
QC_k(M,\omega; h, g;\mbox{\Bb Q})\to {QC}_{k-1}(M,\omega; h,
g;\mbox{\Bb Q})$ is given by
$$\left.\begin{array}{ll}\partial^Q_k(\langle a, A\rangle)=\sum_{\mu(b)=\mu(a)-1}
n(a, b)\langle b, A\rangle,\end{array}\right.\leqno(1.6)$$
 where $n(a, b)$ is the oriented number of the solutions of (1.3) from $a$ to
$b$. It is easily checked that $\;\lambda\star\zeta\in
QC_\ast(M,\omega; h, g; \mbox{\Bb Q})$ and that the latter is a
graded vector space over $\Lambda_{\omega}(\mbox{\Bb Q})$
according to the multiplication defined in (1.5), and that
$\partial^Q$ is a boundary operator and also
$\Lambda_{\omega}(\mbox{\Bb Q})$-linear with respect to the
multiplication. Consequently, $(QC_\ast(M,\omega; h, g;\mbox{\Bb
Q}),\partial^Q)$ is a chain complex. Let us denote its homology
by\vspace{-1mm}
$$QH_\ast(h, g;\mbox{\Bb Q}):
=H_\ast(QC_\ast(M,\omega; h, g;\mbox{\Bb Q}),\;\;\partial^Q).$$ It
is easy to derive from \cite{Sch4} that there exists an explicit
graded $\Lambda_\omega$-module isomorphism between $QH_\ast(h,
g;\mbox{\Bb Q})$ and $H_\ast(M;\mbox{\Bb
Q})\otimes\Lambda_\omega$(\cite{Sch2}).

To construct an explicit $\Lambda_\omega$-module isomorphism
between $QH_\ast(h, g;\mbox{\Bb Q})$ and $HF_\ast(M,\omega; H,
J,\nu; \mbox{\Bb Q})$ we shall associated two rational numbers
$n_+^{\nu^+}(a, \tilde x\sharp(-A))$ and $n_-^{\nu^-}(a, \tilde
x\sharp(-A))$ in (3.2) to both $\langle a, A\rangle\in{\rm
Crit}(h)\times\Gamma$ and $\tilde x\in{\tilde{\mathcal  P}}(H)$,
and then mimic \cite{PSSc} to construct formally two maps on the
levels of chains
$$\left.\begin{array}{ll}\Phi(\langle a, A\rangle)=\sum_{\mu(\tilde x)=\mu(\langle a, A\rangle)}
n_+^{\nu^+}(a, \tilde x\sharp(-A))\tilde
x,\end{array}\right.\leqno(1.7)$$\vspace{-2mm}
$$\left.\begin{array}{ll}\Psi(\tilde x)=\sum_{\mu(\langle a,  A\rangle)=\mu(\tilde x)}
n_-^{\nu^-}(a, \tilde x\sharp(-A))\langle a,
A\rangle.\end{array}\right.\leqno(1.8)$$ Indeed, in Remark 3.7 we
shall show that $\Phi$ and $\Psi$ are $\Lambda_\omega$-module
chain homomorphisms from $QC_\ast(M, \omega; h, g;\mbox{\Bb Q})$
to $C_\ast(H, J,\nu;\mbox{\Bb Q})$ and from $C_\ast(H,
J,\nu;\mbox{\Bb Q})$ to $QC_\ast(M,\omega; h, g;\mbox{\Bb Q})$
respectively. Our main result is:

\begin{thm}\label{thm:1.1}  $\Phi$ induces a
$\Lambda_\omega$-module isomorphism
$$\Phi_\ast:QH_\ast(h, g;\mbox{\Bb Q})\to
HF_\ast(M,\omega; H, J,\nu;\mbox{\Bb Q})$$ with an inverse
$\Psi_\ast$.
\end{thm}

\noindent{\bf Remark 1.2.}\hspace{1mm} As concluding remarks we
point out that Theorems 3.1, 3.7 and 5.1 in \cite{PSSc} may easily
be extended to any closed symplectic manifold. Such an extension
of Theorem 3.1 was actually carried out in \cite{LiuT3}. Following
the lines in \cite{PSSc} and combing the methods in \cite{LiuT3}
with ones in this paper we easily complete the extension of
Theorem 5.1 in \cite{PSSc} to arbitrary closed symplectic
manifolds (in fact, a long exercise). As a consequence we get that
the isomorphism in Theorem 1.1 is also the ring
isomorphism.\vspace{1.5mm}

Section 2 introduces the moduli spaces of stable disks and
constructs the virtual moduli cycles associated with them such
that they are compatible with Liu-Tian's relative ones. Section 3
deals with the intersections of these virtual moduli cycles with
stable and unstable manifolds. The proof of Theorem~\ref{thm:1.1}
is completed in \S4.

 {\bf Acknowledgements}.\hspace{2mm}The author is
very grateful to Professors Gang Tian and Dietmar Salamon for
their helps in my understanding for Floer homology in past years.
He would also like to thank Professors Dusa McDuff, Yong-Geun Oh,
Yuli B. Rudyak, Matthias Schwarz and Claude Viterbo for sending me
their preprints on Floer homology and Arnold conjecture. Finally,
he is also grateful to referees for their many very good
improvement suggestions.

\section{The disk solution spaces and virtual moduli cycles}\label{sec:2}

\subsection{Moduli space of stable disks}\label{2.1}
We begin with the disk solution spaces introduced in \cite{PSSc}.
 For $J\in{\mathcal  J}(M,\omega)$  and $[x, v]\in\tilde{\mathcal  P}(H)$ let
${\mathcal  M}_+([x, v]; H, J)$ be the set of all smooth maps $u:
\mbox{\Bb R}\times S^1\to M$ such that
$$\partial_s u(s,t)+ J(u)
(\partial_t u-\beta_+(s)X_{H}(t, u))
 =0,\leqno(2.1)$$\vspace{-4mm}
$$\left.\begin{array}{ll}
u(+\infty)=x\quad{\rm and}\quad
E_+(u):=\int^{\infty}_{-\infty}\int^1_0 |\partial_s
u|^2_{g_J}dsdt<+\infty,\end{array}\right.$$\vspace{-4mm}
 $$u\sharp(-v): S^2\to M\;{\rm represents\;a\;
torsion\;homology\;class\;in}\; H_2(M;\mbox{\Bb
Z}).\leqno(2.2)$$\vspace{-1mm}
 Here a smooth cut-off function
$\beta_+:\mbox{\Bb R}\to [0,1]$ is given by
\begin{eqnarray*}
\beta_+(s)=\left\{\begin{array}{ll}0 \;&{\rm as}\; s\le 0\\
1\;&{\rm as}\;s\ge 1\end{array}\right.,\;2>\beta_+^\prime(s)>0,\;{\rm
for}\;0<s<1.
\end{eqnarray*}
Since $E_+(u)<+\infty$,  $u$ extends over the end of $s=-\infty$
by removable singularity theorem. So the connected union
$u\sharp(-v)$ in (2.2) is well-defined. If
 ${\mathcal  M}_+([x, v]; H, J)\ne\emptyset$ its virtual dimension is
${\dim M}-\mu([x, v])$.

Correspondingly, for $[x, v]\in\tilde{\mathcal  P}(H)$ we denote
by ${\mathcal  M}_-([x, v]; H, J)$ the set of all smooth maps
$u:\mbox{\Bb R}\times S^1\to M$ such that
$$\partial_s u(s,t)+ J(u)
(\partial_t u-\beta_+(-s)X_{H}(t, u))
 =0,\leqno(2.3)$$\vspace{-4mm}
$$\left.\begin{array}{ll}u(-\infty)=x\quad{\rm and}\quad
E_-(u):=\int^{\infty}_{-\infty}\int^1_0|\partial_s
u|^2_{g_J}dsdt<+\infty,\end{array}\right.$$\vspace{-4mm}
 $$v\sharp u:  S^2\to M\;{\rm represents\;a\;torsion\;homology\;class\;in}\;
 H_2(M;\mbox{\Bb Z}).\leqno(2.4)$$
Similarly, if ${\mathcal  M}_-([x, v]; H, J)\ne\emptyset$ its
virtual dimension is $\mu([x, v])$. As above we here have extended
$u$ over the end of $s=+\infty$. Recall that: \vspace{2mm}

\noindent{\bf Definition 2.1}(\cite{LiuT1}).\hspace{2mm} Let
$(\Sigma,\underline l)$ be a semistable $\mathcal  F$-curve  with
$z_-=z_1,\cdots, \\z_{N_p+1}=z_+$, as those double points
connecting the principal components(cf. Def.3.1 in [LiuT2]). A
continuous map $f:\Sigma\setminus\{z_1,\cdots,z_{N_p+1}\}\to M$ is
called a {\bf stable $(J, H)$-map} if there exist $[x_i,
v_i]\in\tilde{\mathcal P}(H)$, $i=1,\cdots,N_p+1$, such that:
\begin{enumerate}
\item[(1)] on each principal component $P_i$ with cylindrical
coordinate $(s,t)$ (obtained by the identification
$(P_i\setminus\{z_i, z_{i+1}\};l_i)\equiv (\mbox{\Bb R}\times
S^1;\{t=0\})$), $f^P_i=f|_{P_i-\{z_i,z_{i+1}\}}$ satisfies:
{\begin{enumerate}
 \item[(i)] $\partial_s f^P_i+
 J(f^P_i)\partial_t f^P_i+\nabla H(f^P_i)
=0$, and
 \item[(ii)] $\lim_{s\to -\infty}f^P_i(s,t)=x_i(t)$ and
 $\lim_{s\to +\infty}f^P_i(s,t)=x_{i+1}(t)$.
\end{enumerate}}
\item[(2)] The restriction $f^B_j$ of $f$ to
each bubble component $B_j$ is $J$-holomorphic; \item[(3)]
$[v_{i+1}]=[v_i]+[f^P_i]+\sum_{j}[f^B_{i,j}]$ as relative homology
class of $(M, x_{i+1})$,
 where the domain of $f^B_{i,j}$ may be joined to $P_i$ by a
chain of the bubble components not intersecting with other
principal components;
\item[(4)] All homotopically trivial  principal components
or homologically trivial  bubble components are not free.
 \end{enumerate}\vspace{2mm}

 The equivalence class of  $f=(f,\Sigma,\underline l)$ is denoted by
 $\langle f,\Sigma,\underline l\rangle$ or simply $\langle f\rangle$( see \cite{LiuT1}
 for its definition).  To compactify the disk solution spaces we introduce: \vspace{2mm}

\noindent{\bf Definition 2.2.}\hspace{2mm}  Given a $[x, v]\in
\tilde{\mathcal  P}(H)$ and a  semistable ${\mathcal F}$-curve
$(\Sigma,\underline l)$,  a continuous map
$f:\Sigma\setminus\{z_2,\cdots,z_{N_p+1}\}\to M$ is called a {\bf
stable  $(J, H)_+$-disk} with cap $[x, v]$ if there exist $[x_i,
v_i]\in{\tilde{\mathcal  P}}(H)$, $i=2,\cdots,N_p+1$, with $[x,
v]=[x_{N_p+1},  v_{N_p+1}]$, such that:
\begin{enumerate}
\item[(1)] on each principal component $P_i$ $(i>1)$ with
cylindrical coordinate $(s,t)$, $f^P_i=f|_{P_i-\{z_i,z_{i+1}\}}$
satisfies {\rm (i) (ii)} in Definition 2.1(1), but $f^P_{1}$ does
(2.1) and $f_{1}^P(+\infty)=x_2$ and $E_+(f_1^P)<+\infty$.
\item[(2)] If $N_p>1$, for each $i>1$ as relative homology classes
of $(M, x_{i+1})$, $[v_{i+1}]=[v_i]+[f^P_i]+\sum_{j}[f^B_{i,j}]$
and $[v_{2}]=[f^P_1]+ \sum_{j}[f^B_{1,j}]$.
 \item[(3)] All
requirements for the bubble components $f^B_{i,j}$ in Definition
2.1(2)(4) still hold. Moreover, all homotopically trivial
principal components $f^P_i\;(i>1)$ are not free and do not appear
in the next way.
 \end{enumerate}\vspace{2mm}

\noindent{\bf Remark 2.3.}\quad
 Actually, $f^P_1$ may not be constant. So, if $f$ has at
least two components then there exist at least two {\it
nonconstant} components. This also holds for $f^P_{N_p}$ in the
following Definition 2.4.  The energy of such a map is defined by
 $$
  E_+(f)=\sum_{i}\int\!\!\int_{\mbox{\Bb R}\times S^1}
 |\partial_s f^P_i|^2_{g_J} ds dt +
\sum_{i,j}\int_{B_{i,j}}(f^B_{i,j})^\ast\omega.
\leqno(2.5)$$

\noindent{\bf Definition 2.4.}\hspace{2mm} Given a $[x, v]\in
\tilde{\mathcal P}(H)$ and a $(\Sigma,\underline l)$ as before, a
continuous map $f:\Sigma\setminus\{z_1, \cdots, z_{N_p}\}\to M$ is
called a {\bf stable $(J, H)_-$-disk} with cap $[x, v]$  if there
exist $[x_i, v_i]\in{\tilde{\mathcal P}}(H)$, $i=1,\cdots, N_p$,
with $[x, v]=[x_1, v_1]$, such that:
\begin{enumerate}
\item[(1)] On each principal component $P_i$($i<N_p$) with
cylindrical coordinate $(s,t)$, $f^P_i=f|_{P_i-\{z_i,z_{i+1}\}}$
satisfies {\rm (i) (ii)} in Definition 2.2(1), but $f^P_{N_p}$
does (2.3) and $f^P_{N_p}(-\infty)=x_{N_p}$ and
$E_-(f^P_{N_p})<+\infty$.
 \item[(2)] If $N_p> 1$, for each $1\le
i\le N_p-1$, as relative
 homology class of $(M, x_{i+1})$,
$[v_i]=[v_{i+1}]+[f^P_i]+\sum_{j}[f^B_{i,j}]$ and $
[v_{N_p}]=[f^P_{N_p}]+\sum_{j}[f^B_{N_p,j}]$.
 \item[(3)] All
assertions for the bubble components $f^B_{i,j}$ in Definition
2.1(2)(4) still hold. Moreover, all homotopically trivial
principal components $f^P_i\;(i<N_p)$ are not free and do not
appear in the next way.
 \end{enumerate}\vspace{2mm}

We still define the energy $E_-(f)$ of such a map by the right
side of (2.5). As before we may define their equivalence classes.
 Let us denote by  $\langle f,\Sigma,\underline l\rangle$ or simply
 $\langle f\rangle$ the  equivalence class of  $(f,\Sigma,\underline l)$.
The energy of $\langle f\rangle$ is defined by that of any
representative of it. The direct computation shows
that\vspace{-1mm}
$$E_\pm(f_\pm)\le \mp{\mathcal  F}_H([x, v])+
\max|H|.\leqno(2.6)$$
 The notions of the (effective) dual graph for the
 stable $(J, H)_\pm$-disks may also be defined with
the same way as in [LiuT1]. Let
 $\overline{\mathcal M}_\pm(\tilde x; H, J)$ be the spaces of the equivalence classes of all
$(J, H)_\pm$-stable disks with cap $\tilde x$ respectively.
 By (2.6),
 $$
 \pm\int_{D^2}v^\ast\omega\ge-2\max |H|\;\;{\rm as}\;\;
\overline{\mathcal  M}_\pm([x, v]; H, J)\ne\emptyset.
\leqno(2.7)$$
 One may equip the weak $C^\infty$ topology on
$\overline{\mathcal  M}_\pm(\tilde x; H, J)$ according to the
definition given by (i)(ii)(iii) above Proposition 4.1 of
\cite{LiuT1} unless we allow the compact set $K$ in (ii) to be
able to contain the  double point $z_-$ (resp. $z_+$) on the chain
of principal components of the domain of $\langle
u_\infty\rangle\in\overline{\mathcal M}_+(\tilde x; H, J)$ (resp.
$\overline{\mathcal  M}_-(\tilde x; H, J)$). Carefully checking
proof of Proposition 4.1 in \cite{LiuT1} we have:\vspace{2mm}

\noindent{\bf Proposition 2.5.}\hspace{2mm}{\it The spaces
$\overline{\mathcal  M}_\pm(\tilde x; H, J)\supseteq {\mathcal
M}_\pm(\tilde x; H, J)$ are compact and Hausdorff with respect to
the weak $C^\infty$-topology.}\vspace{2mm}

 Notice that we have two natural
continuous maps
$${\rm EV}_\pm(\tilde x):
\overline{\mathcal  M}_\pm(\tilde x; H, J)\to M,\; \langle
f_\pm\rangle\mapsto f_\pm(z_\mp)\leqno(2.8)$$ with respect to the
 weak $C^\infty$-topology. Correspondingly, the notions of the
dual graphs of the stable $(J, H)_\pm$-disks with a cap $\tilde x$
may be introduced in the similar ways. Let  ${\mathcal
D}_\pm(\tilde x)$ be the sets of their dual graphs respectively.
Both are finite.

Now we are in the position to introduce the notions of positive
and negative $L^p_k$-stable disks with cap $\tilde x$. As usual it
is always assumed that $k-\frac{2}{p}\ge 1$. In principle, we may
proceed as in \cite{LiuT1}. For example, for each $D^+\in{\mathcal
D}_+(\tilde x)$ let $(f,\Sigma,\underline l)$ be a stable $(J,
H)_+$-disk as in Definition 2.2 and with the dual graph $D^+$.
Then a positive $L^p_k$-stable disk with cap $\tilde x$ and of
$D^+$ is a tuple $(\bar f,\Sigma,\underline l)$, where $\bar
f:\Sigma\setminus\{z_2,\cdots,z_{N_p+1}\}\to M$ is locally $L^p_k$
map such that (3) (4) in Definition 2.2 and the following are
satisfied:

 ${\bf (1)^\prime}$\hspace{2mm}$\bar f^P_i=\bar
f|_{P_i-\{z_i,z_{i+1}\}}$($i>1$) satisfy (i) in Definition 2.2(1)
and suitable exponential decay condition along ends $z_i$ and
$z_{i+1}$ as in [LiuT1],  but $\bar f^P_1$ only satisfies
$\lim_{s\to +\infty}\bar f^P_1(s,t)=x_2(t)$ and the exponential
decay condition along the end $z_2$.

We still define its energy $E_+(f)$ by (2.5). The equivalence
class of it can also be defined similarly. Denote by ${\mathcal
B}^{p, k}_\pm(\tilde x; H)$ the sets of equivalence classes of
those $L^p_k$-stable disks with cap $\tilde x$ and of the dual
graphs in ${\mathcal D}_\pm(\tilde x)$, and whose energy are less
than $\mp{\mathcal F}_H(\tilde x)+\max|H|+1$ (because of (2.6)).
As in \cite{LiuT1} one can equip the strong $L^p_k$-topology on
small neighborhoods ${\mathcal W}_\pm(\tilde x; H, J)$ of
 $\overline{\mathcal M}_\pm(\tilde x; H, J)$ in ${\mathcal B}^{p,k}_\pm(\tilde x, H)$ and
prove that it is equivalent to the above weak $C^\infty$-topology
on $\overline{\mathcal  M}_\pm(\tilde x; H, J)$. Once these are
well defined we may use Liu-Tian's method to construct the virtual
moduli cycles\vspace{-1mm}
$$\overline{\mathcal M}_+^{\nu^+}(\tilde x; H, J)=\!\!\!\!\!
\sum_{I\in{\mathcal N}(\tilde
x)_+}\!\frac{1}{|\Gamma_I|}\{\pi^+_I: {\mathcal
M}_+^{\nu^+_I}(\tilde x; H, J)\to {\mathcal W}_+(\tilde x; H,
J)\}\leqno(2.9)$$\vspace{-1mm}
 of dimension $\dim M-\mu(\tilde x)$ in ${\mathcal W}_+(\tilde x; H, J)$, and
$$\overline{\mathcal M}_-^{\nu_-}(\tilde x; H, J)=\!\!\!\!\!
\sum_{I\in{\mathcal N}(\tilde
x)_-}\!\frac{1}{|\Gamma_I|}\{\pi^-_I: {\mathcal
M}_-^{\nu^-_I}(\tilde x; H, J)\to {\mathcal W}_-(\tilde x; H,
J)\}\leqno(2.10)$$
 of dimension $\mu(\tilde x)$ in ${\mathcal
W}_-(\tilde x; H, J)$ (cf.\;\cite{LiuT1,LiuT2} and \cite{LiuT3}).
It should also be pointed out that the maps ${\rm EV}_\pm$ in
(2.8) can naturally extend onto the spaces ${\mathcal
W}_\pm(\tilde x; H, J)$. By composing them with the obvious
finite-to-one maps from $\overline{\mathcal M}_\pm^\nu(\tilde x;
H, J)$ to ${\mathcal W}_\pm(\tilde x; H, J)$ that forget the
parameterization we get two continuous and stratawise smooth
evaluations\vspace{-1mm}
$${\rm EV}_\pm^{\nu^\pm}(\tilde x): \overline{\mathcal M}_\pm^{\nu^\pm}(\tilde x; H, J)\to M.
\leqno(2.11)$$
Notice that  the choices of different small
$\nu^\pm$ give the cobordant virtual moduli cycles
 $\overline{\mathcal M}_\pm^{\nu^\pm}(\tilde x; H, J)$ and
 corresponding evaluations ${\rm EV}_\pm^{\nu^\pm}$.

However, notice that the boundary operator $\partial^F$ in (1.2)
depends on the choice of $\nu$. As pointed out below Remark 4.2 of
\cite{LiuT1}, for  $\partial^F$ in (1.2) being indeed a boundary
operator the choices of $\nu$ in all relative virtual  moduli
cycles $\widetilde{\mathcal  M}^\nu(\tilde x,\tilde y)$
($\mu(\tilde x)-\mu(\tilde y)\le 2$) must satisfy suitable
compatible conditions, i.e.,
$$\partial(C(\overline{\mathcal  M}^{\nu(\tilde x\tilde y)}(\tilde x,\tilde y))=
\sum_{\mu(\tilde z)=\mu(\tilde x)-1} C(\overline{\mathcal
M}^{\nu(\tilde x\tilde z)}(\tilde x,\tilde z))\times
C(\overline{\mathcal  M}^{\nu(\tilde z\tilde y)}(\tilde z,\tilde
y))\leqno(2.12)$$
 holds for all pairs $(\tilde x, \tilde y)$ with $\mu(\tilde y)-\mu(\tilde
x)=2$. Now in order to guarantee that the maps $\Phi$ and $\Psi$
constructed in (1.7) and (1.8) commute with $\partial^F$ and the
boundary operator $\partial^Q$ in (1.6), we also need carefully
choose $\nu^\pm$ in (2.9) and (2.10) so that they are compatible
with all $\nu$ chosen in the definition of $\partial^F$ in (1.2) .

\subsection{Compatible  virtual moduli cycles}\label{2.2}

 In this subsection we shall first complete Liu-Tian's arguments in
 detail,i.e., proving (2.12),   and then outline how to construct all virtual moduli
cycles  $\overline{\mathcal  M}_\pm^{\nu^\pm}(\tilde x; H, J)$
  compatible with all relative virtual moduli cycles in
(1.2).  For convenience of the later proof we need to recall
briefly the construction of the relative virtual cycle
  $C(\overline{\mathcal  M}^\nu(\tilde x, \tilde y))$  in [LiuT1].

\noindent{\it Step 1}: {\bf Local construction.}\quad For a
representative $f$ of $\langle f\rangle\in\overline{\mathcal
M}(\tilde x, \tilde y; H, J)$ one may construct a stratified
Banach orbifold chart $(\widetilde W(f),\Gamma_f, \pi_f)$ around
$\langle f\rangle$ in ${\mathcal B}(\tilde x, \tilde y)={\mathcal
B}^p_k(\tilde x, \tilde y)$, where $\Gamma_f={\rm Aut}(f)$. There
exists a natural stratified Banach bundle $\widetilde {\mathcal
L}(f)\to\widetilde W(f)$ with a stratawise smooth right
$\Gamma_f$-action  such that the usual $\bar\partial_{J,
H}$-operator gives rise to a $\Gamma_f$-equivariant stratawise
smooth section of this bundle, still denoted by $\bar\partial_{J,
H}$. Let $I$ denote  the dual graph of $f$ and $R(f)\subset
(\widetilde {\mathcal  L}(f)^I)_f$ be the cokernal ${\rm coker}(
D\bar\partial_{J,H}(f))$. Take a smooth cut-off function
$\beta_\epsilon(f)$ supported outside of the
$\epsilon$-neighborhood of double points of the domain
 $\Sigma_f$. Then each $\nu\in R_\epsilon(f):=\{\beta_\epsilon(f)\cdot\xi\,|\, \xi\in R(f)\,\}$
 may naturally determine a section of
the bundle $\widetilde {\mathcal  L}(f)\to \widetilde W(f)$,
denoted by $\tilde\nu$, such that for each $g\in\widetilde W(f)$
the support of $\tilde\nu(g)\in (\widetilde {\mathcal L}(f))_g$ is
away from the gluing region of the domain $\Sigma_g$ of $g$.

\noindent{\it Step 2}: {\bf Global construction.}\quad By the
compactness of $\overline{\mathcal M}(\tilde x, \tilde y; H, J)$
one can choose finite points $\langle f_1\rangle,\cdots, \langle
f_m\rangle$ such that the union $W\!:=\cup^m_{i=1}W(\langle
f_i\rangle)$ is an open neighborhood of $\overline{\mathcal
M}(\tilde x, \tilde y; H, J)$ and ${\mathcal
L}\!:=\cup^m_{i=1}{\mathcal L}(\langle f_i\rangle)$ is an orbifold
bundle over $W$. Let ${\mathcal N}(\tilde x,\tilde y)$ be the set
of all subsets $I=\{i_1,\cdots, i_l\}$ of $\{1,\cdots, m\}$ with
$W_I\!:=\cap_{i\in I}W(\langle f_i\rangle)\ne\emptyset$.
 Let $\pi_i: \widetilde W(f_i)\to W(\langle f_i\rangle)$
be the natural projections. For each $I\in {\mathcal N}(\tilde
x,\tilde y)$ they defined the group $\Gamma_I:=\prod_{i\in
I}\Gamma_{f_i}$ and the fiber product
$$
\widetilde W_I^{\Gamma_I}=\Bigl\{(u_i)_{i\in I}\in\prod_{i\in I}
\widetilde W(f_i)\,\Bigm|\, \pi_i(u_i)=\pi_j(u_j)\;\forall i, j\in
I\Bigr\}.
 \leqno(2.13)
 $$
 Then the projection  $\pi_I: \widetilde
W_I^{\Gamma_I}\to W_I$ has covering group $\Gamma_I$.
 Moreover,
for $J\subset I\in {\mathcal  N}(\tilde x,\tilde y)$ there is an
obvious projection $\pi^I_J: \widetilde W_I^{\Gamma_I}\to
\widetilde W_J^{\Gamma_J}$ satisfying the relation
$\pi_J\circ\pi^I_J=\pi_I$.  Repeating the same construction from
$\widetilde {\mathcal  L}(f_i)$ one obtains the bundles
$\widetilde{\mathcal  L}_I^{\Gamma_I}$ and thus a system of
bundles $(\widetilde{\mathcal  L}^\Gamma, \widetilde W^\Gamma)=
\{(\widetilde {\mathcal  L}_I^{\Gamma_I}, \widetilde
W_I^{\Gamma_I}), \pi_I,  \pi^I_J\,|\, J\subset I\in{\mathcal
N}(\tilde x,\tilde y)\}$.  Take open sets $W_i^1\subset\subset
W(\langle f_i\rangle)$, $i=1,\cdots, m$ and the pairs of open sets
$W^j_i\subset\subset U^j_i$, $i, j=1,\cdots, m$,  such that
$\cup^m_{i=1}W_i^1$ still contains $\overline{\mathcal  M}(\tilde
x,\tilde y; J, H)$ and that
$$W^1_i\subset\subset U^1_i\subset\subset W^2_i\subset\subset U^2_i\cdots
\subset\subset W^m_i\subset\subset W(\langle
f_i\rangle).\leqno(2.14)$$
In \cite{LiuT1}, for each
$I\in{\mathcal N}(\tilde x,\tilde y)$ with cardinal number
$|I|=k$, they  defined
$$V_I=\cap_{i\in I}W_i^k\setminus\cup_{|J|>k}Cl(\cap^k_{j\in J}U^k_j)\leqno(2.15)$$
and proved that  the open covering $\{V_I\,|\, I\in{\mathcal
N}(\tilde x,\tilde y)\}$ of $\overline{\mathcal M}(\tilde x,
\tilde y; H, J)$  satisfies
$$\left.\begin{array}{ll}
V_I\subset W_I\,\forall I\in {\mathcal N}(\tilde x,\tilde
y),\quad{\rm and}\\
 Cl(V_I)\cap Cl(V_J)\ne\emptyset\,{\rm only\,
if}\,I\subset J\,{\rm or}\,J\subset
I.\end{array}\right.\leqno(2.16)$$

Set $\widetilde V_I=(\pi_I)^{-1}(V_I)$ and $\widetilde
E_I=(\pi_I)^{-1}({\mathcal  L}|_{V_I})$, one gets a system of
bundles
$$(\widetilde E^\Gamma, \widetilde V^\Gamma)=
\Bigl\{(\widetilde E_I^{\Gamma_I}, \widetilde V_I^{\Gamma_I}),
\pi_I, \Gamma_I, \pi^I_J\,\Bigm|\, J\subset I\in{\mathcal
N}(\tilde x,\tilde y)\Bigr\}. \leqno(2.17)$$
Taking
$\Gamma_{f_i}$-invariant stratawise smooth cut-off function
$\gamma(f_i)$ on $\widetilde W(f_i)$ such that\vspace{-1mm}
$$\gamma(f_i)=1\;\;{\rm in}\;\;\pi_i^{-1}( W^m_i),\leqno(2.18)$$
then  each $\nu_i\in R_\epsilon(f_i)$ determines a smooth global
section $\bar\nu_i$ of  $(\widetilde E^{\Gamma}, \widetilde
V^{\Gamma})$. Set
$$R^\epsilon_\delta(\{f_i\})=\{\nu\in \oplus^m_{i=1}R_\epsilon(f_i)\,|\,|\nu|<\delta\}\leqno(2.19)$$
for a small $\delta>0$. The bundle system  $(\widetilde
E^\Gamma\times R^\epsilon_\delta(\{f_i\}), \widetilde
V^\Gamma\times R^\epsilon_\delta(\{f_i\}))$ has a well-defined
global section\vspace{-1mm}
$$
\bar\partial_{J,H}+e:\, (u_I, \nu)\mapsto \bar\partial_{J,H}u_I+
\sum^m_{i=1}(\bar\nu_i)_I(u_I).
\leqno(2.20)$$
 for any $(u_I, \nu)=(u_I,
(\nu_1,\cdots,\nu_m))\in\widetilde V_I\times R^\epsilon_\delta(\{f_i\})$.
 Moreover, each  $\nu\in R_\delta^\epsilon(\{f_i\})$ yields a smooth section
$\bar\partial_{J,H}^\nu=\{(\bar\partial_{J,H}^\nu)_I\,|\,
I\in{\mathcal  N}(\tilde x, \tilde y)\}$ of  $(\widetilde
E^\Gamma,\widetilde V^\Gamma)$,
$$
(\bar\partial_{J,H}^\nu)_I(u_I)= \bar\partial_{J,H}u_I+
\sum^m_{i=1}(\bar\nu_i)_I(u_I)\quad \forall u_I\in\widetilde
V_I.
\leqno(2.21)$$

\noindent{\bf Theorem 2.6}(\cite{LiuT1}).\hspace{2mm}{\it The
section in (2.20) is smooth and transversal to the zero section.
Therefore when $\delta>0$ is small enough, for a generic choice of
$\nu\in R_\delta^\epsilon(\{f_i\})$ the section
$\bar\partial_{J,H}^\nu$ is transversal to the zero section. Thus
the family of perturbed moduli spaces $\widetilde{\mathcal
M}^\nu=\widetilde{\mathcal  M}^\nu(\tilde x,\tilde y)=
\{\widetilde{\mathcal  M}_I^\nu=
(\bar\partial_{J,H}^\nu)^{-1}_I(0)\,|\, I\in{\mathcal  N}(\tilde
x, \tilde y)\}$ is compatible in the sense that
$\pi^I_J(\widetilde{\mathcal  M}_I^\nu)=\widetilde{\mathcal
M}_J^\nu\cap ({\rm Im}\pi^I_J)$ for all $J\subset I\in{\mathcal
N}(\tilde x,\tilde y)$.}\vspace{2mm}

Let $W^m_i$ be given by (2.14). For sufficiently small $\nu$ we
can also require that
$$\widetilde{\mathcal  M}^\nu_I\subset\cap_{i\in I}\pi_i^{-1}(W^m_i).\leqno(2.22)$$
 Let ${\widetilde {\mathcal  M}}_I^{\nu, D_T}$ (resp.
 \! ${\widetilde {\mathcal  M}}_I^{\nu, D_B}$)
  be the top strata (resp. the strata of
``broken'' connecting orbits) of ${\widetilde {\mathcal
M}}_I^{\nu}$.
 Then ${\widetilde {\mathcal
M}}_I^{\nu, c}:=
 {\widetilde {\mathcal  M}}_I^{\nu, D_T}\cup{\widetilde {\mathcal  M}}_I^{\nu, D_B}$
is a smooth manifold of dimension $\mu(\tilde x)-\mu(\tilde y)-2$
and with boundary $\partial {\widetilde {\mathcal  M}}_I^{\nu,
c}={\widetilde {\mathcal  M}}_I^{\nu, D_B}$. Formally writing
$$
C(\overline{\mathcal  M}^\nu(\tilde x,\tilde y))
=\sum_{I\in{\mathcal  N}(\tilde x,\tilde y)}
\frac{1}{|\Gamma_I|}{\widetilde {\mathcal  M}}_I^{\nu,
c},\quad\quad\overline{\mathcal  M}^\nu(\tilde x,\tilde y)=
\sum_{I\in{\mathcal N}(\tilde x,\tilde
y)}\frac{1}{|\Gamma_I|}{\widetilde {\mathcal M}}_I^{\nu},
$$
 the former was called the relative moduli cycle in \cite{LiuT1}.
 Later we call the compatible family $\widetilde{\mathcal M}^\nu(\tilde
x,\tilde y)$ in Theorem 2.6
 a {\it derived family} for $C(\overline{\mathcal  M}^\nu(\tilde x,\tilde
 y))$ and $\overline{\mathcal  M}^\nu(\tilde x,\tilde
y)$. It is clear that $C(\overline{\mathcal  M}^\nu(\tilde
x,\tilde y))= \overline{\mathcal M}^\nu(\tilde x,\tilde y)$ in the
case $\mu(\tilde x)-\mu(\tilde y)\le 2$. For the sake of clearness
$C(\overline{\mathcal  M}^\nu(\tilde x,\tilde y))$ will  be
denoted by $C\bigl(\overline{\mathcal M}^{\nu(\tilde x\tilde
y)}(\tilde x,\tilde y)\bigr)$ below.
 If $\mu(\tilde
x)-\mu(\tilde y)=2$, then for every $\tilde z\in\tilde{\mathcal
P}(H)$ with $\mu(\tilde x)-\mu(\tilde z)=1$ one has also the
associated relative moduli cycles $C(\overline{\mathcal
M}^{\nu(\tilde x\tilde z)}(\tilde x,\tilde z))$ and
$C(\overline{\mathcal M}^{\nu(\tilde z\tilde y)}(\tilde z,\tilde
y))$.  The relative virtual moduli cycles satisfying (2.12) are
called {\it compatible}.\vspace{2mm}

 \noindent{\it Proof of (2.12).}\quad To construct such  relative virtual cycles,
 note that there are only finitely many $\tilde z\in\tilde{\mathcal
P}(H)$, saying $\tilde z_1,\cdots, \tilde z_r$, such that
$\mu(\tilde z_q)=\mu(\tilde x)-1$ and
 $\overline{\mathcal M}(\tilde x,\tilde z_q; H, J)\ne\emptyset$,
  $\overline{\mathcal M}(\tilde z_q,\tilde y; H, J)\ne\emptyset$ for $q=1,\cdots, r$.
Let $\langle f^{(1q)}_s\rangle\in\overline{\mathcal  M}(\tilde x,
\tilde z_q; H, J)$ and $\langle
f^{(2q)}_t\rangle\in\overline{\mathcal M}(\tilde z_q, \tilde y; H,
J)$, $s=1,\cdots, m_q$, $t=1,\cdots, n_q$, be finite points from
which one may construct  the relative virtual moduli cycles
$$
C(\overline{\mathcal  M}^{\nu(\tilde x\tilde z_q)}(\tilde x,
\tilde z_q))= \sum_{I\in{\mathcal  N}(\tilde x,\tilde z_q)}
\frac{1}{|\Gamma_I^{(1q)}|}{\widetilde {\mathcal M}}_I^{\nu(\tilde
x\tilde z_q), c},
 \leqno(2.23)$$\vspace{-3mm}
$$
C(\overline{\mathcal  M}^{\nu(\tilde z_q\tilde y)}(\tilde z_q,
\tilde y))= \sum_{I\in{\mathcal  N}(\tilde z_q,\tilde y)}
\frac{1}{|\Gamma_I^{(2q)}|}{\widetilde {\mathcal M}}_I^{\nu(\tilde
z_q\tilde y), c}.
\leqno(2.24)$$
 For the future convenience we assume that for $q=1,\cdots,r$,
$$\begin{array}{ll}
 W^1(\langle
f^{(1q)}_s\rangle)\subset\subset U^1(\langle
f^{(1q)}_s\rangle)\cdots \subset\subset U^{m_q-1}(\langle
f^{(1q)}_s\rangle ) \\
\subset\subset W^{m_q}(\langle f^{(1q)}_s\rangle)\subset\subset
W(\langle f^{(1q)}_s\rangle),\quad s=1,\cdots, m_q,
\end{array}\leqno(2.25)$$
 and
$$
\begin{array}{ll}
  W^1(\langle
f^{(2q)}_t\rangle)\subset\subset U^1(\langle f^{(2q)}_t\rangle)
\cdots\subset\subset U^{n_q-1}(\langle f^{(2q)}_t\rangle)\\
\subset\subset W^{n_q}(\langle f^{(2q)}_t\rangle) \subset\subset
W(\langle f^{(2q)}_t\rangle),\quad t=1,\cdots, n_q,\quad\quad\quad
\end{array}\leqno(2.26)$$
 are respectively the open sets as defined in (2.14) that are used to
construct the relative virtual moduli cycles in (2.23) and (2.24)
above. By (2.15)-(2.17) we get the corresponding bundle systems
$(\widetilde E^{\Gamma^{(1q)}}, \widetilde V^{\Gamma^{(1q)}})$ and
$(\widetilde E^{\Gamma^{(2q)}}, \widetilde V^{\Gamma^{(2q)}})$,
$q=1,\cdots,r$. As in Theorem 2.6 let
$$\left.\begin{array}{ll} \nu(\tilde x\tilde
z_q)=\oplus^{m_q}_{s=1}\nu(\tilde x\tilde z_q)_s
\in R^\epsilon_\delta(\{f^{(1q)}_s\})\quad{\rm and}\quad\quad \\
\nu(\tilde z_q\tilde y)=\oplus^{n_q}_{t=1}\nu(\tilde z_q\tilde
y)_t \in R^\epsilon_\delta(\{f^{(2i)}_t\})
\end{array}\right.\leqno(2.27)$$
be such that the smooth section $\bar\partial_{J,H}^{\nu(\tilde
x\tilde z_q)}$ of  the bundle system $(\widetilde
E^{\Gamma^{(1q)}}, \widetilde V^{\Gamma^{(1q)}})$
 and that $\bar\partial_{J,H}^{\nu(\tilde z_q\tilde y)}$ of
$(\widetilde E^{\Gamma^{(2q)}}, \widetilde V^{\Gamma^{(2q)}})$ are
transversal to the zero section respectively. Here $\nu(\tilde
x\tilde z_q)_s\in R_\epsilon(f^{(1q)}_s)$ and $\nu(\tilde
z_q\tilde y)_t\in R_\epsilon(f^{(2q)}_t)$, $s=1,\cdots, m_q$ and
$t=1,\cdots, n_q$. In fact, as in (2.21) we can also require
$\nu(\tilde x\tilde z_q)$ and $\nu(\tilde z_q\tilde y)$ so small
that
$$\left.\begin{array}{ll}
\widetilde{\mathcal  M}^{\nu(\tilde x\tilde z_q)}_I
\subset\cap_{s\in I}(\pi^{(1q)}_s)^{-1}(W^{m_q}(\langle
f^{(1q)}_s\rangle))\quad \forall I\in{\mathcal  N}(\tilde x,\tilde
z_q),\\
\widetilde{\mathcal  M}^{\nu(\tilde z_q\tilde y)}_I
\subset\cap_{t\in I}(\pi^{(2q)}_t)^{-1}(W^{n_q}(\langle
f^{(2q)}_t\rangle))\quad \forall I\in{\mathcal  N}(\tilde
z_q,\tilde y).\end{array}\right.$$
 For $q=1,\cdots, r$ we set
$f^{(q)}_{st}:=f^{(1q)}_s\sharp_{\tilde z_q}f^{(2q)}_t$,
$s=1,\cdots, m_q,\, t=1,\cdots, n_q$.
 Clearly, these $\langle f^{(q)}_{st}\rangle$ belong to $\overline{\mathcal  M}(\tilde
x,\tilde y; H, J)$. Moreover, the automorphism group
$\Gamma^{(q)}_{st}$ of $f^{(q)}_{st}$ may be identified with the
product $\Gamma^{(1q)}_s\times\Gamma^{(2q)}_t$ of the automorphism
group $\Gamma^{(1q)}_s$ of $f^{(1q)}_s$ and that $\Gamma^{(2q)}_t$
of $f^{(2q)}_t$. By the construction of the local uniformizer in
\S2 of \cite{LiuT1} we easily construct a uniformizer $\widetilde
W(f^{(q)}_{st})$ such that
$$\widetilde W(f^{(1q)}_s)\sharp\widetilde W(f^{(2q)}_t)\subset
\widetilde W(f^{(q)}_{st})\leqno(2.28)$$
 and  that the restriction of $\Gamma^{(q)}_{st}$-action over $\widetilde
W(f^{(q)}_{st})$ to $\widetilde W(f^{(1q)}_s)\sharp\widetilde
W(f^{(2q)}_t)$ is exactly that of
$\Gamma^{(1q)}_s\times\Gamma^{(2q)}_t$ over $\widetilde
W(f^{(1q)}_s)\sharp\widetilde W(f^{(2q)}_t)$ in the obvious way,
where $\widetilde W(f^{(1q)}_s)\sharp\widetilde W(f^{(2q)}_t)$
 denotes the set of all join functions at $\tilde z_q$ of functions in
$\widetilde W(f^{(1q)}_s)$ and $\widetilde W(f^{(2q)}_t)$.
 (One may increase $m_q$, $n_q$ and shrink $\widetilde
W(f^{(1q)}_s)$ and $\widetilde W(f^{(2q)}_t)$ if necessary.)
 Since
each $\overline{\mathcal  M}(\tilde x,\tilde z_q; J, H)\sharp
\overline{\mathcal  M}(\tilde z_q,\tilde y; J, H)$ is compact and
$\cup^r_{q=1}\overline{\mathcal  M}(\tilde x,\tilde z_q; J,
H)\sharp \overline{\mathcal  M}(\tilde z_q,\tilde y; J, H)$
 are disjoint unions we can require that the above uniformizers
$\widetilde W(f^{(q)}_{st})$ satisfy
$$\widetilde W(f^{(q)}_{st})\cap \widetilde W(f^{(q^\prime)}_{s^\prime t^\prime})=
\emptyset\quad\forall q\ne q^\prime.\leqno(2.29)$$
 By the definition of the index set ${\mathcal  N}(\tilde x,\tilde y)$
above (2.13) it easily follows from (2.29)  that the corresponding
index set ${\mathcal  N}^0(\tilde x,\tilde y)$  with the
collection $\{\widetilde W(f^{(q)}_{st})\;|\, 1\le s\le m_q,\;1\le
t\le n_q,\; 1\le q\le r\}$ must have the following form
$${\mathcal  N}^0(\tilde x,\tilde y)=\cup^r_{q=1}
{\mathcal  N}(\tilde x,\tilde z_q)\times {\mathcal  N}(\tilde
z_q,\tilde y).\leqno(2.30)$$\vspace{-2mm}

Notice that every $W(\langle f^{(q)}_{st}\rangle)$ determines an
open neighborhood $W(\langle f^{(q)}_{st}\rangle)_1$ of $\langle
f^{(1q)}_s\rangle$ in ${\mathcal  B}(\tilde x,\tilde z_q)$ by
$$\Big\{\,\langle g_1\rangle\in{\mathcal  B}(\tilde x,\tilde z_q)\,\Bigm|\,
\exists \langle g_2\rangle\in{\mathcal  B}(\tilde z_q, \tilde y)\;
s.t.\; \langle g_1\sharp g_2\rangle\in {\mathcal  B}(\tilde
x,\tilde y)\cap W(f^{(q)}_{st})\,\Bigr\}$$
  and that $W(\langle f^{(q)}_{st}\rangle)_2$ of $\langle f^{(2q)}_t\rangle$ in
 ${\mathcal B}(\tilde z_q, \tilde y)$ by
$$\Bigl\{\,\langle g_2\rangle\in
{\mathcal  B}(\tilde z_q,\tilde x)\,\Bigm|\; \exists \langle
g_1\rangle\in{\mathcal  B}(\tilde x, \tilde z_q)\; s.t.\; \langle
g_1\sharp g_2\rangle\in {\mathcal  B}(\tilde x,\tilde y)\cap
W(f^{(q)}_{st})\,\Bigr\}.$$
 Thus (2.28) implies that
$W(\langle f^{(1q)}_s\rangle)\subset W(\langle
f^{(q)}_{st}\rangle)_1$ and $W(\langle f^{(2q)}_t\rangle)\subset
W(\langle f^{(q)}_{st}\rangle)_2$ for all $q, s, t$. By this and
(2.25)(2.26) we can, as in (2.14), choose pairs of open sets
$W^j(\langle f^{(q)}_{st}\rangle)\subset\subset U^j(\langle
f^{(q)}_{st}\rangle)$,
 $j=1,\cdots, m=\sum^r_{q=1}m_q n_q$,
 such that
 $$\left.\begin{array}{ll}W^1(\langle
f^{(q)}_{st}\rangle)\subset\subset U^1(\langle
f^{(q)}_{st}\rangle) \subset\subset W^2(\langle
f^{(q)}_{st}\rangle)\cdots  \\
\subset\subset W^m(\langle f^{(q)}_{st}\rangle)\subset\subset
W(\langle f^{(q)}_{st}\rangle)
\end{array}\right.\leqno(2.31)$$
 and  that
$W^{m_q}(\langle f^{(1q)}_s\rangle)\subset W^1(\langle
f^{(q)}_{st}\rangle)_1$ and $W^{n_q}(\langle
f^{(2q)}_t\rangle)\subset W^1(\langle f^{(q)}_{st}\rangle)_2$ for
all  $q, s, t$.
 As in (2.15) and (2.17) we use (2.30) and (2.31) to construct a bundle system
$$(\widetilde E, \widetilde V)=\{(\widetilde E_{I_1\times I_2},
 \widetilde V_{I_1\times I_2})\;|\; I_1\times I_2\in{\mathcal  N}^0(\tilde x,\tilde y)\}.
$$
As in (2.18) we take $\Gamma_{st}^{(q)}$-invariant smooth cut-off
function $\gamma(f^{(q)}_{st})$ on $\widetilde W(\langle
f^{(q)}_{st}\rangle)$  such that $\gamma(f^{(q)}_{st})=1$ in
$(\pi^q_{st})^{-1}(W^m(\langle f^{(q)}_{st}\rangle))$. Note that
the supports of $\nu(\tilde x\tilde z_q)_s$ and $\nu(\tilde
z_q\tilde y)_t$ in (2.27) are away from double points on their
domains. So the support of $\nu(\tilde x\tilde y)^q_{st}:=
\nu(\tilde x\tilde z_q)_s\sharp\nu(\tilde z_q\tilde y)_t$ is away
from all double points of the domain of $f^{(q)}_{st}$ and thus
$\nu(\tilde x\tilde y)^q_{st}$ sits in the fibre of
$\widetilde{\mathcal L}(f^{(q)}_{st})$ at $f^{(q)}_{st}$,
$L^p_{k-1}(\wedge^{0,1}((f^{(q)}_{st})^\ast TM))$. As before we
can use $\gamma(f^{(q)}_{st})$ and $\nu(\tilde x\tilde y)^q_{st}$
to get a global section
$$\overline{\nu(\tilde x\tilde y)^q_{st}}=\{
(\overline{\nu(\tilde x\tilde y)^q_{st}})_{I_1\times I_2}\;|\;
I_1\times I_2\in{\mathcal  N}^0(\tilde x,\tilde y)\}$$ of the
bundle system $(\widetilde E, \widetilde V)$. Let us set
$$\nu(\tilde x\tilde y)^0:=\sum^r_{q=1}\sum^{m_q}_{s=1}\sum^{n_q}_{t=1}
\nu(\tilde x\tilde y)^q_{st}\quad{\rm
and}\quad\overline{\nu(\tilde x\tilde
y)^0}:=\sum^r_{q=1}\sum^{m_q}_{s=1}\sum^{n_q}_{t=1}
\overline{\nu(\tilde x\tilde y)^q_{st}}.$$
 As in (2.21) we get a global section of the bundle system
 $(\widetilde E, \widetilde V)$,
$$\bar\partial_{J,H}^{\nu(\tilde x\tilde y)^0}=
\biggl\{\Bigl(\bar\partial_{J,H}^{\nu(\tilde x\tilde
y)^0}\Bigr)_{I_1\times I_2}=
\bar\partial_{J,H}+(\overline{\nu(\tilde x\tilde y)^0})_{I_1\times
I_2}\; \Bigm|\; I_1\times I_2\in{\mathcal N}^0(\tilde x,\tilde
y)\biggr\}.$$
 Then it easily follows from (2.29) and the choice of this section  that
$$\widetilde{\mathcal  M}^{\nu(\tilde x\tilde z_q)}_{I_1}\sharp
\widetilde{\mathcal  M}^{\nu(\tilde z_q\tilde y)}_{I_2}\subset
\widetilde{\mathcal  M}^{\nu(\tilde x\tilde y)^0}_{I_1\times
I_2}:= \biggl(\Bigl(\bar\partial_{J,H}^{\nu(\tilde x\tilde
y)^0}\Bigr)_{I_1\times I_2}\biggr)^{-1}(0) $$
for $I_1\times
I_2\in {\mathcal  N}(\tilde x,\tilde z_q)\times {\mathcal
N}(\tilde z_q,\tilde y)$, and that each section
$(\bar\partial_{J,H}^{\nu(\tilde x\tilde y)^0})_{I_1\times I_2}$
is transversal to the zero section at all points of
$$\widetilde{\mathcal  M}^{\nu(\tilde x\tilde z_q)}_{I_1}\sharp
\widetilde{\mathcal  M}^{\nu(\tilde z_q\tilde y)}_{I_2}\subset
\widetilde W^{\nu(\tilde x\tilde z_q)}_{I_1}\sharp \widetilde
W^{\nu(\tilde z_q\tilde y)}_{I_2}.$$
 The last set consists of all
points $(u_s\sharp v_t)_{(s,t)\in I_1\times I_2}$ in
$$
\prod_{(s,t)\in I_1\times I_2} \widetilde
W(f^{(1q)}_s)\sharp\widetilde W(f^{(2q)}_t)
$$ such that
$\pi^{(q)}_{st}(u_s\sharp v_t)=\pi^{(q)}_{s^\prime t^\prime}(
u_{s^\prime}\sharp v_{t^\prime})$
 for any $(s,t)$ and $(s^\prime,
t^\prime)$ in $I_1\times I_2$. By the stability of a surjective
map under small perturbation we can show that there exists a
$\Gamma^{(q)}_{I_1\times I_2}$-invariant open neighborhoods of
$\widetilde{\mathcal  M}^{\nu(\tilde x\tilde z_q)}_{I_1}\sharp
\widetilde{\mathcal  M}^{\nu(\tilde z_q\tilde y)}_{I_2}$ in
$\widetilde V_{I_1\times I_2}$,
$$\widetilde V^0_{I_1\times I_2}\subset\subset\widetilde V^1_{I_1\times I_2}
\leqno(2.32)$$
such that the section
$(\bar\partial_{J,H}^{\nu(\tilde x\tilde y)^0})_{I_1\times I_2}$
is transversal to the zero section at points of
$$\widetilde V^1_{I_1\times I_2}\cap\widetilde{\mathcal  M}^{\nu(\tilde x\tilde y)^0}_{
I_1\times I_2}.\leqno(2.33)$$
So the space in (2.33) is a cornered
and stratified Banach variety that has the dimension given by the
Index Theorem on all of its strata. Moreover, the collection
$$\Bigl\{ \widetilde{\mathcal  M}^{\nu(\tilde x\tilde y)^0}_{I_1\times I_2}\;|\;
I_1\times I_2\in {\mathcal  N}(\tilde x,\tilde z_q)\times
{\mathcal N}(\tilde z_q,\tilde y)\Bigr\}
$$
also satisfy the
compatibility as in Theorem 2.6. In particular the open
neighborhoods in (2.32) may be chosen to guarantee that the
collection
$$\Bigl\{\widetilde V^1_{I_1\times I_2}\cap\widetilde{\mathcal
 M}^{\nu(\tilde x\tilde y)^0}_{
I_1\times I_2}\;|\; I_1\times I_2\in {\mathcal  N}(\tilde x,\tilde
z_q)\times {\mathcal  N}(\tilde z_q,\tilde y)\Bigr\}
$$
 satisfies
the compatibility. Note that the projection image of
 $\widetilde V^1_{I_1\times I_2}\cap\widetilde{\mathcal  M}^{\nu(\tilde x\tilde y)^0}_{
I_1\times I_2}$ in ${\mathcal  B}(\tilde x,\tilde y)$ is not
compact in general (unlike in Theorem 2.6). We denote by
$$\Bigl(\widetilde V^1_{I_1\times I_2}\cap
\widetilde{\mathcal  M}^{\nu(\tilde x\tilde y)^0}_{I_1\times
I_2}\Bigr)^c,$$ the union of the top and $1$-codimensional strata
of $\widetilde V^1_{I_1\times I_2}\cap \widetilde{\mathcal
M}^{\nu(\tilde x\tilde y)^0}_{I_1\times I_2}$, then it is a smooth
manifold with boundary of dimension $\mu(\tilde x)-\mu(\tilde
y)-1$ that is contained in the smooth locus of $\widetilde
V^1_{I_1\times I_2}\cap\widetilde{\mathcal  M}^{\nu(\tilde x\tilde
y)^0}_{ I_1\times I_2}$ and that has the boundary
$$\partial\Bigl(\widetilde V^1_{I_1\times I_2}\cap
\widetilde{\mathcal  M}^{\nu(\tilde x\tilde y)^0}_{I_1\times I_2}
\Bigr)^c=\widetilde{\mathcal  M}^{\nu(\tilde x\tilde z_q),
D_T}_{I_1}\sharp \widetilde{\mathcal  M}^{\nu(\tilde z_q\tilde y),
D_T}_{I_2}$$
 if $I_1\times I_2\in{\mathcal  N}(\tilde x,\tilde z_q)\times
{\mathcal  N}(\tilde z_q,\tilde y)$. In the following we shall
extend
$$\Bigl\{\widetilde V^1_{I_1\times I_2}\cap
\widetilde{\mathcal  M}^{\nu(\tilde x\tilde y)^0}_{I_1\times
I_2}\;|\; I_1\times I_2\in{\mathcal  N}(\tilde x,\tilde
y)^0\Bigr\}
$$
into a virtual moduli cycle for $\overline{\mathcal
M}(\tilde x,\tilde y; J, H)$ under condition that
$$\Bigl\{\widetilde V^0_{I_1\times I_2}\cap
\widetilde{\mathcal  M}^{\nu(\tilde x\tilde y)^0}_{I_1\times
I_2}\;|\; I_1\times I_2\in{\mathcal  N}(\tilde x,\tilde
y)^0\Bigr\}$$
 is not changed. Since both
 $\overline{\mathcal  M}(\tilde x,\tilde y; J, H)$ and
$$\cup^r_{q=1} \overline{\mathcal  M}(\tilde x,\tilde z_q; J,
H)\sharp \overline{\mathcal  M}(\tilde z_q,\tilde y; J, H)$$
 are compact  we can take points $\langle h_1\rangle,\cdots, \langle
h_m\rangle$ in
$$\overline{\mathcal  M}(\tilde x,\tilde y; J, H)\setminus\cup^r_{q=1}
\overline{\mathcal  M}(\tilde x,\tilde z_q; J, H)\sharp
\overline{\mathcal  M}(\tilde z_q,\tilde y; J, H)$$ and their open
neighborhoods $W(\langle h_j\rangle)=\pi_j^{\tilde x\tilde y}(
\widetilde W(h_j))$, $j=1,\cdots, n$, in ${\mathcal  B}(\tilde
x,\tilde y)$ such that all $W(\langle h_j\rangle)$, $W(\langle
f^{(q)}_{st}\rangle)$ form an open covering of $\overline{\mathcal
M}(\tilde x,\tilde y; J, H)$ satisfying the conditions for the
construction of the virtual moduli cycle above. Let us choose
$\Gamma_{h_j}$-invariant smooth cut-off functions $\gamma(h_j)$ on
$\widetilde W(h_j)$ such that
$$\pi_j({\rm suppt}(\gamma(h_j)))\cap V^0_{I_1\times I_2}=\emptyset\;\;\forall 1\le j\le n,
I_1\times I_2\in{\mathcal  N}(\tilde x,\tilde y)^0$$
 where $V^0_{I_1\times I_2}$ is the projection image of $\widetilde
V^0_{I_1\times I_2}$ in ${\mathcal  B}(\tilde x,\tilde y)$.

Let $(\widetilde E(\tilde x\tilde y), \widetilde V(\tilde x\tilde
y))$ and ${\mathcal  N}(\tilde x,\tilde y)$ be the corresponding
bundle system and index set to this covering. Assume that
$$\overline{\overline{\nu(\tilde x\tilde y)^q_{st}}}=\Bigl\{
(\overline{\overline{\nu(\tilde x\tilde y)^q_{st}}})_I\,|\,
I\in{\mathcal  N}(\tilde x,\tilde y)\Bigr\}$$
 is a global section of this bundle system obtained by the cut-off function
$\gamma(f^q_{st})$ and $\nu(\tilde x\tilde y)^q_{st}$
 as in Step 2. Of course, for each $\nu_j\in R_\epsilon(h_j)$
 we still denote by
$\bar\nu_j=\{(\bar\nu_j)_I\,|\, I\in{\mathcal  N}(\tilde x,\tilde
y)\} $ the global section of $(\widetilde E(\tilde x\tilde y),
\widetilde V(\tilde x\tilde y))$ obtained from $\gamma(h_j)$ and
$\nu_j$ as below (2.18).

For $\delta>0$ we assume that $Z^\epsilon_\delta(\{h_j\})$ is a
$\delta$-neighborhood of zero of $\oplus^n_{j=1}R_\epsilon(h_j)$.
As before we have a bundle system
$(\widetilde E(\tilde x\tilde y)\times Z^\epsilon_\delta(\{h_j\}),
\widetilde V(\tilde x\tilde y)\times Z^\epsilon_\delta(\{h_j\}))$
and a well-defined global section of it
$$\bar\partial_{J,H}+\hat e:\;(u_I, \nu)\mapsto
\bar\partial_{J,H}u_I+\sum^n_{j=1}(\bar\nu_j)_I(u_I)+
\sum^r_{q=1}\sum^{m_q}_{s=1}\sum^{n_q}_{t=1}
 \Bigl(\overline{\overline{\nu(\tilde x\tilde y)^q_{st}}}\Bigr)_I(u_I)
$$ for $(u_I, \nu)=(u_I, (\nu_1,\cdots,\nu_n))\in \widetilde
V(\tilde x\tilde y)\times Z^\epsilon_\delta(\{h_j\})$. Notice that
our choices above imply this section to be transversal to the zero
section for $\delta>0$ small enough. It follows that for a generic
choice of $\nu\in Z^\epsilon_\delta(\{h_j\})$ the section of
$(\widetilde E(\tilde x\tilde y), \widetilde V(\tilde x\tilde y))$
given by\vspace{-2mm}
$$u_I\mapsto \bar\partial_{J,H}u_I+\sum^n_{j=1}(\bar\nu_j)_I(u_I)+
\sum^r_{q=1}\sum^{m_q}_{s=1}\sum^{n_q}_{t=1}
 (\overline{\overline{\nu(\tilde x\tilde y)^q_{st}}})_I(u_I)
\leqno(2.34)$$
for $u_I\in\widetilde V(\tilde x\tilde y)_I$, is
transversal to the zero section.

Setting $\nu(\tilde x\tilde
y):=(\oplus^r_{q=1}\oplus^{m_q}_{s=1}\oplus^{n_q}_{t=1}
 \nu(\tilde x\tilde y)^q_{st})\oplus(\oplus^n_{j=1}\nu_j)$
then it belongs to
$(\oplus^r_{q=1}\oplus^{m_q}_{s=1}\oplus^{n_q}_{t=1}
 R_\epsilon(f^{(1q)}_s)\sharp R_\epsilon(f^{(2q)}_t))
\oplus(\oplus^n_{j=1}R_\epsilon(h_j))$, and the section
$$\bar\partial_{J,H}^{\nu(\tilde x\tilde y)}=
\Bigl\{\,\Bigl(\bar\partial_{J,H}^{\nu(\tilde x\tilde
y)}\Bigr)_I\,|\, I\in{\mathcal  N}(\tilde x,\tilde y)\Bigr\}$$
 of $(\widetilde E(\tilde x\tilde y), \widetilde V(\tilde x\tilde y))$
is transversal to the zero section. Here
$(\bar\partial_{J,H}^{\nu(\tilde x\tilde y)})_I$ is given by
(2.34). Denote by $C(\overline{\mathcal  M}^{\nu(\tilde x\tilde
y)}(\tilde x,\tilde y))$ the virtual moduli cycle constructed from
this section.  It is easy to see that its boundary is given by
$$\sum^r_{q=1}\sum_{I_1\in{\mathcal  N}(\tilde x,\tilde z_q)}
\sum_{I_2\in{\mathcal  N}(\tilde z_q,\tilde y)}
\frac{1}{|\Gamma_{I_1}^{(1q)}|\cdot |\Gamma_{I_2}^{(2q)}|}
\widetilde{\mathcal  M}_{I_1}^{\nu(\tilde x\tilde z_q), D_T}\sharp
\widetilde{\mathcal  M}_{I_2}^{\nu(\tilde z_q\tilde y), D_T},$$
which may be identified with
$$\sum^r_{q=1}\sum_{I_1\in{\mathcal  N}(\tilde x,\tilde z_q)}
\sum_{I_2\in{\mathcal  N}(\tilde z_q,\tilde y)}
\frac{1}{|\Gamma_{I_1}^{(1q)}|\cdot |\Gamma_{I_2}^{(2q)}|}
\widetilde{\mathcal  M}_{I_1}^{\nu(\tilde x\tilde z_q), D_T}\times
\widetilde{\mathcal  M}_{I_2}^{\nu(\tilde z_q\tilde y), D_T}.$$
But the latter is just $\sum^r_{q=1}C(\overline{\mathcal
M}^{\nu(\tilde x\tilde z_q)}(\tilde x,\tilde z_q)) \times
C(\overline{\mathcal  M}^{\nu(\tilde z_q\tilde y)}(\tilde
z_q,\tilde y))$. (2.12) is proved.\hfill$\Box$

Now we first construct all relative virtual moduli cycles
$C(\overline{\mathcal  M}^{\nu(\tilde x\tilde y)}(\tilde x,\tilde
y))$ for all $\tilde x,\tilde y\in\tilde{\mathcal  P}(H)$ with
$\mu(\tilde x)-\mu(\tilde y)=1$,  such that
$$\sharp(C(\overline{\mathcal  M}^\nu(\tilde x, \tilde y)))
 =\sharp(C(\overline{\mathcal  M}^\nu(\tilde x\sharp(-A), \tilde y\sharp(-A))))\leqno(2.35)$$
 for any $A\in\Gamma$.
Then we
 follow the above methods
to construct all relative virtual moduli cycles
$C(\overline{\mathcal M}^{\nu(\tilde x\tilde y)}(\tilde x,\tilde
y))$ for all $\tilde x,\tilde y\in\tilde{\mathcal  P}(H)$ with
$\mu(\tilde x)-\mu(\tilde y)=2$. Such the family of the relative
virtual moduli cycles is compatible and thus satisfies the
requirements in the definition of $\partial^F$.\vspace{2mm}

\noindent{\bf Remark 2.7.}\hspace{2mm} As in \cite{LiuT2} and
\cite{LiuT3} we need the virtual moduli cycles of dimension more
than one in this paper, and can use the notion of local components
in \cite{LiuT3}
 to construct a desingularization of the bundle system
$(\widetilde{\mathcal L}^\Gamma, \widetilde W^\Gamma)$, a new
bundle system $(\widehat{\mathcal  L}^\Gamma, \widehat
W^\Gamma)=\{(\widehat{\mathcal  L}_I^{\Gamma_I}, \widehat
W_I^{\Gamma_I})\,|\, I\in{\mathcal  N}(\tilde x,\tilde y)\}$ such
that each $\widehat{\mathcal  L}_I^{\Gamma_I}$ (resp. $\widehat
W_I^{\Gamma_I}$) is a stratified Banach manifold (resp. bundle).
Then replacing $(\widetilde{\mathcal  L}^\Gamma, \widetilde
W^\Gamma)$ everywhere by $(\widehat{\mathcal L}^\Gamma, \widehat
W^\Gamma)$ in the previous construction of the virtual moduli
cycles one can get a virtual moduli cycle
$$
\sum_{I\in{\mathcal  N}(\tilde x,\tilde
y)}\frac{1}{|\Gamma_I|}\Bigl\{ \pi_I: \widehat{\mathcal
M}_I^\nu\to W\Bigr\},\;\;{\rm still}\;{\rm denoted}\;{\rm
by}\;\overline{\mathcal M}^\nu(\tilde x, \tilde y),$$ such that
each $\widehat{\mathcal M}_I^\nu$ is a cornered smooth manifold.
We can use the same method to  make suitable modifications for the
above arguments and in the case $\mu(\tilde x)-\mu(\tilde y)=2$
obtain the conclusion corresponding with (2.12), i.e.,
$$\left.\begin{array}{ll}\partial\overline{\mathcal  M}^{\nu(\tilde x\tilde y)}(\tilde x,\tilde y)=
\sum_{\mu(\tilde z)=\mu(\tilde x)-1} \overline{\mathcal
M}^{\nu(\tilde x\tilde z)}(\tilde x,\tilde z)\times
\overline{\mathcal M}^{\nu(\tilde z\tilde y)}(\tilde z,\tilde
y).\end{array}\right.$$
  Hence we can replace $\sharp
C(\overline{\mathcal  M}^{\nu(\tilde x\tilde y)}(\tilde x,\tilde
y))$ by $\sharp\overline{\mathcal  M}^{\nu(\tilde x\tilde
y)}(\tilde x,\tilde y)$ in (1.2).\vspace{2mm}

Now we begin to construct all
 virtual moduli cycles $\overline{\mathcal  M}_\pm^{\nu^\pm}(\tilde x; H, J)$
in (2.9) and (2.10)  which are compatible with the relative
virtual cycles used in (1.2). We only consider $\overline{\mathcal
M}_+^{\nu^+}(\tilde x; H, J)$. Denote by $T\overline{\mathcal
M}^{\nu^+}_+(\tilde x; H, J)$ and $B\overline{\mathcal
M}^{\nu^+}_+(\tilde x; H, J)$ its top strata and $1$-codimensional
strata.
 If $\dim\overline{\mathcal  M}_+^{\nu^+}(\tilde x; H, J)=\dim M-\mu(\tilde
x)>0$, then any element $f_+$ of $B{\mathcal M}_+^{\nu^+_I}(\tilde
x; H, J)$ must have a form $f_+=(h_+, g)$. Here $g\in
\widetilde{\mathcal  M}^{\nu_2}_{I_2}(\tilde y, \tilde x; H, J)$
and $\mu(\tilde y)-\mu(\tilde x)=1$, and $h_+\in T{\mathcal
M}_+^{\nu^+_{I_1}}(\tilde y; H, J)$. So we can construct each
$\overline{\mathcal  M}_+^{\nu^+}(\tilde x; H, J)$ inductively
with respect to $\dim M-\mu(\tilde x)$. In fact, if these have
been constructed for $\dim M-\mu(\tilde x)=0$, then for each
$\tilde x\in\tilde{\mathcal  P}(H)$ with $\dim M-\mu(\tilde x)=1$,
 we can use them and
 $\widetilde{\mathcal  M}^\nu(\tilde y, \tilde x; H, J)$ used in (1.2) to construct
 $$B{\mathcal  M}_+^{\nu^+}(\tilde x; H, J):=
 \cup_{\mu(\tilde y)=\mu(\tilde x)-1}T{\mathcal  M}_+^{\nu^+}(\tilde y; H, J)\times
\widetilde{\mathcal  M}^\nu(\tilde y, \tilde x; H, J).$$
 Next as done in the proof of (2.12)
 we can extend $B{\mathcal  M}_+^{\nu^+}(\tilde x; H, J)$ into
a derived family ${\mathcal  M}^{\nu^+}_+(\tilde x; H,
J)=\{{\mathcal M}^{\nu^+_I}_+(\tilde x; H, J)\,|\, I\in{\mathcal
N}^+(\tilde x)\}$ for a virtual moduli cycle $\overline{\mathcal
M}^{\nu^+}_+(\tilde x; H, J)$. Let $\{\overline{\mathcal
M}^{\nu^+}_+(\tilde x; H, J)\,|\,
 \tilde x\in\tilde{\mathcal  P}(H)\}$
be all virtual moduli cycles constructed by induction. Then it holds that
$$B{\mathcal  M}_+^{\nu^+}(\tilde x; H, J)=\cup_{\mu(\tilde y)=\mu(\tilde x)-1}
T{\mathcal  M}_+^{\nu^+}(\tilde y; H, J)\times\widetilde{\mathcal
M}^\nu(\tilde y, \tilde x; H, J).$$ Consider the natural
orientations on them again we can write
$$B\overline{\mathcal  M}^{\nu^+}_+(\tilde x; H, J)=
\sum_{\mu(\tilde y)=\mu(\tilde x)+1} \sharp(C(\overline{\mathcal
M}^\nu(\tilde y, \tilde x)))\cdot T\overline{\mathcal
M}^{\nu^+}_+(\tilde y; H, J).\leqno(2.36)$$
 Here the
rational numbers $\sharp(C(\overline{\mathcal  M}^\nu(\tilde y,
\tilde x)))$ is as in (1.2), and the identity in (2.36) is
understood as follows: Since the sums at the right side of (2.36)
are finite we can take $L>0$ to be the smallest common multiple of
denominators of all rational numbers $\sharp(C(\overline{\mathcal
M}^\nu(\tilde y, \tilde x)))$. Then (2.36) is equivalent to
$$
L\cdot B\overline{\mathcal
M}^{\nu^+}_+(\tilde x; H, J)=\sum_{\mu(\tilde y)=\mu(\tilde x)+1}
(L\cdot\sharp(C(\overline{\mathcal  M}^\nu(\tilde y, \tilde
x))))\cdot T\overline{\mathcal  M}^{\nu^+}_+(\tilde y; H, J).
$$
For this identity the left side is understood as the disjoint
union of $L$ copies of $B\overline{\mathcal M}^{\nu^+}_+(\tilde x;
H, J)$, and the right side is also the disjoint union that
contains, for each $\tilde y=[y, u]$ with $\mu(\tilde
y)=\mu(\tilde x)+1$, $\;L\cdot\sharp(C(\overline{\mathcal
M}^\nu(\tilde y, \tilde x)))\;$ copies of $\;T\overline{\mathcal
M}^{\nu^+}_+(\tilde y; H, J)\;$ as $\;\sharp(C(\overline{\mathcal
M}^\nu(\tilde y, \tilde x)))>0$, and
$\;(-L)\cdot\sharp(C(\overline{\mathcal M}^\nu(\tilde y, \tilde
x)))\;$ copies of $\;T\overline{\mathcal M}^{\nu^+}_+(\tilde y; H,
J)^\ast\;$ as $\;\sharp(C(\overline{\mathcal  M}^\nu(\tilde y,
\tilde x)))<0$, where $T\overline{\mathcal  M}^{\nu^+}_+(\tilde y;
H, J)^\ast$ is $T\overline{\mathcal  M}^{\nu^+}_+(\tilde y; H, J)$
with the orientation reversed. In other words, it is understood as
the topological sum. With the same methods we can construct the
compatible $\overline{\mathcal  M}_-^{\nu^-}(\tilde x; H, J)$ and
obtain
$$B\overline{\mathcal  M}^{\nu^-}_-(\tilde x; H, J)=
\sum_{\mu(\tilde z)=\mu(\tilde x)-1} \sharp(C(\overline{\mathcal
M}^\nu(\tilde x, \tilde z)))\cdot T\overline{\mathcal
M}^{\nu^-}_-(\tilde z; H, J).\leqno(2.37)$$

\section{The intersections  of virtual moduli cycles with stable and
unstable manifolds}\label{3}

For the materials on Morse homology the readers may refer to
\cite{AuB, Sch1} and \cite{Sch4}. Given a Morse-Smale pair $(h,
g)$ on $M$ with $h\in{\mathcal  O}(h_0)$ the stable and unstable
manifolds of a critical point $a\in{\rm Crit}(h)$ are given by
\begin{eqnarray*}
&&W^s(a, h, g):=\{\gamma: [0, \infty)\to M\,|\,
\dot\gamma(s)+\nabla^gh(\gamma(s))=0, \;\gamma(\infty)=a\},\\
&&W^u(a, h, g):=\{\gamma:(-\infty,0]\to M\,|\,
\dot\gamma(s)+\nabla^gh(\gamma(s))=0, \;\gamma(-\infty)=a\}.
\end{eqnarray*}
There are two obvious evaluations
$$\left.\begin{array}{ll}
 E^s_{a}: W^s(a, h, g)\to M,\;\gamma\mapsto \gamma(0)\quad{\rm and}\\
E^u_a: W^u(a, h,  g)\to M,\;\gamma\mapsto \gamma(0).
\end{array}\right.\leqno(3.1)$$
Both are also smooth embeddings into $M$. Throughout this paper we
fix orientations for all unstable manifolds $W^u(a, h, g)$, then
the orientation of $M$ induces orientations for $W^s(a, h, g)$ and
$W^u(a, h, g)\cap W^s(b, h, g)$. We wish to study the intersection
numbers of the maps in (2.11) and (3.1) under some conditions. As
usual we need their compactifications of the following versions.
Consider  the disjoint union ${\overline W}^u(a, h, g):=W^u(a, h,
g)\cup S{\overline W}^u(a, h, g)$, where $S{\overline W}^u(a, h,
g)$ is the disjoint unions
$$\cup\widehat M_{a, a_1}(h, g)\times\cdots\times
\widehat M_{a_{i-1}, a_i}(h, g)\times W^u(a_i, h, g)$$ for all
critical points $a_0=a,\cdots, a_i$ with the Morse indexes
$\mu(a_0)>\cdots>\mu(a_i)$. Similarly, in the disjoint union
${\overline W}^s(a, h, g):=W^s(a, h, g)\cup S{\overline W}^s(a, h,
g)$ the second set $S{\overline W}^s(a, h,  g)$ is the disjoint
unions
$$\cup W^s(a_i, h, g)\times\widehat M_{a_i, a_{i-1}}(h, g)\times\cdots\times
\widehat M_{a_1, a}(h, g)$$ for critical points $a_0=a,\cdots,
a_i$ with the Morse indexes $\mu(a_0)<\cdots<\mu(a_i)$. As usual
the compactness and gluing arguments in Morse homology (see
\cite{AuB} and \cite{Sch1}) give:\vspace{2mm}

\noindent{\bf Lemma 3.1}.\hspace{2mm}{\it The sets ${\overline
W}^u(a, h,  g)$ and ${\overline W}^s(a, h, g)$ may be topologized
with a natural way so that they are the compactifications of
$W^u(a, h, g)$ and $W^s(a, h, g)$ respectively, and that
$\partial{\overline W}^u(a, h, g)=S{\overline W}^u(a, h,  g)$ and
$\partial{\overline W}^s(a, h, g)=S{\overline W}^s(a, \\h, g)$.
Moreover these compactified spaces both have the structure of a
manifold with corners, and maps $E^u_a$ and $E^s_a$ may smoothly
extend to them, denoted by $\bar E^u_a$ and $\bar E^s_a$, which
also give natural injective immersions from  these two
compactified spaces into $M$.}\vspace{2mm}

The following lemma, slightly different from Theorem 4.9 in
\cite{Sch4},  may be easily proved (cf. \cite{Lu2}).\vspace{2mm}

\noindent{\bf Lemma 3.2.}\hspace{2mm}{\it Let ${\mathcal  R}$ be
the set of all Riemannian metrics on $M$. Then for any smooth map
$\chi$ from a smooth manifold $V$ to $M$ there exists a Baire
subset $({\mathcal  O}\times{\mathcal  R})_{\rm
reg}\subset{\mathcal O}(h_0)\times{\mathcal  R}$ such that for
each pair $(h, g)\in ({\mathcal  O}(h_0)\times{\mathcal  R})_{\rm
reg}$ and $a\in{\rm Crit}(h)$ the maps $E^s_a$ and $E^u_a$ are
transverse to $\chi$. Consequently, the spaces
\begin{eqnarray*}
&&{\mathcal  M}^s_{\chi, a}(h, g):=\{(p,\gamma)\in V\times W^s(a,
h, g)\,
|\, \chi(p)=E^s_a(\gamma)=\gamma(0)\},\\
&&{\mathcal  M}^u_{\chi, a}(h, g):=\{(p,\gamma)\in V\times W^u(a,
h, g)\, |\, \chi(p)=E^u_a(\gamma)=\gamma(0)\}
\end{eqnarray*}
are respectively smooth manifolds of $\dim V-\mu(a)$ and $\dim V+
\mu(a)-2n$.}\vspace{2mm}

As usual, if $V$, $M$, $W^s(a, h, g)$ and $W^u(a, h, g)$  are
oriented then  ${\mathcal  M}^s_{\chi;a}(h, g)$ and ${\mathcal
M}^u_{\chi;a}(h, g)$ have the natural induced orientations.
Specially, if ${\mathcal  M}^s_{\chi;a}(h, g)$ and ${\mathcal
M}^u_{\chi;a}(h, g)$ are of dimension zero  we have the oriented
intersection numbers $\chi\cdot E^s_a$ and $\chi\cdot E^u_a$ which
counts the algebraic sum of the oriented points in ${\mathcal
M}^s_{\chi;a}(h, g)$ and ${\mathcal  M}^u_{\chi;a}(h, g)$
respectively. Note that the $1$-codimensional stratum of
${\overline W}^u(a, h, g)$ (resp. ${\overline W}^s(a, h, g)$) is
given by
$$\sum_{\mu(b)=\mu(a)-1}n(a, b)\cdot W^u(b, h, g)\quad \biggl({\rm
resp.} \sum_{\mu(c)=\mu(a)+1}n(c, a)\cdot W^s(c, h, g)\biggr),$$
where $n(a, b)$ and $n(c, a)$ are as in (1.6).\vspace{2mm}

 Since there are only countable manifolds involved
in our arguments using Claim A.1.11 in \cite{LO} we may always fix
a $h\in{\mathcal O}(h_0)$ such that all transversal arguments hold
for a generic Riemannian metric $g$ on $M$. By lemmas 3.1 and 3.2
and these remarks we immediately obtain:\vspace{2mm}

\noindent{\bf Proposition 3.3.}\hspace{2mm}{\it If $\mu(a)=
\mu(\tilde x)$ then for a generic Riemannian metric $g$ on $M$ the
maps $\bar E^u_{a}$ and ${\rm EV}_+^{\nu^+}(\tilde x)$, and $\bar
E^s_{a}$ and ${\rm EV}_-^{\nu^-}(\tilde x)$ are transversal, and
their intersection numbers are the well-defined rational numbers
$\bar E^u_a\cdot{\rm EV}_+^{\nu^+}(\tilde x)$ and $\bar
E^s_a\cdot{\rm EV}_-^{\nu^-}(\tilde x)$.}\vspace{2mm}

From now on, for the sake of clearness we denote by
$$\left.\begin{array}{ll}
n_+^{\nu^+}(a, \tilde x)\equiv n_+^{\nu^+}(a, \tilde x; H, J; h,
g):=
\bar E^u_a\cdot{\rm EV}_+^{\nu^+}(\tilde x),\\
n_-^{\nu^-}(a, \tilde x)\equiv n_-^{\nu^-}(a,\tilde x; H, J; h,
g):= \bar E^s_a\cdot{\rm EV}_-^{\nu^-}(\tilde x).
\end{array}\right.\leqno(3.2)$$
The standard arguments show that these numbers are independent of
small regular $\nu^\pm$. Later, we state no longer this and often
omit $\nu^\pm$ without occurrence of confusions. It easily follows
from (2.7) that: \vspace{2mm}

\noindent{\bf Proposition 3.4.}\hspace{2mm}{\it If $n_\pm(a, [x,
v])\ne 0$  then $\pm\int_{D^2}v^\ast\omega\ge-2\max
|H|$.}\vspace{2mm}

Using (2.36), (2.37) and lemmas 3.1 and 3.2 we may obtain the
following two results:\vspace{2mm}

\noindent{\bf Proposition 3.5.}\hspace{2mm}{\it If
$\mu(a)-\mu(\tilde x)=1$ then for a generic  Riemannian metric $g$
on $M$ the fibre product\vspace{-1mm}
$${\overline W}^u(a, h, g)\times_{\bar E^u_a={\rm EV}^{\nu^+}_+(\tilde x)}
\overline{\mathcal  M}^{\nu^+}_+(\tilde x; H, J)$$ is still a
collection of compatible local cornered smooth manifolds of
dimension $1$ and with the natural orientations. Its boundary is
given by\vspace{-1mm}
\begin{eqnarray*}\left.\begin{array}{ll}
\bigl(\cup_{\mu(b)=\mu(a)-1}n(a, b)\cdot \bigl(W^u(b, h, g)
\times_{E^u_b={\rm EV}^{\nu^+}_+(\tilde x)}\overline{\mathcal
M}^{\nu^+}_+(\tilde x; H, J)\bigr)
\bigr)\cup\\
\!\!\!\bigl(\!-\cup_{\mu(\tilde y)=\mu(\tilde x)+1}
\sharp(C(\overline{\mathcal  M}^\nu(\tilde y,\tilde x)))\cdot
\bigr({\overline W}^u(a, h, g)\times_{\bar E^u_a={\rm
EV}^{\nu^+}_+(\tilde y)} \overline{\mathcal M}^{\nu^+}_+(\tilde y;
H, J)\bigr)\!\bigr).
\end{array}\right.
\end{eqnarray*}
Notice that the projection onto $M\times{\mathcal W}_+(\tilde x;
H, J)$ of this set is a finite set. Let us denote by
$\sharp(\partial({\overline W}^u(a, h, g) \times_{\bar E^u_a={\rm
EV}^{\nu^+}_+(\tilde x)} \overline{\mathcal M}^{\nu^+}_+(\tilde x;
H, J)))$ the number of elements of this finite set counted with
appropriate signs and rational weights. Then this number must be
zero, and it follow that
$$\sum_{\mu(b)=\mu(a)-1}n(a, b)\cdot
n_+^{\nu^+}(b, \tilde x)\\
=\sum_{\mu(\tilde y)=\mu(\tilde x)+1}n^{\nu^+}_+(a, \tilde y)\cdot
\sharp(C(\overline{\mathcal  M}^\nu(\tilde y, \tilde x))).$$ }

Remark that  the fibre product in Proposition 3.5 has boundary.
But we assume its dimension to be $1$. Hence its boundary agrees
with the $1$-codimensional stratum of the fibre product. This
remark is still valid for the following proposition:\vspace{2mm}

\noindent{\bf Proposition 3.6.}\hspace{2mm}{\it If $\mu(\tilde
x)-\mu(a)=1$
 then for  a generic Riemannian metric $g$ on $M$ the fibre product\vspace{-1mm}
$${\overline W}^s(a, h, g)\times_{\bar E^s_a={\rm EV}^{\nu^-}_-(\tilde x)}
\overline{\mathcal  M}^{\nu^-}_-(\tilde x; H, J)$$ is a collection
of compatible local cornered smooth manifolds of dimension $1$ and
with the natural orientations. Its boundary is given
by\vspace{-1mm}
\begin{eqnarray*}\left.\begin{array}{ll}
\bigl(\cup_{\mu(a)=\mu(b)-1}n(b, a)\cdot\bigl(W^s(b, h, g)
\times_{E^u_b={\rm EV}^{\nu^-}_-(\tilde x)}\overline{\mathcal
M}^{\nu^-}_-(\tilde x; H, J)\bigr)\bigr)\cup\\
\!\!\!\bigl(\!-\cup_{\mu(\tilde z)=\mu(\tilde x)-1}
\sharp(C(\overline{\mathcal  M}^\nu(\tilde x, \tilde z)))\cdot
\bigl({\overline W}^u(a, h, g)\times_{\bar E^u_a={\rm
EV}^{\nu^-}_-(\tilde z)} \overline{\mathcal M}^{\nu^-}_-(\tilde z;
H, J)\bigr)\!\bigr).\end{array}\right.
\end{eqnarray*}
Consequently, $\sharp(\partial({\overline W}^s(a, h, g)
\times_{\bar E^u_a={\rm EV}^{\nu^-}_-(\tilde x)}
\overline{\mathcal M}^{\nu^-}_-(\tilde x; H, J)))=0$ implies that
$$\sum_{\mu(\tilde z)=\mu(\tilde x)-1}
\sharp(C(\overline{\mathcal  M}^\nu(\tilde x, \tilde z)))\cdot n^{\nu^-}_-(a, \tilde z)\\
=\sum_{\mu(b)=\mu(a)+1}n_-^{\nu^-}(b, \tilde x)\cdot n(b, a).$$}

\noindent{\bf Remark 3.7.}\hspace{2mm}Using Proposition 3.4 it is
 easily  proved  that $\Phi$ and $\Psi$ are
  $\Lambda_\omega$-module homomorphisms( see \cite{Lu2}).
Moreover, using Propositions 3.4, 3.5, 3.6 and (2.35) we may prove
that $\Phi\circ\partial^Q_k=\partial^F_k\circ\Phi$ and
$\Psi\circ\partial^F_k=\partial^Q_k\circ\Psi$ for every integer
$k$. That is, $\Phi$ and $\Psi$ also induce the homomorphisms
between two  corresponding homology groups( see \cite{Lu2}).

\section{Proof of  Theorem 1.1}\label{4}

Theorem 1.1 may follow from the following Theorems 4.1 and 4.9
immediately.\vspace{2mm}

\noindent{\bf Theorem 4.1.}\hspace{2mm}{\it $\Psi\circ\Phi$ is
chain homotopy equivalent to the identity. Consequently, $\Phi$
induces an injective $\Lambda_\omega$-module homomorphism from
$QH_\ast(h, g;\mbox{\Bb Q})$ to $HF_\ast(M,\omega;\\ H,
J,\nu;\mbox{\Bb Q})$.}\vspace{2mm}

\noindent{\it Proof}.\hspace{2mm}We shall prove it in four steps.

\noindent{\it Step 1.}\hspace{2mm} For every $\langle a,
A\rangle\in{\rm Crit}(h)\times\Gamma$ with  $\mu(\langle a,
A\rangle)=k$ we have
$$
\Psi\circ\Phi(\langle a, A\rangle)=\sum_{\mu(\langle b,
B\rangle)=k} m^\nu_{+,-}(\langle a, A\rangle;\langle b, B\rangle)
\langle b, B\rangle,
$$
 where $m^\nu_{+,-}(\langle a,
A\rangle;\langle b, B\rangle)$ is given by
$$
\sum_{\mu(\tilde x)=k}n^{\nu^+}_+(a, \tilde x\sharp(-A))\cdot
n^{\nu^-}_-(b, \tilde x\sharp(-B)).
\leqno(4.1)$$
Notice here that
$\langle a, A\rangle$ and $\langle b, B\rangle$ satisfy
$$\mu(a)-\mu(b)+ 2c_1(B-A)=0.\leqno(4.2)$$
Firstly, we claim that the sum in (4.1) is finite. In fact, if
 $\tilde x\in {\tilde{\mathcal  P}}_k(H)$ is such that
$$n^{\nu^+}_+(a, \tilde x\sharp(-A))\cdot n^{\nu^-}_-(b, \tilde x\sharp(-B))\ne 0,
\leqno(4.3)$$
by Proposition 3.4 we have
$$
\omega(A)-2\max|H|\le\int_{D^2}v^\ast\omega\le\omega(B)+ 2\max|H|.
\leqno(4.4)$$
It easily follows  that the number of such $\tilde
x\in\tilde{\mathcal P}(H)$ is finite. Let them be $\tilde
x_1,\cdots, \tilde x_s$. Then for a generic $(h,g)\in{\mathcal
O}(h_0)\times{\mathcal  R}$, it holds that\vspace{-1mm}
$$
m^\nu_{+,-}(\langle a, A\rangle;\langle b, B\rangle)= \sum_{i=1}^s
n^{\nu^+}_+(a, \tilde x_i\sharp(-A))\cdot n^{\nu^-}_-(b,  \tilde
x_i\sharp(-B)).
\leqno(4.5)$$
 Note that for
each $\tilde x\in\tilde{\mathcal  P}_k(H)$ the product
$n^{\nu^+}_+(a, \tilde x\sharp(-A))\cdot n^{\nu^-}_-(b, \tilde
x\sharp(-B))$
 can be explained as the rational intersection number
$$(\bar E^u_a\times\bar E^s_b)\cdot
({\rm EV}^{\nu^+}_+(\tilde x\sharp(-A))\times {\rm
EV}^{\nu^-}_-(\tilde x\sharp(-B)))\leqno(4.6)$$
 of the product map
$${\rm EV}^{\nu^+}_+(\tilde x\sharp(-A))\times
{\rm EV}^{\nu^-}_-(\tilde x\sharp(-B))$$
from $\overline{\mathcal
M}_+^{\nu^+}(\tilde x\sharp(-A); H, J)\times \overline{\mathcal
M}_-^{\nu^-}(\tilde x\sharp(-B); H, J)$ to $M\times M$ and the
evaluation
$$\bar E^u_a\times\bar E^s_b: {\overline W}^u(a, h, g)\times
{\overline W}^s(b, h, g)\to M\times M\leqno(4.7)$$
 for a generic
$(h, g)\in{\mathcal  O}(h_0)\times{\mathcal  R}$. But the
intersection number in (4.6) is the number of the oriented points
with rational weights in the fibre product
$$({\overline W}^u(a, h, g)\times
{\overline W}^s(b, h, g))\times_R (\overline{\mathcal
M}_+^{\nu^+}(\tilde x\sharp(-A))\times \overline{\mathcal
M}_-^{\nu^-}(\tilde x\sharp(-B))\!). \leqno(4.8)$$
 Here $R$ is representing
 ${\rm EV}^{\nu^+}_+\times{\rm EV}^{\nu^-}_-=
\bar E^u_a\times\bar E^s_b$, and we have omitted $H, J$ in
$\overline{\mathcal M}_\pm^{\nu^\pm}(\cdot\,; H, J)$. However,
because of dimension relations the fibre product in (4.8) is an
empty set for a generic
 $(h, g)\in{\mathcal O}(h_0)\times{\mathcal R}$  if
 $\mu([x, v])\ne k$. As usual we still understand the intersection number
being zero in this case. Hence (4.1) and (4.5) become
$$\left.\begin{array}{ll}\quad m^\nu_{+,-}(\langle a,
A\rangle;\langle b, B\rangle)\\
=\!\! \sum_{\tilde x\in\tilde{\mathcal  P}(H)} \!(\bar
E^u_{a}\times\bar E^s_{b})\cdot ({\rm EV}^{\nu^+}_+(\tilde
x\sharp(-A))\times {\rm EV}^{\nu^-}_-(\tilde x\sharp(-B))\!).
\end{array}\right.\leqno(4.9)$$

\noindent{\it Step 2.}\hspace{3mm}We need to give another
explanation of the intersection number in (4.6). To this goal we
note that  for a given $D\in\Gamma$ a positive disk
$u_+\in{\mathcal M}_+(\tilde y; H, J)$ and a negative disk
$u_-\in{\mathcal  M}_-(\tilde y\sharp(-D); H, J)$  may be glued
into a sphere in $M$ along $y$, denoted by $u_+\sharp u_-$. Using
(2.2) and (2.4) one easily checks that it is a representative of
$D$. Let us denote by ${\mathcal M}_D(\tilde y; H, J)$ by the
space of all such $u_+\sharp u_-$.
 Clearly, it can be identified with the product space
${\mathcal  M}_+(\tilde y; H, J)\times{\mathcal  M}_-(\tilde
y\sharp(-D); H, J)$. Therefore, its virtual dimension is equal to
$\dim M+ 2c_1(D)$.  We denote by
$${\mathcal  M}_D(H, J)=\cup_{\tilde y\in
\tilde{\mathcal  P}(H)}{\mathcal  M}_D(\tilde y; H, J).
$$
As in (4.3) and (4.4) it is easy to prove that the union at the
right side is only finite union. That is, there only exist
finitely many $\tilde y_1,\cdots, \tilde y_r$ in $\tilde{\mathcal
P}(H)$ such that ${\mathcal  M}_D(\tilde y_i; H, J)\ne\emptyset$
for $i=1,\cdots, r$. Thus
$${\mathcal  M}_D(H, J)=\cup^r_{i=1}{\mathcal  M}_D(\tilde y_i; H, J).\leqno(4.10)$$
In particular,  (4.3) and (4.4) imply that
$${\mathcal  M}_{B-A}(H, J)=\cup^s_{i=1}{\mathcal  M}_{B-A}(\tilde x_i\sharp(-A); H, J).
\leqno(4.11)$$
Note that the right sides of both (4.10) and (4.11)
are all disjoint unions. In order to compactify
 ${\mathcal  M}_D(H, J)$ we introduce:\vspace{2mm}

\noindent{\bf Definition 4.2.}\hspace{2mm} Given  $\tilde
y\in\tilde{\mathcal  P}(H)$, $D\in\Gamma$, a semistable ${\mathcal
F}$-curve $(\Sigma,\underline l)$ with at least two principal
components (cf. Def.2.1),  a continuous map
$f:\Sigma\setminus\{z_2,\cdots, z_{N_p}\\
\}\to M$ is called a {\bf
stable $(J, H)$-broken solution with a joint $\tilde y$ and of
class} $D$ if we may divide $(\Sigma,\underline l)$ into two
semistable ${\mathcal F}$-curve $(\Sigma_+,\underline l^+)$ and
$(\Sigma_-,\underline l^-)$ at some double point $z_{i_0}$ between
two principal components, $2\le i_0\le N_p$, such that
$f|_{\Sigma_+}$ and $f|_{\Sigma_-}$ are the stable $(J, H)_+$-disk
with cap $\tilde y$ and stable $(J, H)_-$-disk with cap $\tilde
y\sharp(-D)$ respectively. Furthermore, a tuple $(f,\Sigma,
\underline l)$ is called a {\bf stable $(J, H)$-broken solution of
class $D$} if it is stable $(J, H)$-broken solution with a joint
$\tilde y$ and of class $D$ for some $\tilde y\in\tilde{\mathcal
P}(H)$.\vspace{2mm}

It is easily seen that if a stable $(J, H)$-broken solution $(f,
\Sigma,\underline l)$ with a joint $\tilde y$ and of a class $D$
has at least three principal components then it is also such a
broken solution with another joint different from $\tilde y$ and
of class $D$. As usual we may define the equivalence class of such
a broken solution in an obvious way. But it should be noted that
two equivalence classes of a given stable $(J, H)$-broken solution
$(f,\Sigma,\underline l)$ with a joint $\tilde y$ as a stable $(J,
H)$-broken solution with a joint $\tilde y$ and as a stable $(J,
H)$-broken solution are same. We still denote by $\langle
f\rangle$ the equivalence class of $(f,\Sigma,\underline l)$.
Define the energy of such a $\langle f\rangle$ by $E_D(\langle
f\rangle)=E_D(f)=E_+(f|_{\Sigma_+})+ E_-(f|_{\Sigma_-})$. Then
(2.6)  yields that $E_D(\langle f\rangle)\le\omega(D)+ 2\max|H|$.
Let us denote by $\overline{\mathcal  M}_D(\tilde y; H, J)$ and
$\overline{\mathcal  M}_D(H, J)$ the spaces of all equivalence
classes of the two kinds of maps respectively. Then
$$\overline{\mathcal  M}_D(H, J)=\cup_{\tilde y\in
\tilde{\mathcal  P}(H)}\overline{\mathcal  M}_D(\tilde y; H,
J).\leqno(4.12)$$
 However, one should also note that
$\overline{\mathcal  M}_D(\tilde y; H, J)$ and
 $\overline{\mathcal M}_D(\tilde z; H, J)$ have probably nonempty intersection for two
different $\tilde y$ and $\tilde z$. Thus we mayn't affirm the
above unions to be the disjoint unions. But each stable
 $(J, H)$-broken solution has at most $2n=\dim M$ joints. As in
(4.3) and (4.4) we can prove the union at the right side of (4.12)
is actually a finite unions, i.e., other
 $\overline{\mathcal M}_D(\tilde y; H, J)$ are empty except finitely many
 $\tilde y\in\tilde{\mathcal  P}(H)$, saying $\tilde y_1,\cdots, \tilde y_t, t\ge r$
  because of (4.10) and the fact that
  ${\mathcal M}_D(\tilde y; H, J)\subset\overline{\mathcal  M}_D(\tilde y; H, J)$.
Then we get that
$$\overline{\mathcal  M}_D(H, J)=\cup^t_{i=1}
\overline{\mathcal  M}_D(\tilde y_i; H, J).\leqno(4.13)$$
 Note that by the above assumptions if $t>r$ then
 ${\mathcal  M}_D(\tilde y_i; H,J)=\emptyset$ for $r<i\le t$.
 Unlike in (4.10) the union at the right side is not necessarily disjoint.
 However, as in Proposition 2.5 one easily shows that the spaces
 $\overline{\mathcal  M}_D(\tilde y; H, J)$ and
$\overline{\mathcal  M}_D(H, J)$ equipped with the weak
$C^\infty$-topology are the Hausdorff compactifications of
 ${\mathcal  M}_D(\tilde y; H, J)$ and ${\mathcal  M}_D(H, J)$ respectively
(cf.[LiuT1]). Moreover, the notions of the corresponding stable
$L^p_k$-broken maps may also be introduced in the similar way. Let
$\overline{\mathcal  M}^\nu_D(\tilde y; H, J)$ and $
\overline{\mathcal  M}^\nu_D(H, J)$ be the corresponding virtual
moduli cycles to them respectively. Consider the evaluation maps
$${\rm EV}^\nu_D(\tilde y): \overline{\mathcal  M}^\nu_D(\tilde y;H, J)\to M\times M
\;\;{\rm and}\;\;{\rm EV}^\nu_D: \overline{\mathcal  M}^\nu_D(H,
J)\to M\times M $$ given by
 ${\rm EV}^\nu_D(\tilde y)(f)=(f(z_-), f(z_+))$ and ${\rm EV}^\nu_D(f)=(f(z_-), f(z_+))$,
where $z_-$ and $z_+$ two double points on the chain of the
principal components of the domain of $f$ (cf. Def. 2.1).  For the
map $\bar E^u_a\times\bar E^s_b$ in (4.7), if
$$\mu(a)-\mu(b)+ 2c_1(D)=0,\leqno(4.14)$$
then for a generic $(h,g)\in {\mathcal  O}(h_0)\times{\mathcal R}$
we get two (rational) intersection numbers
$$\left.\begin{array}{ll} n_D^\nu(a, b, \tilde y):=(\bar
E^u_a\times\bar
E^s_b)\cdot {\rm EV}^\nu_D(\tilde y)\quad{\rm and}\\
n_D^\nu(a, b):=(\bar E^u_a\times\bar E^s_b)\cdot {\rm EV}^\nu_D,
\end{array}\right.\leqno(4.15)$$
 which are independent of generic $\nu$ and $(h,g)\in
{\mathcal O}(h_0)\times{\mathcal  R}$.\vspace{2mm}

\noindent{\bf Proposition 4.3.}\hspace{2mm}{\it Under the above
assumptions it
 holds that
$$n_D^\nu(a, b, \tilde
y)=n^{\nu^+}_+(a, \tilde y)\cdot n^{\nu^-}_-(b, \tilde
y\sharp(-D)),\leqno(4.16)$$\vspace{-5mm}
$$\left.\begin{array}{ll}n_D^\nu(a, b)=\sum^r_{i=1}n_D^\nu(a, b, \tilde
y_i).\end{array}\right.\leqno(4.17)$$}

From these it easily follows that (4.5) and (4.9) become
$$m^\nu_{+,-}(\langle a,
A\rangle;\langle b, B\rangle)=(\bar E^u_a\times\bar E^s_b)\cdot
{\rm EV}^\nu_{B-A}.\leqno(4.18)$$

\noindent{\it Proof of Proposition 4.3.}\hspace{2mm}We first prove
(4.16).
 By the definition in (3.2) the
rational numbers $n^{\nu^+}_+(a, \tilde y)$ and $n^{\nu^-}_-(b,
\tilde y\sharp(-D))$ are independent of the choices of generic
small $\nu^+$ and $\nu^-$. Therefore, we only need to prove (4.16)
for suitable regular $\nu$ and $\nu^\pm$. Notice that
$\overline{\mathcal M}_D(\tilde y; H, J)$ may be identified with
the product
$$\overline{\mathcal  M}_+(\tilde y; H, J)\times\overline{\mathcal  M}_-(\tilde y\sharp(-D); H, J)$$
by the map $\langle f, \Sigma,\underline l\rangle\mapsto (\langle
f|_{\Sigma_+}, \Sigma_+, \underline l^+\rangle, \langle
f|_{\Sigma_-}, \Sigma_-, \underline l^-\rangle)$. More generally,
let\vspace{-1mm}
$${\mathcal  B}^{p,k}_D(\tilde y,;H),\quad {\mathcal  B}^{p,k}_+(\tilde y;H)\quad{\rm and}\quad
{\mathcal  B}^{p,k}_-(\tilde y\sharp(-D);H)$$ be the corresponding
$L^p_k$-stable map spaces then the first one can be identified
with the product space
 ${\mathcal  B}^{p,k}_+(\tilde y;H)\times {\mathcal  B}^{p,k}_-(\tilde y\sharp(-D);H)$
 by the natural map\vspace{-1mm}
$${\rm GL}: (\langle f_+,\Sigma_+,\underline l^+\rangle, \langle f_-,\Sigma_-,
\underline l^-\rangle)\mapsto \langle f_+\cup_y f_-,
\Sigma_+\cup\Sigma_-, \underline l^+\cup\bar
l^-\rangle.\leqno(4.19)
$$
We shall prove that a virtual moduli cycle  $\overline{\mathcal
M}^\nu_D(\tilde y; H, J)$ for $\overline{\mathcal M}_D(\tilde y;
H, J)$ naturally induces the virtual moduli cycle
$\overline{\mathcal M}_+^{\nu^+}(\tilde y; H, J)$ for
$\overline{\mathcal  M}_+(\tilde y; H, J)$ and that
$\overline{\mathcal M}_-^{\nu^-}(\tilde y\sharp(-D); H, J)$ for
$\overline{\mathcal M}_-(\tilde y\sharp(-D); H, J)$ such that
\begin{eqnarray*}\left.\begin{array}{ll}
\hspace{-5mm}(4.20)\hspace{15mm}\overline{\mathcal  M}^\nu_D(\tilde y; H,J)=\\
\quad\quad\quad\{f_+\sharp_y f_- | (f_+, f_-)\in
\overline{\mathcal M}^{\nu^+}_+(\tilde y; H, J)\times
\overline{\mathcal M}^{\nu^-}_-(\tilde y\sharp(-D); H, J)\}.
\end{array}\right.
\end{eqnarray*}
Once it is proved. Note that for $f=f_+\sharp_y f_-\equiv (f_+,
f_-)\in \overline{\mathcal M}^\nu_D(\tilde y; H, J)$,
$${\rm EV}^\nu_D(\tilde y)(f)=({\rm EV}^{\nu^+}_+(\tilde y)(f_+),\;
{\rm EV}^{\nu^-}_-(\tilde y\sharp(-D))(f_-)).
$$
By (4.15) and (4.20),  for a generic choice of $(h, g)\in
{\mathcal O}(h_0)\times{\mathcal  R}$
\begin{eqnarray*}
 n_D^\nu(a, b, \tilde y)\hspace{-5mm}&&=(\bar E^u_a\times\bar E^s_b)\cdot
{\rm EV}^\nu_D(\tilde y)\\
&&=(\bar E^u_a\times\bar E^s_b)\cdot({\rm EV}^{\nu^+}_+(\tilde
y)\times
{\rm EV}^{\nu^-}_-(\tilde y\sharp(-D)))\\
&&=(\bar E^u_a\cdot {\rm EV}^{\nu^+}_+(\tilde y))\times
(\bar E^s_b\cdot {\rm EV}^{\nu^-}_-(\tilde y\sharp(-D)))\\
&&=n_+^{\nu^+}(a, \tilde y)\cdot n_-^{\nu^-}(b, \tilde
y\sharp(-D)).
\end{eqnarray*}
Namely, (4.16) holds. In order to prove (4.20) we follow
\cite{LiuT1, LiuT2} and \cite{LiuT3}  to choose finitely many
points $\langle f^+_i\rangle=\langle f^+_i,\Sigma^+_i, \bar
l^+_i\rangle\in\overline{\mathcal M}_+(\tilde y; H, J)$ and
$\langle f^-_j\rangle=\langle f^-_j,\Sigma^-_j, \bar
l^-_j\rangle\in \overline{\mathcal M}_-(\tilde y\sharp(-D); H,
J)$, and their open neighborhoods $\;W^+_i=W^+_i(\langle
f^+_i\rangle)\;$ in $\;{\mathcal B}_+^{p,k}(\tilde y;H)\;$ and
those $\;W^-_j=W^-_j(\langle f^-_j\rangle)\;$ in
\linebreak
${\mathcal B}_-^{p,k}(\tilde y\sharp(-D);H)$, $i=1,\cdots, m_+$,
$j=1,\cdots, m_-$, such that  the following hold:

\begin{enumerate}
\item[(i)] ${\mathcal  W}^+(\tilde y):=\cup^{m_+}_{i=1}W^+_i$ and
${\mathcal  W}^-(\tilde y\sharp(-D)):=\cup^{m_-}_{j=1}W^-_j$ are
respectively open neighborhoods of $\overline{\mathcal M}_+(\tilde
y; H, J)$ in ${\mathcal  B}_+^{p,k}(\tilde y;H)$ and those
of\linebreak
 $\overline{\mathcal M}_-(\tilde y\sharp(-D); H, J)$
in
 ${\mathcal B}_-^{p,k}(\tilde y\sharp(-D);H)$.

\item[(ii)] $W^+_i$ and $W^-_j$ have  the uniformizers
$(\widetilde W^+_i, \Gamma^+_i, \pi^+_i)$ and $(\widetilde W^-_j,
\Gamma^-_j, \pi^-_j)$ respectively, $i=1,\cdots, m_+$,
$j=1,\cdots, m_-$.

\item[(iii)] $\{(\widetilde W^+_i, \Gamma^+_i,
\pi^+_i)\}^{m_+}_{i=1}$ and $\{(\widetilde W^-_j, \Gamma^-_j,
\pi^-_j)\}^{m_-}_{j=1}$ constitute the orbifold atlases on
${\mathcal W}^+(\tilde y)$ and
 ${\mathcal  W}^-(\tilde y\sharp(-D))$ respectively.

\item[(iv)] There exist the local orbifold bundles
$\widetilde{\mathcal  L}^+_i\to\widetilde W^+_i$ and
$\widetilde{\mathcal  L}^-_j\to\widetilde W^-_j$ with uniformizing
groups $\Gamma^+_i$ and $\Gamma^-_j$ respectively.

\item[(v)] There exist two obvious smooth sections
\begin{eqnarray*}
&&\bar\partial^+_{J,H}: {\mathcal  W}^+(\tilde y)\to{\mathcal
L}^+=
\cup^{m_+}_{i=1}{\mathcal  L}^+_i\;\;{\rm and}\\
&&\bar\partial^-_{J,H}: {\mathcal  W}^-(\tilde
y\sharp(-D))\to{\mathcal L}^- =\cup^{m_-}_{j=1}{\mathcal L}^-_j
\end{eqnarray*}
 such that their zero sets are
$\overline{\mathcal M}_+(\tilde y; H, J)$ and $\overline{\mathcal
M}_-(\tilde y\sharp(-D); H, J)$, respectively, and that they may
be lifted to the collections
$\{\bar\partial^{+,i}_{J,H}\}^{m_+}_{i=1}$ of
$\Gamma^+_i$-equivariant sections of $\widetilde{\mathcal
L}^+_i\to\widetilde W^+_i$ and those
$\{\bar\partial^{-,j}_{J,H}\}^{m_-}_{j=1}$ of
$\Gamma^-_j$-equivariant sections of
 $\widetilde{\mathcal L}^-_j\to\widetilde W^-_j$ respectively.
\end{enumerate}

 As defined above (2.13) we have corresponding index sets
$${\mathcal N}^\pm:=\{I\subset\{1,\cdots,
m^\pm\}\,|\, W^\pm_I= \cap_{k\in I}W^\pm_k\ne\emptyset\}.$$
 For
$I=\{i_1,\cdots, i_p\}\in{\mathcal  N}^+$ and $J=\{j_1,\cdots,
j_q\}\in{\mathcal  N}^-$ we denote by
$\Gamma_I^+=\prod^p_{k=1}\Gamma^+_{i_k}$ and $\Gamma_J^-=
\prod^q_{l=1}\Gamma^-_{j_l}$. Correspondingly, we have also the
fibre products
\begin{eqnarray*}
&&\widetilde W_I^{\Gamma_I^+}=\{ x^+_I=( x^+_i)_{i\in I}\,|\,
\pi^+_i(x^+_i)=\pi^+_{i^\prime}(x_{i^\prime})\in W^+_I,\;\forall
i, i^\prime\in I\},\\
&&\widetilde W_J^{\Gamma_J^-}=\{y^-_J=(y^-_j)_{j\in J}\,|\,
\pi^-_j(y^-_i)=\pi^-_{j^\prime}(y_{j^\prime})\in W^-_J,\;\forall
j, j^\prime\in J\}
\end{eqnarray*}
and the projections $\pi^+_I:\widetilde  W_I^{\Gamma_I^+}\to
W_I^+$ and $\pi^-_J:\widetilde  W_J^{\Gamma_J^-}\to W_J^-$. For
$I, I^\prime\in{\mathcal  N}^+$ with $I\subset I^\prime$ and $J,
J^\prime\in{\mathcal  N}^-$ with $J\subset J^\prime$ we have also
obvious projections $\pi^+_{I, I^\prime}$ and $\pi^-_{J,
J^\prime}$ satisfying
$\pi^+_{I^\prime}\circ\pi^+_{I,I^\prime}=\pi^+_I$ and
$\pi^-_{J^\prime}\circ\pi^-_{J,J^\prime}=\pi^-_J$. Following
[LiuT3] we may construct the desingularizations of $\widetilde
W^{\Gamma^+}$ and $\widetilde W^{\Gamma^-}$ as follows:
\begin{eqnarray*}
&&\widehat W^{\Gamma^+}\!=\Bigl\{\widehat W_I^{\Gamma^+},
\pi^+_{I, I^\prime}\, |\,I\!\subset I^\prime\!\in{\mathcal
N}^+\Bigr\}\;\;{\rm
and}\\
&& \widehat W^{\Gamma^-}\!=\Bigl\{\widehat W_J^{\Gamma^-},
\pi^+_{J, J^\prime}\,| \, J\!\subset J^\prime\!\in{\mathcal
N}^-\Bigr\}.
\end{eqnarray*}
 By (2.15) and (2.16), for every $I\subset{\mathcal N}^+$ (resp.
$J\subset{\mathcal  N}^-$) there exist subsets $V^+_I\subset
W^+_I$ (resp. $V^-_I\subset W^-_I$) such that:

\noindent{\rm (a)}\;\; $\overline{\mathcal  M}_+(\tilde y; H,
J)\subset\cup_{I\in{\mathcal  N}^+}V^+_I $ (resp.
$\overline{\mathcal  M}_-(\tilde y\sharp(-D); H, J)\subset
\cup_{\in{\mathcal  N}^-}V^-_J$);\vspace{1mm}

\noindent{\rm (b)}\;\; $V^+_{I^\prime}\cap V^+_I\ne\emptyset$
(resp. $V^-_{J^\prime}\cap V^-_J\ne\emptyset$) only when $I\subset
I^\prime$ (resp. $J\subset J^\prime$).

\noindent{Then} one can obtain  from $\widetilde W_I^{\Gamma_I^+}$
(resp. $\widetilde W_J^{\Gamma_J^-}$) the desingularization
$\widehat V_I^{\Gamma^+_I}$ (resp. $\widehat V_J^{\Gamma^-_J}$) of
the restriction of $\widetilde W_I^{\Gamma_I^+}$ (resp.
$\widetilde W_J^{\Gamma_J^-}$)  to $V_I^+$ (resp. $V_J^-$). Let us
denote by
\begin{eqnarray*}
&&\widehat V^{\Gamma^+}=\Bigl\{\widehat V_I^{\Gamma_I^+},
\pi^+_{I, I^\prime}\, \Bigm|\,I\subset
I^\prime\in{\mathcal N}^+\Bigr\}\;\;{\rm and}\\
&& \widehat V^{\Gamma^-}=\Bigl\{\widehat V_J^{\Gamma_J^-},
\pi^+_{J, J^\prime}\, \Bigm|\,J\subset J^\prime\in{\mathcal
N}^-\Bigr\}.
\end{eqnarray*}
With the
same way we may use $\widetilde{\mathcal L}^+_i$ and
$\widetilde{\mathcal L}^-_j$ to construct $\widetilde{\mathcal
L}_I^{\Gamma^+_I}$, $\widehat{\mathcal L}_I^{\Gamma^+_I}$ and
$\widetilde{\mathcal L}_J^{\Gamma^-_J}$, $\widehat{\mathcal
L}_J^{\Gamma^-_J}$ for $I\in{\mathcal  N}^+$ and $J\in{\mathcal
N}^-$, and obtain the following systems of the bundles
\begin{eqnarray*}
&& \bigl(\widetilde{\mathcal L}^{\Gamma^+}, \widetilde
W^{\Gamma^+}\bigr)= \biggl\{\Bigl(\widetilde{\mathcal
L}^{\Gamma_I^+}, \widetilde W_I^{\Gamma_I^+}\Bigr)\Bigm|
 I\in{\mathcal  N}^+\biggr\}\quad{\rm and}\\
&&\bigl(\widehat{\mathcal  L}^{\Gamma^+}, \widehat
V^{\Gamma^+}\bigr)=\biggl\{\Bigl(\widehat{\mathcal
L}^{\Gamma_I^+}, \widehat V_I^{\Gamma_I^+}\Bigr)\Bigm|
I\in{\mathcal
N}^+\biggr\}, \\
&& \bigl(\widetilde{\mathcal L}^{\Gamma^-}, \widetilde
W^{\Gamma^-}\bigr)= \bigg\{\Bigl(\widetilde{\mathcal
L}^{\Gamma_J^-}, \widetilde W_J^{\Gamma_J^-}\Bigr)\Bigm|
 J\in{\mathcal  N}^-\biggr\}\quad{\rm and}\\
 &&\bigl(\widehat{\mathcal  L}^{\Gamma^+},
\widehat V^{\Gamma^-}\bigr)=\biggl\{\Bigl(\widehat{\mathcal
L}_J^{\Gamma_J^+}, \widehat V_J^{\Gamma_J^-}\Bigr)\Bigm|
J\in{\mathcal N}^-\biggr\}.
\end{eqnarray*}
 In particular,
$(\widehat{\mathcal L}^{\Gamma^+}, \widehat V^{\Gamma^+})$ and
$(\widehat{\mathcal L}^{\Gamma^-}, \widehat V^{\Gamma^-})$ are
also the systems of the stratified Banach bundles. For
$i=1,\cdots, m_+$, $j=1,\cdots, m_-$ , let $\langle
f_{ij}\rangle=\langle f^+_i\sharp_y f^-_j\rangle$ and
\begin{eqnarray*}\left.\begin{array}{ll}
W_{ij}=\{\langle g^+\sharp_y g^-\rangle\,|\, \langle g^+\rangle\in
W^+_i,\; \langle g^-\rangle\in W^-_j\}\;\;{\rm and}\\
\widetilde W_{ij}=\{g^+\sharp_y g^-\,|\, g^+\in\widetilde W^+_i,\;
g^-\in\widetilde W^-_j\}.
\end{array}\right.
\end{eqnarray*}
 Then $\langle f_{ij}\rangle\in\overline{\mathcal  M}_D(\tilde y; H, J)$ and
$\{W_{ij}| 1\le i\le m_+, 1\le j\le m_-\}$ is a covering of
$\overline{\mathcal  M}_D(\tilde y; H, J)$. Moreover, each
$W_{ij}$ has a uniformizer $(\widetilde W_{ij},
\Gamma_{ij},\pi_{ij})$ with
$$\Gamma_{ij}=\Gamma^+_i\times\Gamma^-_j\quad{\rm and}\quad
\pi_{ij}(g^+\sharp_y g^-)={\rm GL}(\pi^+_i(g^+), \pi^-_j(g^-)).$$
 The action of $\Gamma_{ij}$ on
$\widetilde W_{ij}$ is given by $(\phi^+,\phi^-)\cdot(g^+,  g^-)=
(\phi^+\cdot g^+, \phi^-\cdot g^-)$. As in \S2 we may construct
the local orbifold bundle ${\mathcal  L}_{ij}\to W_{ij}$
uniformized by $\widetilde{\mathcal  L}_{ij}\to \widetilde W_{ij}$
with uniformizing group $\Gamma_{ij}$. Moreover, there exists a
natural section $\bar\partial^D_{H,J}:{\mathcal  W}\to {\mathcal
L} =\cup_{i,j}W_{ij}$ whose zero set is $\overline{\mathcal
M}_D(\tilde y;H, J)$ and that lifts to a collection
$\{\bar\partial_{J,H}^{D,ij}| 1\le i\le m_+, 1\le j\le m_-\}$ of
$\Gamma_{ij}$-equivariant sections of $\widetilde{\mathcal
L}_{ij}\to\widetilde W_{ij}$. Here $\bar\partial_{J,H}^{D,ij}$,
$\bar\partial_{J,H}^{+,i}$ and $\bar\partial_{J, H}^{-,j}$
satisfy: for $g=g^+\sharp_y g^-\in\widetilde W_{ij}$ and $\xi\in
(\widetilde{\mathcal  L}_{ij})_g=L^p_{k-1}(\wedge^{0,1}(
 g^\ast TM))$ it must hold that
$$\bar\partial_{J,H}^{D,ij}\xi=\bigl(\bar\partial_{J,H}^{+,i}(
\xi|_{\Sigma^+})\bigr)\sharp\bigl(\bar\partial_{J,H}^{-,j}(\xi|_{\Sigma^-})
\bigr)$$
 because two terms at the right side are same when
restricted to $\Sigma^+\cap\Sigma^-$. Here $\Sigma^+$ and
$\Sigma^-$ are  the domains of $g^+$ and $g^-$ respectively. For
given $f_+\in{\mathcal  B}^{p,k}_+(\tilde y; H)$ and
$f_-\in{\mathcal B}^{p,k}_-(\tilde y\sharp(-D); H)$ the elements
$\xi^+\in L^p_{k-1}(\wedge^{0,1}(f_+^\ast TM))$ and $\xi^-\in
L^p_{k-1}(\wedge^{0,1}(f_-^\ast TM))$ might not be glued into one
of $L^p_{k-1}(\wedge^{0,1}((f_+\sharp_yf_-)^\ast TM))$. Thus we
might not glue ${\mathcal  L}^+$ and ${\mathcal L}^-$ in general.

Also denote by ${\mathcal  N}={\mathcal  N}^+\times {\mathcal
N}^-$. For $(I, J)$ and $(I^\prime, J^\prime)$ in ${\mathcal  N}$
we say $(I, J)\subset (I^\prime, J^\prime)$ if $I\subset I^\prime$
and $J\subset J^\prime$. Corresponding to this covering we have
$W_{(I,J)}={\rm GL}(W_I\times W_J)$, $\Gamma_{(I,J)}=
\Gamma^+_I\oplus\Gamma^-_J$
 and the system of bundles
\begin{eqnarray*}
\hspace{-20mm}&&\hspace{15mm}(\widetilde{\mathcal  L}^{\Gamma},
\widetilde W^{\Gamma})=\Bigl\{\Bigl(\widetilde{\mathcal
L}_{(I,J)}^{\Gamma_{(I,J)}}, \widetilde
W_{(I,J)}^{\Gamma_{(I,J)}}\Bigr), \pi_{(I, J)},
\pi_{(I,J)}^{(I^\prime, J^\prime)}\,\bigm|\,
(I, J)\subset (I^\prime, J^\prime)\in{\mathcal  N}\Bigr\},\\
\hspace{-20mm}&&\hspace{15mm}(\widehat{\mathcal  L}^{\Gamma},
\widehat V^{\Gamma})= \Bigl\{\Bigl(\widehat{\mathcal
L}_{(I,J)}^{\Gamma_{(I,J)}}, \widehat
V_{(I,J)}^{\Gamma_{(I,J)}}\Bigr), \pi_{(I, J)},
\pi_{(I,J)}^{(I^\prime, J^\prime)}\,\bigm|\, (I, J)\subset
(I^\prime, J^\prime)\in{\mathcal  N}\Bigr\}.
\end{eqnarray*}

Let $f_{ij}=f^+_i\sharp_y f^-_j$ and $R(f_{ij})$ be the cokernel
of $D(\bar\partial^{D, ij}_{J,H})(f_{ij})$ in
$(\widetilde{\mathcal  L}_{ij})_{f_{ij}}$ as before.
 Take smooth cut-off functions
$\beta_\epsilon(f^+_i)$ on the domain $\Sigma^+_i$ of $f^+_i$
and $\beta_\epsilon(f^-_j)$ on that $\Sigma^-_j$ of $f^-_j$
supported outside of the $\epsilon$-neighborhood of their double points.
Then $\beta_\epsilon(f^+_i)$ and $\beta_\epsilon(f^-_j)$ naturally determine
a smooth cut-off function, denoted by
$\beta_\epsilon(f_{ij})=\beta_\epsilon(f^+_i)\sharp\beta_\epsilon(f^-_j)$,
 on the domain $\Sigma_{ij}=\Sigma^+_i\sharp\Sigma^-_j$
of $f_{ij}$ supported outside of the $\epsilon$-neighborhood of their double points.
 As in Step 1 in $\S2.2$ and (2.19) we may use these to define
corresponding spaces
$$R_\epsilon(f_{ij})\quad{\rm and}\quad R^\epsilon(\{f_{ij}\})=
\oplus^{m_+}_{i=1}\oplus^{m_-}_{j=1}R_\epsilon(f_{ij}).$$
 Now
as in (2.18) we take the smooth $\Gamma^+_i$-invariant cut-off
functions $\gamma(f^+_i)$ on $\widetilde W^+_i$ and
$\Gamma^-_j$-invariant cut-off functions $\gamma(f_j^-)$ on
 $\widetilde W^-_j$ such that for each $\nu_{ij}\in R_\epsilon(f_{ij})$
(which may be written as $\nu^+_i\sharp\nu^-_j$ for some
$\nu^+_i\in R_\epsilon(f_i^+)$ with domain $\Sigma^+_i$ and
$\nu^-_j\in R_\epsilon(f_j^-)$ with domain $\Sigma^-_j$),
$$\gamma(f^+_i)\cdot\widetilde{\nu_{ij}|_{\Sigma^+_i}}=
\gamma(f^+_i)\cdot\widetilde{\nu^+_i}\quad{\rm and}\quad
\gamma(f_j^-)\cdot\widetilde{\nu_{ij}|_{\Sigma^-_j}}=
\gamma(f^-_j)\cdot\widetilde{\nu^-_j}$$
 will give rise to the global section
$$\overline{\nu^+_i}=\overline{\nu_{ij}|_{\Sigma^+_i}}
=\Bigl\{\bigl(\overline{\nu_{ij}|_{\Sigma^+_i}}\bigr)_I=
\bigl(\overline{\nu^+_i}\bigr)_I\,|\, I\in{\mathcal  N}^+\Bigr\}$$
of $(\widehat{\mathcal  L}^+, \widehat V^+)$ and that
$$\overline{\nu^-_j}= \overline{\nu_{ij}|_{\Sigma^-_j}}=
\Bigl\{\bigl(\overline{\nu_{ij}|_{\Sigma^-_j}}\bigr)_J=\bigl(\overline{\nu^-_j}\bigr)_J\,|\,
J\in{\mathcal  N}^-\Bigl\}$$
 of $(\widehat{\mathcal  L}^-, \widehat
V^-)$ respectively. Then
$\gamma(f^+_i)\sharp\gamma(f^-_j)\cdot\widetilde{\nu_{ij}}$ (
defined by
$(\gamma(f^+_i)\sharp\gamma(f^-_j)\cdot\widetilde{\nu_{ij}})(g^+_i\sharp
g^-_j)=\gamma(f^+_i)(g^+_i)\cdot\gamma(f^-_j)(g^-_j)\cdot
\widetilde{\nu_{ij}}(g^+_i\sharp g^-_j)$) can yield a global
section $\overline{\nu_{ij}}=\{(\overline{\nu_{ij}})_{I, J}\,|\,
(I, J)\in{\mathcal  N}\}$ of $(\widehat{\mathcal  L}, \widehat
V)$. For $\delta>0$ small enough, it follows as before that for  a
generic choice of $\nu\in R_\delta^\epsilon(\{f_{ij}\})$  the
global section transversal to the zero section
$$S_{J, H}^{D,\nu}=\bigl\{ S_D^{\nu_{(I,J)}}=\bar\partial^D_{J,H}+
\bar\nu_{(I,J)}\bigm|(I, J)\in{\mathcal  N}\bigr\},$$ and thus get
the virtual moduli cycle
$$\left.\begin{array}{ll}\overline{\mathcal  M}^\nu_D(\tilde y;H, J)=\sum_{(I,J)\in{\mathcal  N}}
\frac{1}{|\Gamma_{(I,J)}|}\bigl\{\pi_{(I,J)}:{\mathcal
M}^{\nu_{(I,J)}}_D(\tilde y;H, J)\to {\mathcal
W}\bigr\}.\end{array}\right.$$
 Here ${\mathcal M}^{\nu_{(I,J)}}_D(\tilde y;H, J)=(S_D^{\nu_{(I,J)}})^{-1}(0)$. Note that
 $\nu$ can be expressed as $\oplus^{m_+}_{i=1}\oplus^{m_-}_{j=1}\nu_{ij}$ with
$\nu_{ij}=\nu^+_i\sharp\nu^-_j\in R_\epsilon(f_{ij})$,
$i=1,\cdots, m_+$ and $j=1,\cdots, m_-$. We put
$$\bar\nu^+=\sum^{m_+}_{i=1}\overline{\nu_{ij}|_{\Sigma^+_i}}=\sum^{m_+}_{i=1}
\overline{\nu^+_i}\quad{\rm and}\quad\bar\nu^-=
\sum^{m_-}_{j=1}\overline{\nu_{ij}|_{\Sigma^-_j}}=
\sum^{m_-}_{j=1}\overline{\nu^-_j}.\leqno(4.21)$$ They are the
global section of $(\widehat{\mathcal  L}^+, \widehat V^+)$ and
that of $(\widehat{\mathcal  L}^-, \widehat V^-)$ respectively. We
assert that $\bar\partial^+_{J,H}+\bar\nu^+$ (resp.
$\bar\partial^-_{J,H}+\bar\nu^-$) is also  transversal to the zero
section of $(\widehat{\mathcal L}^+, \widehat{\mathcal V}^+)$
(resp.$(\widehat{\mathcal L}^-, \widehat{\mathcal V}^-)$). We only
prove the assertion for $\bar\partial^+_{J,H}+\bar\nu^+$. Firstly,
using (4.19) it is not hard to check that $f\in {\mathcal
M}^{\nu_{(I,J)}}_D(\tilde y;H, J)$ if and only if  $f=f^+\sharp_y
f^-$ for some $f^+\in (\bar\partial^+_{J,H}+ \bar\nu^+_I)^{-1}(0)$
and $f^-\in (\bar\partial^-_{J,H}+ \bar\nu^-_J)^{-1}(0)$. Next, by
(4.21) we have
$$\bar\nu^+=\sum^{m_+}_{i=1}\overline{\nu_{ij}|_{\Sigma^+_i}}=\Bigl\{
\sum^{m_+}_{i=1}(\overline{\nu_{ij}|_{\Sigma^+_i}})_I\Bigm|
I\in{\mathcal  N}^+\Bigr\}=
\Bigl\{\sum^{m_+}_{i=1}(\overline{\nu^+_i})_I\Bigm| I\in{\mathcal
N}^+\Bigr\}.$$ For a given $f^+\in (\bar\partial^+_{J,H}+
\bar\nu^+_I)^{-1}(0)$ we  choose any $f^-\in
(\bar\partial^-_{J,H}+ \bar\nu^-_J)^{-1}(0)$ and obtain
 a $f=f^+\sharp_y f^-\in{\mathcal  M}^{\nu_{(I,J)}}_D(\tilde y;H, J)$. Note that
we can always extend any $\xi^+\in (\widehat{\mathcal
L}^+_I)_{f^+}$ into an element $\xi\in (\widehat{\mathcal
L}_{(I,J)})_f$. Since $\bar\partial^D_{J,H}+\bar\nu_{(I,J)}:
\widehat V_{(I,J)}\to\widehat{\mathcal  L}_{(I,J)}$ is transversal
to the zero section we have $\eta\in T_f\widehat V_{(I,J)}$ such
that $D(\bar\partial^D_{J,H}+\bar\nu_{(I,J)})(\eta)=\xi$. This
implies that
$D(\bar\partial^+_{J,H}+\bar\nu^+_I)(\eta|_{\Sigma^+})=\xi^+$. The
assertion is proved. In particular we get two virtual moduli
cycles
\begin{eqnarray*}
&& \overline{\mathcal M}^{\nu^+}_+(\tilde y;H,
J)=\sum_{I\in{\mathcal N}^+}
\frac{1}{|\Gamma^+_I|}\{\pi^+_I:{\mathcal M}_+^{\nu^+_I}(\tilde
y;H, J)\to
{\mathcal  W}^+\},\\
&&\overline{\mathcal  M}^{\nu^-}_-(\tilde y\sharp(-D);H, J)=
\sum_{J\in{\mathcal N}^-}\frac{1}{|\Gamma^-_J|}\{\pi^-_J:{\mathcal
M}_-^{\nu^-_J}(\tilde y\sharp(-D);H, J)\to{\mathcal  W}^-\}.
\end{eqnarray*}
Here ${\mathcal  M}_+^{\nu^+_I}(\tilde y;H,
J)=(\bar\partial^+_{J,H}+ \bar\nu^+_I)^{-1}(0)$ and ${\mathcal
M}_-^{\nu^-_J}(\tilde y\sharp(-D);H, J)=
(\bar\partial^-_{J,H}+\bar\nu^-_J)^{-1}(0)$. Now (4.19) and the
facts that $\pi_{(I,J)}=\pi^+_I\times\pi^-_J$ and $|\Gamma_{I,
J}|=|\Gamma_I^+\times\Gamma_J^-|= |\Gamma_I^+||\Gamma_J^-|$
together lead to\vspace{-2mm}
$$\overline{\mathcal  M}^\nu_D(\tilde y;H, J)=
\overline{\mathcal  M}^{\nu^+}_+(\tilde y;H, J)
\times\overline{\mathcal  M}^{\nu^-}_-(\tilde y\sharp(-D);H, J).$$
This is equivalent to (4.20).

Next we prove (4.17). The ideas are similar to the proof of
(4.16). Denote by $\tilde y_i=[y_i, w_i]$. By (4.13), for
$i=1,\cdots, t$, we take small neighborhoods $W^{(i)}_l$ centred
at $\langle f^{(i)}_l\rangle\in\overline{\mathcal M}^\nu_D(\tilde
y_i; H, J)$ in ${\mathcal  B}_D^{p,k}(\tilde y_i; H)$,
$l=l_{i-1}+1,\cdots, l_{i-1}+ l_i$ with $l_0=1$, such that
$\{W^{(i)}_l\}^{l_i}_{l=l_{i-1}+1}$ constitutes a finite covering
of $\overline{\mathcal  M}^\nu_D(\tilde y_i; H, J)$ satisfying the
requirements to construct the virtual moduli cycle. Let ${\mathcal
L}^{(i)}_l\to W^{(i)}_l$ be the local orbifold bundle uniformized
by $\widetilde{\mathcal L}^{(i)}_l\to\widetilde W^{(i)}_l$ with
uniformizing group $\Gamma^{(i)}_l$ and projection $\pi^{(i)}_l$,
and
\begin{eqnarray*}
&&\bar\partial_{J,H, D}^{(i)}: {\mathcal
W}^{(i)}=\bigcup^{l_i}_{l=l_{i-1}+1}W^{(i)}\to {\mathcal
L}^{(i)}=\bigcup^{l_i}_{l=l_{i-1}+1}{\mathcal L}^{(i)}_l,\\
&&\bar\partial_{J,H}^D: {\mathcal
W}=\bigcup^t_{i=1}\bigcup^{l_i}_{l=l_{i-1}+1}W^{(i)}\to {\mathcal
L}=\bigcup^t_{i=1}\bigcup^{l_i}_{l=l_{i-1}+1}{\mathcal L}^{(i)}_l
\end{eqnarray*}
the obvious sections whose zero sets are $\overline{\mathcal
M}^\nu_D(\tilde y_i; H, J)$ and $\overline{\mathcal M}^\nu_D(H,
J)$ respectively. Let ${\mathcal N}^{(i)}$ be the nerve of the
covering $\{W^{(i)}_l\}^{l_i}_{l=l_{i-1}+1}$ of
$\overline{\mathcal M}^\nu_D(\tilde y_i; H, J)$, $i=1,\cdots, t$,
and ${\mathcal N}$ be that of the covering
$\cup^t_{i=1}\{W^{(i)}_l\}^{l_i}_{l=l_{i-1}+1}$ of
$\overline{\mathcal M}^\nu_D(H, J)$  as before. The elements of
${\mathcal N}^{(i)}$ and those of ${\mathcal  N}$ are denoted by
$I^{(i)}$ and $I$ respectively. Correspondingly, we have the
bundle systems
\begin{eqnarray*}
&&\bigl(\widehat{\mathcal  L}^{\Gamma^{(i)}}, \widehat
V^{\Gamma^{(i)}}\bigr)= \Biggl\{\Biggl(\widehat{\mathcal
L}^{\Gamma^{(i)}_{I^{(i)}}}_{I^{(i)}}, \widehat
V^{\Gamma^{(i)}_{I^{(i)}}}_{I^{(i)}}\Biggr)\Bigm|
I^{(i)}\in{\mathcal N}^{(i)}\Biggr\},\\
&&(\widehat{\mathcal  L}^{\Gamma}, \widehat V^{\Gamma})=
\bigl\{\bigl(\widehat{\mathcal  L}^{\Gamma_I}_I, \widehat
V^{\Gamma_I}_I\bigr)\bigm|
 I\in{\mathcal  N}\bigl\}.
 \end{eqnarray*}
 As in (2.18) let $R^\epsilon(\{f^{(i)}_l\})
=\oplus^t_{i=1}\oplus^{l_i}_{l=l_{i-1}+1}R_\epsilon(f^{(i)}_l)$ be
the corresponding finite dimensional space that is used to
construct the virtual moduli cycle from $(\widehat{\mathcal
L}^\Gamma, \widehat V^\Gamma)$ and $\bar\partial^D_{J,H}$. Then
for a generic small $\nu\in R^\epsilon(\{f^{(i)}_l\})$ we obtain a
global section $\bar\nu: \widehat V^{\Gamma}\to\widehat{\mathcal
L}^{\Gamma}$ such that $S_D^\nu=\{
S_D^{\nu_I}=\bar\partial_{J,H}^D+\bar\nu_I\,|\, I\in{\mathcal
N}\}$ is transversal to the zero section. As above we may prove
that $\bar\nu$ induces a global section $\bar\nu^{(i)}: \widehat
V^{\Gamma^{(i)}}\to \widehat{\mathcal L}^{\Gamma^{(i)}}$ such that
$$S_D^{\nu^{(i)}}=\Biggl\{ S_D^{\nu^{(i)}_{I^{(i)}}}=\bar\partial_{J,H, D}^{(i)}+
\bar\nu^{(i)}_{I^{(i)}}\;\Bigm|\; I^{(i)}\in{\mathcal
N}^{(i)}\Biggr\}$$ is also transversal to the zero section for
each $i$. Let
\begin{eqnarray*}
&&\overline{\mathcal  M}^\nu_D(H, J)=\sum_{I\in{\mathcal
N}}\frac{1}{|\Gamma_I|} \{\pi_I: {\mathcal  M}_D^{\nu_I}(H,
J)\to{\mathcal W}\}\quad{\rm and}\\
&&\overline{\mathcal  M}_D^{\nu^{(i)}}(\tilde y_i; H,
J)=\sum_{I^{(i)}\in{\mathcal  N}^{(i)}}
\frac{1}{|\Gamma^{(i)}_{I^{(i)}}|} \Biggl\{\pi^{(i)}_{I^{(i)}}:
{\mathcal  M}_D^{\nu^{(i)}_{I^{(i)}}}(\tilde y_i; H, J)\to
{\mathcal W}^{(i)}\Biggr\},
\end{eqnarray*}
 $i=1,\cdots, t$, be the corresponding virtual
moduli cycles, where
$${\mathcal  M}_D^{\nu_I}(H, J)=(S_D^{\nu_I})^{-1}(0)\quad{\rm and}\quad
{\mathcal  M}_D^{\nu^{(i)}_{I^{(i)}}}(\tilde y_i; H, J)=
\Bigl(S_D^{\nu^{(i)}_{I^{(i)}}}\Bigr)^{-1}(0).$$ By (4.10) the top
strata of $\overline{\mathcal  M}^\nu_D(H, J)$ can only contain
those of  $\overline{\mathcal  M}_D^{\nu^{(i)}}(\tilde y_i; H,
J)$, $i=1,\cdots, r$. Other top strata are all empty. Thus
$$T\overline{\mathcal  M}^\nu_D(H, J)=
\cup^r_{i=1}T\overline{\mathcal  M}_D^{\nu^{(i)}}(\tilde y_i; H,
J).$$
 The union is also disjoint because $\tilde y_i$,
$i=1,\cdots, r$, are different. Since  (4.14) implies that the
intersections in (4.15) may only occur in the top strata we arrive
at
$$(\bar E_a^u\times\bar E_b^s)\cdot {\rm EV}_D^{\nu}=\sum^r_{i=1}
 (\bar E_a^u\times\bar E_b^s)\cdot {\rm EV}_D^{\nu^{(i)}}(\tilde y_i).$$
This completes the proof of (4.17). \hfill$\Box$\vspace{2mm}

\noindent{\it Step 3.}\hspace{2mm}Now we need to introduce a kind
of deformation
 spaces for understanding the right side of (4.18). For $\rho\ge 0$ and
$f:\mbox{\Bb R}\times S^1\to M$ we define $\bar\partial_{J,
H_\rho}f$ by
\begin{eqnarray*}
(4.22)\quad\left.\begin{array}{ll}
 \bar\partial_{J, H_\rho}f(s,t)=
\partial_s f(s,t)+ \\
\quad\quad J(f) (\partial_t f-(\beta_+(s+\rho+1)\cdot
\beta_+(\rho+1-s))X_{H}(t, f))=0, \end{array}\right.
\end{eqnarray*}
where $\beta_+$ is as in \S 2. For such a map $f$ we define the
energy of it by
$$
E_\rho(f):=\int^{\infty}_{-\infty}\int^1_0 |\partial_s
f|^2_{g_J}dsdt. \leqno(4.23)$$ By the removable singularity
theorem, $E_\rho(f)<+\infty$ implies that $f$ can be extended into
a smooth map from $\mbox{\Bb C}P^1\approx\{-\infty\}\cup\mbox{\Bb
R}\times S^1\cup \{\infty\}$ to $M$ which is also $J$-holomorphic
near $0$ and $\infty$. For a given $D\in\Gamma$ we  denote by
${\mathcal M}_D(H_\rho, J)$ the space of all maps  satisfying
(4.22), having finite energy and representing the class $D$. To
compactify it we introduce:\vspace{2mm}

\noindent{\bf Definition 4.4.}\hspace{2mm} Given $D\in\Gamma$ and
a  semistable ${\mathcal  F}$-curve $(\Sigma,\underline l)$ with a
unique principal component(cf. Def. 2.1), a continuous map
$f:\Sigma\to M$ is called a {\bf stable $(J, H_\rho)$-map of a
class $D$} if  it represents $D$ in the usual sense and also
satisfies:
\begin{enumerate}
\item[(1)] On the unique principal component $P$ with cylindrical coordinate
$(s,t)$, $f^P=f|_{P-\{z_-,z_+\}}$ satisfies: $\bar\partial_{J,
H_\rho} f^P=0$ and $ E_\rho(f^P)<+\infty$.
\item[(2)] The restriction $f^B_i$ of $f$ to each bubble
component $B_i$ is $J$-holomorphic, and the domain of each
homologically trivial bubble component is stable.
 \end{enumerate}\vspace{2mm}

We may also define the equivalence class of such a map. Denote by
$\overline{\mathcal  M}_D(H_\rho,  J)$ the space of all
equivalence classes of such maps. The energy of $f$ is defined by
$$
E_\rho(f)=\int^{\infty}_{-\infty}\int^1_0 |\partial_s
f^P|^2_{g_J}dsdt+ \sum_i\int_{B_i}(f^B_i)^\ast\omega.
$$
 From the
arguments in [Sch3, 4.2] it follows that
$$E_\rho(f)\le\omega(D)+ 2\max|H|.\leqno(4.24)$$
 As usual, for any $\rho\ge 0$ we may construct the virtual moduli cycle
$\overline{\mathcal  M}_{B-A}^{\nu_\rho}(H_\rho, J)$ of dimension
$2n+ 2c_1(B-A)$ corresponding to the space $\overline{\mathcal
M}_{B-A}(H_\rho,  J)$.
 Consider the evaluation
$${\rm E}^{\nu_\rho}_{B-A}:\overline{\mathcal  M}^{\nu_\rho}_{B-A}(H_\rho,  J)
 \to M\times M,\;  f\mapsto (f(z_-), f(z_+)).$$
By (4.2), for a generic $(h, g)\in{\mathcal O}(h_0)\times{\mathcal
R}$ we get a rational intersection number
$$n_{B-A}^{\nu_\rho}(a, b, h, g;  H_\rho,  J):=
( \bar E^u_a\times\bar E^s_b)\cdot {\rm
E}^{\nu_\rho}_{B-A}.\leqno(4.25)
$$
Define $(\Psi\circ\Phi)_\rho: QC_\ast(M,\omega; h, g;\mbox{\Bb
Q})\to QC_\ast(M,\omega; h, g;\mbox{\Bb Q})$ by
$$(\Psi\circ\Phi)_\rho(\langle a, A\rangle)=\sum_{\mu(\langle b, B\rangle)=
\mu(\langle a, A\rangle)}n_{B-A}^{\nu_\rho}(a, b, h, g;
 H_\rho,  J)\cdot\langle b, B\rangle.$$

\noindent{\bf Proposition 4.5.}\hspace{3mm}{\it For every $\rho\ge
0$, $(\Psi\circ\Phi)_\rho$ is a chain homomorphism.}\vspace{2mm}

\noindent{\it Proof}.\hspace{2mm}Firstly, note that if
$n_{B-A}^{\nu_\rho}(a, b, h, g;  H_\rho,  J)\ne 0$ then
$\overline{\mathcal  M}_{B-A}(H_\rho,  J)\ne\emptyset$ and thus it
follows from (4.24) that $0\le\omega(B-A)+ 2\max|H|$. Using the
fact,  as in Step1 we can prove that
 $(\Psi\circ\Phi)_\rho$ indeed maps $QC_\ast(M,\omega; h, g;\mbox{\Bb Q})$ into
itself.

Next, we show that $(\Psi\circ\Phi)_\rho$ commutes with the
boundary operator $\partial^Q$ in (1.6). It suffices to prove that
$\partial\circ (\Psi\circ\Phi)_\rho(\langle a, A\rangle)=
(\Psi\circ\Phi)_\rho\circ\partial^Q(\langle a, A\rangle)$ for each
$\langle a, A\rangle\in ({\rm Crit}(h)\times\Gamma)_k$. The direct
computation shows that the left  equals
$$
\sum_{\mu(\langle c, B\rangle)=k-1}\Biggl[
\sum_{\mu(d)=\mu(c)+1}n(d,c)n_{B-A}^{\nu_\rho}(a, d, h, g; H_\rho,
J)\Biggr] \cdot\langle c, B\rangle,$$
 and the right does
$$\sum_{\mu(\langle c, B\rangle)=k-1}\Biggl[
\sum_{\mu(b)=\mu(a)-1}n(a,b)n_{B-A}^{\nu_\rho}(b, c, h, g; H_\rho,
J)\Biggr] \cdot\langle c, B\rangle.$$
 Therefore
we only need to prove
\begin{eqnarray*}
&&\sum_{\mu(d)=\mu(c)+1}n(d,c)n_{B-A}^{\nu_\rho}(a, d,
h, g;  H_\rho,  J)\\
&&= \sum_{\mu(b)=\mu(a)-1}n(a,b)n_{B-A}^{\nu_\rho}(b, c, h, g;
H_\rho,  J).
\end{eqnarray*}
 This can be proved as
Propositions 3.5 and 3.6. In fact, by
$$\mu(a)-\mu(c)+ 2c_1(B-A)=\mu(\langle a, A\rangle)-\mu(\langle c, B\rangle)=1,$$
for a generic $(h, g)\in{\mathcal  O}(h_0)\times{\mathcal  R}$ the
fibre product
$$({\overline W}^u(a, h, g)\times{\overline W}^s(c, h, g))
\times_{\bar E^u_a\times\bar E^s_c={\rm E}_{B-A}^{\nu_\rho}}
\overline{\mathcal  M}_{B-A}^{\nu_\rho}(H_\rho, J)$$ is a
collection of compatible local cornered smooth manifolds of
dimension $1$ and with the natural orientations. Its boundary is
given by the union
\begin{eqnarray*}
&(\cup_{\mu(b)=\mu(a)-1}n(a, b)\cdot (W^u(b, h, g)
\times{\overline W}^s(c, h, g)) \times_{R_{bc}}
\overline{\mathcal  M}_{B-A}^{\nu_\rho}(H_\rho, J))\cup\\
&\hspace{-3mm}(-\cup_{\mu(d)=\mu(c)+1}n(d,c)\cdot (W^u(a, h,
g)\times{\overline W}^s(d, h, g)) \times_{R_{ad}}
\overline{\mathcal M}_{B-A}^{\nu_\rho}(H_\rho, J))
\end{eqnarray*}
because $\overline{\mathcal  M}_{B-A}^{\nu_\rho}(H_\rho, J)$ has
no components of codimension $1$. Here $R_{ad}$ and $R_{bc}$ are
representing $\bar E^u_a\times\bar E^s_d={\rm E}_{B-A}^{\nu_\rho}$
and $\bar E^u_b\times\bar E^s_c={\rm E}_{B-A}^{\nu_\rho}$.
 Hence as before the conclusions can follow from
$$\sharp\partial(({\overline W}^u(a, h, g)\times{\overline W}^s(c, h,
 g))\times_{\bar E^u_a\times\bar E^s_c={\rm E}_{B-A}^{\nu_\rho}}
\overline{\mathcal  M}_{B-A}^{\nu_\rho}(H_\rho, J))=0.$$
\hfill$\Box$

 Our purpose is to prove that $(\Psi\circ\Phi)_0$ is chain homotopy equivalent to
$\Psi\circ\Phi$. To this goal let us consider the space
$${\mathcal  M}_{B-A}(\{H_\rho\}, J)=\cup_{\rho\ge 0}\{\rho\}\times
{\mathcal  M}_{B-A}(H_\rho, J).$$
 From the arguments in
\cite{HS1, S1, S2} and \cite{Sch2} it is not hard to derive that
for any sequence $u^m\in {\mathcal  M}_{B-A}(H_{\rho_m}, J)$ with
$\rho_m\to +\infty$ there must exist a subsequence (still denoted
by $u^m$),  finitely many elements $\tilde x_1, \cdots, \tilde
x_k\in\tilde{\mathcal P}( H)$,  and $u_0\in {\mathcal  M}_+(\tilde
x_1;\theta, J)$, $u_j\in {\mathcal M}(\tilde x_j, \tilde x_{j+1};
H, J)$, $j=1,\cdots, k-1$, $u_k\in {\mathcal  M}_-(\tilde x_k; H,
J)$, and sequences $-\rho_m-1\equiv s^0_m <\cdots<
s^k_m\equiv\rho_m+1$,
 such that $u^m(s+s^i_m,t)$ {\bf converge modulo bubbling } to $u_i(s,t)$ for
$i=0,\cdots, k$ (see \cite{S2} for the precise definition of this
term).
 Moreover, if $w^i_l$, $i=1,\cdots, i_l$, are all bubbles attached to
$u_i$, then the connected sum of all  $u_i$, $w^i_l$, $l=1,\cdots,
i_l$, $i=0,\cdots, k$, represents the class $B-A$. This
convergence result  shows that the stabilized space of ${\mathcal
M}_{B-A}(\{H_\rho\}, J)$ is given by
$$\overline{\mathcal  M}_{B-A}(\{H_\rho\}, J)=
\cup_{\rho\in [0, +\infty]}\{\rho\}\times \overline{\mathcal
M}_{B-A}(H_\rho, J),\leqno(4.26)$$
where $\overline{\mathcal
M}_{B-A}(H_{+\infty}, J)$ is understood as $\overline{\mathcal
M}_{B-A}(H, J)$ in (4.11) with $D$ replaced by $B-A$, and $[0,
+\infty]=[0, +\infty)\cup\{+\infty\}$ is the compactification of
$[0, +\infty)$ equipped with the structure of a bounded manifold
obtained by  requiring that
$$h: [0, +\infty]\to [0, 1],\;t\mapsto t/\sqrt{1+t^2}$$
is a diffeomorphism. Then by the standard arguments we can prove
that the space in (4.26) is compact and Hausdorff with respect to
the weak $C^\infty$-topology.
 In addition, one has the obvious continuous evaluation map
$$\Xi_{B-A}: \overline{\mathcal  M}_{B-A}(\{H_\rho\}, J)\to
 M\times M\leqno(4.27)$$
given by $\Xi_{B-A}(\rho, \langle f\rangle)=(f(z_-), f(z_+))$ for
$\langle f\rangle\in\overline{\mathcal  M}_{B-A}(H_\rho, J)$ with
$\rho\in [0, +\infty)$, and  $\Xi_{B-A}(+\infty, \langle
f\rangle)=(f(z_-), f(z_+))$ for $\langle
f\rangle\in\overline{\mathcal  M}_{B-A}(H_{+\infty}, J)$.

 As above we can construct an associated virtual moduli
cycle $\overline{\mathcal M}^{\nu}_{B-A}(\{H_\rho\},  J) $ of
dimension $2n+ 2c_1(B-A)+1$ with the boundary
$$\partial\overline{\mathcal  M}^{\nu}_{B-A}(\{H_\rho\},  J)=
(-\overline{\mathcal  M}^{\nu_0}_{B-A}(H_0,  J))\cup
\overline{\mathcal  M}^{\nu_{+\infty}}_{B-A}(H_{+\infty},  J).$$
 Moreover, the
evaluation map in (4.27) can naturally be extended onto the
virtual module cycle, denoted by
$$\Xi^{\nu}_{B-A}:\overline{\mathcal  M}^{\nu}_{B-A}(\{H_\rho\}, J)\to
M\times M.$$
 Let $n_{B-A}^{\nu_0}(a, b, h, g;  H_0, J)$ (resp.
$n_{B-A}^{\nu_{+\infty}}(a, b, h, g;
 H_{+\infty},  J)$) denote the intersection number of $\bar E^u_a\times\bar E_b^s$
and the restriction of $\Xi^{\nu}_{B-A}$ to $\overline{\mathcal
M}^{\nu_0}_{B-A}(H_0,  J)$ (resp. $\overline{\mathcal
M}^{\nu_{+\infty}}_{B-A}(H_{+\infty},  J)$). Note that
$\overline{\mathcal  M}^{\nu_{+\infty}}_{B-A}(H_{+\infty}, J)$ is
just a virtual moduli cycle associated with the space
$\overline{\mathcal  M}_{B-A}(H, J)$ defined in (4.12). By (4.18)
we get
$$\Psi\circ\Phi(\langle a, A\rangle)=\sum_{\mu(\langle b, B\rangle)=
\mu(\langle a, A\rangle)}n_{B-A}^{\nu_{+\infty}}(a, b, h, g;
 H_{+\infty},  J)\cdot\langle b, B\rangle.\leqno(4.28)$$
As in \cite{F, SZ} and \cite{S2} we wish to define a homomorphism
$\varphi: QC_\ast(M,\omega; h, g;\mbox{\Bb Q})\to
QC_\ast(M,\omega; h, g;\mbox{\Bb Q})$ such that
$$\Psi\circ\Phi-(\Psi\circ\Phi)_0=\partial^Q\varphi +\varphi\partial^Q.\leqno(4.29)$$
For $h\in{\mathcal  O}(h_0)$, $\langle a, A\rangle\in ({\rm
Crit}(h)\times\Gamma)_k$ and $\langle d, D\rangle\in ({\rm
Crit}(h)\times\Gamma)_{k+1}$, the equality $\mu(a)-\mu(d)+
2c_1(D-A)+1=0$ implies that for a generic $(h, g)\in{\mathcal
O}(h_0)\times{\mathcal R}$ the evaluations $\bar E^u_a\times\bar
E^s_d$ and $\Xi^{\nu}_{D-A}$ are intersecting transversally. So
the rational intersection number
$$n_{D-A}(a, d, h, g;  \{H_\rho\},  J):=
( \bar E^u_a\times\bar E^s_d)\cdot \Xi^{\nu}_{D-A}$$ is
well-defined. Then it is not difficult to check that $\varphi$
defined by
$$\varphi(\langle a, A\rangle)=\sum_{\mu(\langle d, D\rangle)=\mu(\langle a, A\rangle)
+1}n_{D-A}(a, d, h, g;  \{H_\rho\},  J)\cdot\langle d,
D\rangle\leqno(4.30)$$
is an endomorphism of $QC_\ast(M,\omega; h,
g;\mbox{\Bb Q})$ and satisfies (4.29). That is,
$$\Psi\circ\Phi(\langle a, A\rangle)-(\Psi\circ\Phi)_0(\langle a, A\rangle)=
\partial^Q\varphi(\langle a, A\rangle) +\varphi\partial^Q(\langle a, A\rangle)$$
for each $\langle a, A\rangle\in{\rm Crit}(h)\times\Gamma$.
In fact, by the direct computation it suffice to prove
$$\left.\begin{array}{ll} \quad n_{B-A}^{\nu_{+\infty}}(a, b,
h,g; H_{+\infty},J)-n^{\nu_0}_{B-A}(a, b, h, g; H_0, J)\vspace{2mm}\\
=\sum_{\mu(c)=\mu(b)+1}n_{B-A}(a, c, h, g;\{H_\rho\},
J)n(c, b)\vspace{2mm}\\
- \sum_{\mu(d)=\mu(a)-1}n(a, d)n_{B-A}(d, b, h, g;\{H_\rho\},
J)\end{array} \right.\leqno(4.31)$$
 for each $\langle b, B\rangle\in{\rm
Crit}(h)\times\Gamma$ with $\mu(\langle b, B\rangle)=\mu(\langle
a, A\rangle)$. To prove it we take a generic $(h, g)\in{\mathcal
O}(h_0)\times{\mathcal  R}$ so that the evaluations $\bar
E^u_a\times\bar E^s_b$, $\bar E^u_a\times\bar E^s_c$ and $\bar
E^u_d\times\bar E^s_b$ are transversal to $\Xi^{\nu}_{B-A}$.
 By Lemmas 3.1 and 3.2, for a generic $(h, g)$ the fibre product
$$({\overline W}^u(a, h, g)\times{\overline W}^s(b, h, g))
\times_{\bar E^u_{a}\times\bar E^s_{b}=\Xi^{\nu}_{B-A}}
\overline{\mathcal  M}^{\nu}_{B-A}(\{H_\rho\},  J)$$
 is a collection of compatible local cornered smooth manifolds of
dimension $1$ and with the natural orientations. Its boundary is
given by
\begin{eqnarray*}\left.\begin{array}{ll}
\! (-({\overline W}^u(a, h, g)\times{\overline W}^s(b, h, g)
  )\times_{\bar E^u_{a}\times\bar E^s_{b}={\rm E}_{B-A}^{\nu_0}}
\overline{\mathcal  M}^{\nu_0}_{B-A}(H_0,  J))\vspace{2mm}\\
\!\bigcup(({\overline W}^u(a, h, g)\times{\overline
 W}^s(b,  h, g))\times_{\bar E^u_{a}\times\bar E^s_{b}={\rm EV}_{B-A}^\nu}
\overline{\mathcal  M}^{\nu_{+\infty}}_{B-A}(H_{+\infty},  J))\vspace{2mm}\\
\!\bigcup\bigcup_{\mu(d)=\mu(a)-1}n(a,d)\hspace{-1mm}\cdot
(({\overline W}^u(a, h, g)\times W^s(d, h, g))\!
\times_{R_{ad}}\overline{\mathcal  M}^{\nu}_{B-A}(\{H_\rho\},  J))\vspace{2mm}\\
\!\bigcup\bigcup_{\mu(c)=\mu(b)+1}n(c,b)\hspace{-1mm}\cdot
(\!-(W^u(c, h, g)\!\times\!{\overline W}^s(b, h, g))\!\!
\times_{R_{cb}} \overline{\mathcal M}^{\nu}_{B-A}(\{H_\rho\}, J)).
\end{array}\right.
\end{eqnarray*}
Here $R_{ad}=\bar E^u_{a}\times\bar E^s_{d}=\Xi^{\nu}_{B-A}$ and
$R_{cb}=\bar E^u_{c}\times\bar E^s_{b}=\Xi^{\nu}_{B-A}$.
 This implies (4.31). To sum up we have proved:\vspace{2mm}

\noindent{\bf Proposition 4.6.}\hspace{2mm}{\it $\Psi\circ\Phi$ is
chain homotopy equivalent to $(\Psi\circ\Phi)_0$.}\vspace{2mm}

\noindent{\it Step 4.}\hspace{2mm}We need to make further
homotopy. For $\tau\in [0, 1]$ we define
$$\bar\partial_{ J, \tau H_0}u(s,t)=\!\partial_s u(s,t)+ J(u)
(\partial_t u-\tau(\beta_+(s+1)\cdot\beta_+(1-s))X_{H}(t, u))
 =0.$$
In Definition 4.4 we replace $\bar\partial_{J,  H_\rho}$ with
$\bar\partial_{J, \tau H_0}$ and define the corresponding stable
$(J, \tau H_0)$-map of  class $D$. Let  $\overline{\mathcal
M}_D(\tau H_0, J)$ be the space of all equivalence classes of such
maps. For $\langle f\rangle\in\overline{\mathcal M}_{B-A}(\tau
H_0, J)$, as in (4.24) we can estimate
$$\left.\begin{array}{ll}E^\tau(f):=\!\int^1_0
|\partial_s f^P|^2_{g_J}dsdt+
\sum_i\int_{B_i}(f^B_i)^\ast\omega\le\omega(B-A)+ 2\tau
\max|H|.\end{array}\right.$$
 As above we can
construct a virtual moduli cycle $\overline{\mathcal
M}^{\nu}_{B-A}(\{\tau H_0\}, J)$ of the compact space
$\cup_{\tau\in[0, 1]}\{\tau\}\times\overline{\mathcal
M}_{B-A}(\tau H_0, J)$
 of dimension $2n+ 2c_1(B-A)+1$ and with boundary
 $$\partial\overline{\mathcal
M}^{\nu}_{B-A}(\{\tau H_0\}, J)=(-\overline{\mathcal
M}^{\nu^0}_{B-A}(0, J))\cup\overline{\mathcal
M}^{\nu^1}_{B-A}(H_0, J).$$
 (Actually, $\overline{\mathcal  M}^{\nu^1}_{B-A}(H_0, J)$ can be
 chosen as $\overline{\mathcal  M}^{\nu_0}_{B-A}(H_0, J)$.)
 For $\tau=0, 1$ and a generic $(h, g)\in{\mathcal  O}(h_0)\times{\mathcal  R}$,
        using the evaluation map
$${\rm E}^{\nu^\tau}_{B-A}:\overline{\mathcal  M}^{\nu^\tau}_{B-A}(
\tau H_0,  J)\to M\times M,\;  f\mapsto (f(z_-), f(z_+)),$$
 we get the well-defined rational
intersection number
$$n_{B-A}^{\nu^\tau}( a, b, h, g; \tau H_0,  J):=
( \bar E^u_a\times\bar E^s_b)\cdot {\rm E}_{B-A}^{\nu^\tau}.$$
  By (4.25) we have
$$n_{B-A}^{\nu^1}(a, b, h, g;  H_0,  J)=n_{B-A}^{\nu_0}(a, b, h, g;  H_0,  J)$$
since they are independent of the generic choices of $\nu$ and
$(g,h)$. Therefore, as in Step 3 we can easily prove that
$(\Psi\circ\Phi)_0$ and thus $\Phi\circ\Phi$ are chain homotopy
equivalent to $(\Psi\circ\Phi)^0$ defined by
$$(\Psi\circ\Phi)^0(\langle a, A\rangle)=
\sum_{\mu(\langle b, B\rangle)=\mu(\langle a, A\rangle)}
n_{B-A}^{\nu^0}( a, b, h, g; 0,  J)\cdot\langle b,
B\rangle.\leqno(4.32)$$
 Here as in Proposition 4.5 it can be
proved that $(\Psi\circ\Phi)^0$ is a chain homomorphism. We omit
it. Now Theorem 4.1 can follow from this and the following result:
\vspace{2mm}

\noindent{\bf Proposition 4.7.}\hspace{2mm}{\it The numbers
$n_{B-A}^{\nu^0}( a, b, h, g; 0,  J)$ at (4.32) satisfies
\begin{eqnarray*}
  n_{B-A}^{\nu^0}(a, b, h, g; 0, J)=
\left\{\begin{array}{ll}1 \; &{\rm if}\;  a= b\;{\rm and}\; A=B,\\
0\;&{\rm otherwise}.
\end{array}\right.
\end{eqnarray*}}

\noindent{\it Proof}.\quad {\it Case 1: $A\ne B$}. In this case
 $\overline{\mathcal  M}_{B-A}(0, J)$ contains no the constant maps.
 Moreover, the domain of elements of $\overline{\mathcal  M}_{B-A}(0, J)$ is only
$0$-pointed semistable ${\mathcal  F}$-curves with at least a
principal component and the operator $\bar\partial_J$ is invariant
under action of the automorphism group of the domain of a
semistable ${\mathcal  F}$-curve. Thus as in [LiuT2] the
associated virtual moduli cycle $\overline{\mathcal
M}_{B-A}^{\nu^0}(0, J)$ can be required to carry a free $\mbox{\Bb
S}^1$-action under which the evaluation is invariant. Hence for a
generic $(h, g)\in{\mathcal  O}(h_0)\times{\mathcal R}$ the fibre
product\vspace{-1mm}
$$({\overline W}^u(a, h, g)\times{\overline W}^s(b, h, g))
\times_{(\bar E^u_a\times\bar E^s_b)={\rm E}_{B-A}^{\nu^0}}
\overline{\mathcal  M}_{B-A}^{\nu^0}(0, J)$$
 must be empty under the condition (4.2). We get
 $n_{B-A}^{\nu^0}(a, b, h, g;  0,  J)=0$.

\noindent{\it Case 2: $A=B$}. Now $\overline{\mathcal M}_{B-A}(0,
J)$ can naturally be identified with $M$ and thus
$\overline{\mathcal M}^{\nu^0}_{B-A}(0, J)$ may be taken as
$\overline{\mathcal M}_{B-A}(0, J)$. Note that $\mu(a)=\mu(b)$ in
the present case and that ${\overline W}^u(a, h, g)\cap{\overline
W}^s(b, h, g)$ carries a free $\mbox{\Bb R}$-action if it is
nonempty. The conclusions follow naturally. \hfill$\Box$

Summing up the above arguments we complete the proof of Theorem
4.1.\vspace{2mm}

\noindent{\bf Remark 4.8.}\hspace{2mm}If the Morse function $h$
has only critical points of even index then $\Phi$ is a right
inverse of $\Psi$ as the chain homomorphisms between $QC_\ast(M,
\omega; h, g;\mbox{\Bb Q})$ and $C_\ast(H, J,\nu;\mbox{\Bb Q})$.
Indeed,  carefully checking the proof of Th.4.1 one will find that
$n_{B-A}^{\nu_{+\infty}}(a, b, h, g; H_{+\infty}, J)-
n^{\nu_0}_{B-A}(a, b, h, g; H_0, J)=0$ in (4.31). The same
reasoning yields $n^{\nu^1}_{B-A}(a, b, h, g; H_0,
J)=n^{\nu^0}_{B-A}(a, b, h, g; 0, J)$. By Proposition 4.7 and
(4.28) we get that $\Psi\circ\Phi(\langle a, A\rangle)=\langle a,
A\rangle$ for each $\langle a, A\rangle$.\vspace{2mm}

\noindent{\bf Theorem 4.9.}\hspace{2mm}{\it $\Phi\circ\Psi$ is
chain homotopy equivalent to the identity. Consequently, $\Phi$
induces a surjective $\Lambda_\omega$-module homomorphism from
$QH_\ast(h, g;\mbox{\Bb Q})$ to $HF_\ast(M,\omega;\\ H,
J,\nu;\mbox{\Bb Q})$.}\vspace{2mm}

\noindent{\it Proof}.\hspace{2mm}The proof is similar to that of
Theorem 4.1. We only give main steps.

\noindent{\it Step 1.}\hspace{2mm} For every
 $\tilde x\in\tilde{\mathcal  P}_k( H)$  the direct computation gives
\begin{eqnarray*}
&& \Phi\circ\Psi(\tilde x)=\sum_{\mu(\tilde y)=k}
m^\nu_{-,+}(\tilde
x,\tilde y)\cdot \tilde y,\\
&& m^\nu_{-,+}(\tilde x,\tilde y)=\sum_{\mu(\langle a,
A\rangle)=k} n^{\nu^-}_-( a, \tilde x\sharp(-A))\cdot n^{\nu^+}_+(
a,
 \tilde y\sharp(-A)).
 \end{eqnarray*}
 As in the proof of Theorem 4.1 using Proposition
3.4 we can easily prove the second sum to be finite. That is,
there are only finitely many $\langle a, A\rangle\in ({\rm
Crit}(h)\times\Gamma)_k$ such that $n^{\nu^-}_-( a, \tilde
x\sharp(-A))\cdot n^{\nu^+}_+( a, \tilde y\sharp(-A))\ne 0$.
Actually, there exist finitely many $A\in\Gamma$, saying
$A_1,\cdots, A_s$,  such that
$$\overline{\mathcal  M}_-(\tilde x\sharp(-A_i); H, J)\ne\emptyset\quad{\rm and}\quad
\overline{\mathcal  M}_+(\tilde y\sharp(-A_i); H, J)\ne\emptyset$$
for $i=1,\cdots, s$. Note that the product $n^{\nu^-}_-( a, \tilde
x\sharp(-A))\cdot n^{\nu^+}_+( a, \tilde y\sharp(-A))$ can be
explained as  the intersection number of the product evaluations
$${\rm EV}^{\nu^-}_-\times{\rm EV}^{\nu^+}_+:
\overline{\mathcal  M}^{\nu^-}_-(\tilde x\sharp(-A);H, J)\times
\overline{\mathcal  M}^{\nu^+}_+(\tilde y\sharp(-A);H, J)\to
M\times M$$ and $\bar E^s_a\times\bar E^u_a: \overline{W}^s( a, h,
g) \times\overline{W}^u( a, h, g)\to M\times M$,
 $$(\bar E^s_a\times\bar E^u_a)\cdot
({\rm EV}^{\nu^-}_-(\tilde x\sharp(-A))\times {\rm
EV}^{\nu^+}_+(\tilde y\sharp(-A)))$$
 for a generic $(h, g)\in{\mathcal  O}(h_0)\times{\mathcal  R}$, and that the fibre
product
$$(\overline{W}^s( a, h, g)\times\overline{W}^u(a, h, g))
\times_{R} (\overline{\mathcal M}^{\nu^-}_-(\tilde x\sharp(-A); H,
J)\times \overline{\mathcal M}^{\nu^+}_+(\tilde y\sharp(-A); H,
J))$$ is an empty set for a generic
 $(h, g)\in{\mathcal O}(h_0)\times{\mathcal  R}$ even if
 $\mu(\langle a, A\rangle)\ne k=\mu(\tilde x)=
\mu(\tilde y)$. Here $R=\bar E^s_{ a}\times\bar E^u_{ a}= {\rm
EV}^{\nu^-}_-\times{\rm EV}^{\nu^+}_+$.  Thus
\begin{eqnarray*}(4.33)\left\{\begin{array}{ll}
m^\nu_{-,+}(\tilde x,\tilde y)\\
\!\!=\hspace{-1mm} \sum_{\langle a, A\rangle\in{\rm
Crit}(h)\times\Gamma} (\bar E^s_{ a}\times\bar E^u_{ a})\cdot
({\rm EV}^{\nu^-}_-(\tilde x\sharp(-A))\times
{\rm EV}^{\nu^+}_+(\tilde y\sharp(-A)))\\
\!\!=\!\!\sum_{a\in{\rm Crit}(h)}\sum^s_{i=1} (\bar E^s_{
a}\times\bar E^u_{ a})\cdot ({\rm EV}^{\nu^-}_-(\tilde
x\sharp(-A_i))\times {\rm EV}^{\nu^+}_+(\tilde
y\sharp(-A_i))).\end{array}\right.
\end{eqnarray*}

\noindent{\it Step 2.}\hspace{2mm}To understand this sum, for $(h,
g)\in{\mathcal  O}(h_0)\times{\mathcal  R}$ and $\rho\ge 0$ we
denote by
$${\mathcal  M}_\rho(h, g):=\{\gamma\in C^\infty([-\rho, \rho], M)\,|\,
 \dot\gamma+\nabla_{ g}h(\gamma)=0\}$$
 It is a compact manifold of dimension $2n$ and ${\mathcal
M}_0(h, g)$ can naturally be identified with $M$. Using the gluing
techniques in the Morse homology (cf. Theorem 6.8 in \cite{Lu2}),
the natural weak compactification of the noncompact manifold
$\cup_{\rho\ge 0}{\mathcal  M}_\rho(h, g)$ of dimension $2n+1$ is
given by
$$\overline{\cup_{\rho\ge 0}{\mathcal  M}_\rho(h, g)}:=
\cup^\infty_{\rho=0}{\mathcal  M}_\rho(h, g),$$
 where
$${\mathcal  M}_\infty(h, g):=\cup_{a\in{\rm Crit}(h)}\overline{W}^s( a, h, g)\times
\overline{W}^u(a, h, g)$$ and the weak convergence of a sequence
$\{\gamma_m\}\subset{\mathcal  M}_{\rho_m}(h, g)$,
$\rho_m\to\infty$, towards a pair $(u, v)\in{\mathcal M}_\infty(h,
g)$ is understood in an obvious way (cf.[AuB],[Sch1]and [Sch4]).
The space has the structure of a manifold with corners and with
boundary
$$\partial\overline{\cup_{\rho\ge 0}{\mathcal M}_\rho(h,
g)}=(-{\mathcal  M}_0(h, g))\cup(\cup_{a\in{\rm
Crit}(h)}\overline{W}^s( a, h, g) \times\overline{W}^u(a, h,
g)).$$
 Moreover, there exists a smooth evaluation
$${\rm ev}_{h, g}: \overline{\cup_{\rho\ge 0}{\mathcal  M}_\rho(h, g)}\to M\times M
$$
 given by ${\rm ev}_{h, g}(\gamma)=(\gamma(-\rho),\gamma(\rho))$ for
 $\gamma\in {\mathcal  M}_\rho(h, g)$ such that
$${\rm ev}_{h, g}|_{\overline{W}^s( a, h, g)\times\overline{W}^u(a, h, g)}=
\bar E^s_a\times\bar E^u_a\;\;{\rm for}\;{\rm each}\; a\in{\rm
Crit}(h).$$ We also denote by ${\rm ev}_{h, g}^\rho$ the
restriction of ${\rm ev}_{h, g}$ to ${\mathcal  M}_\rho(h, g)$ for
$0\le\rho\le\infty$.
 To make further arguments we need to assume that
$${\rm EV}^{\nu^-}_-(\tilde x\sharp(-A_i))\pitchfork
{\rm EV}^{\nu^+}_+(\tilde y\sharp(-A_i)),\; i=1,\cdots, s.
\leqno(4.34)$$
These may actually be obtained for a generic small
$(\nu^-,\nu^+)\in R^-_\varepsilon\times R^+_\varepsilon$ by
increasing some points $f^-_j\in\overline{\mathcal  M}_-(\tilde
x\sharp(-A_i); H, J)$ and $f^+_j\in\overline{\mathcal M}_+(\tilde
y\sharp(-A_i); H, J)$ and enlarging the spaces $R^\pm_\varepsilon$
in the construction of the virtual module cycles.

Using these, for a generic pair $(h, g)\in{\mathcal
O}(h_0)\times{\mathcal  R}$ the fibre product
$$\overline{\cup_{\rho\ge 0}{\mathcal  M}_\rho(h, g)}
\times_{R_3} (\cup^s_{i=1}\overline{\mathcal M}^{\nu^-}_-(\tilde
x\sharp(-A_i); H, J)\times \overline{\mathcal M}^{\nu^+}_+(\tilde
y, w\sharp(-A_i); H, J))$$ is a collection of compatible local
cornered smooth manifolds of dimension $1$ and with the natural
orientations. Here $R_3$ represents $ev(h,g)={\rm
EV}^{\nu^-}_-\times{\rm EV}^{\nu^+}_+ $. Its boundary is the union
of the following four sets
$$
-{\mathcal  M}_0(h, g)\times_{R_1} (\cup^s_{i=1}\overline{\mathcal
M}^{\nu^-}_-(\tilde x\sharp(-A_i))\times \overline{\mathcal
M}^{\nu^+}_+(\tilde y\sharp(-A_i))),$$\vspace{-1mm}
$$ (\hspace{-4mm}\bigcup_{\quad a\in{\rm
Crit}(h)}\hspace{-3mm}\overline{W}^s( a, h, g)
\times\overline{W}^u(a, h, g))\times_{R_2}
(\cup^s_{i=1}\overline{\mathcal M}^{\nu^-}_-(\tilde
x\sharp(-A_i))\times \overline{\mathcal M}^{\nu^+}_+(\tilde
y\sharp(-A_i))),$$\vspace{-1mm}
$$ \bigcup_{\mu(\tilde
z)=\mu(\tilde x)-1}\hspace{-5mm} -m(\tilde x, \tilde z)\cdot
\overline{\bigcup_{\rho\ge 0}{\mathcal M}_\rho(h, g)}
\times_{R_3}\Bigl(\cup^s_{i=1}\overline{\mathcal
M}^{\nu^-}_-(\tilde z\sharp(-A_i))\times\overline{\mathcal
M}^{\nu^+}_+(\tilde y\sharp(-A_i))\Bigr),$$\vspace{-1mm}
$$ \bigcup_{\mu(\tilde
z^\prime)=\mu(\tilde y)+1}\hspace{-6mm} m(\tilde z^\prime,\tilde
y)\cdot \overline{\bigcup_{\rho\ge 0}{\mathcal  M}_\rho(h, g)}
\times_{R_3} \Bigl(\cup^s_{i=1}\overline{\mathcal
M}^{\nu^-}_-(\tilde x\sharp(-A_i))\times \overline{\mathcal
M}^{\nu^+}_+(\tilde z^\prime\sharp(-A_i))\Bigr).$$
 Here $R_1$ and $R_2$   represent
 ${\rm ev}_{h,g}^0={\rm EV}^{\nu^-}_-
\times{\rm EV}^{\nu^+}_+$ and ${\rm ev}_{h,g}^\infty= {\rm
EV}^{\nu^-}_-\times{\rm EV}^{\nu^+}_+$ respectively,  $m(\tilde x,
\tilde z)=\sharp(C(\overline{\mathcal  M}^\nu(\tilde x, \tilde
z)))$ and $m(\tilde z^\prime, \tilde y)=
\sharp(C(\overline{\mathcal M}^\nu(\tilde z^\prime, \tilde y)))$.
Hereafter we omit $H,J$ in $\overline{\mathcal
M}^{\nu^\pm}_\pm(\cdot; H,J)$ without confusions. Note that
$$\sharp\Bigl(\!\!\Bigl(\hspace{-5mm}\bigcup_{\quad\;
a\in{\rm Crit}(h)}\hspace{-6mm}\overline{W}^s( a, h, g)
\times\overline{W}^u(a, h, g)\Bigr)\times_{{\mathcal  R}_2}
\Bigl(\bigcup^s_{i=1}\overline{\mathcal  M}^{\nu^-}_-\!(\tilde
x\sharp(-A_i))\times \overline{\mathcal  M}^{\nu^+}_+\!(\tilde
y\sharp(-A_i))\Bigr)\!\!\Bigr)$$
 is exactly the number in (4.33), and that because of (4.34),
$$\sharp({\mathcal  M}_0(h, g)\times_{{\mathcal  R}_1}
(\cup^s_{i=1}\overline{\mathcal  M}^{\nu^-}_-(\tilde
x\sharp(-A_i))\times \overline{\mathcal  M}^{\nu^+}_+(\tilde
y\sharp(-A_i))))$$ is exactly the sum of the intersection numbers
$$\left.\begin{array}{ll}\sum^s_{i=1}{\rm EV}^{\nu^-}_-(\tilde x\sharp(-A_i))
\cdot{\rm EV}^{\nu^+}_+(\tilde
y\sharp(-A_i))\end{array}\right.\leqno(4.35)$$
 As in the proof of Theorem 4.1  the above boundary relations
lead to:\vspace{2mm}

\noindent{\bf Proposition 4.10.}\hspace{2mm}{\it $\Phi\circ\Psi$
is chain homotopy equivalent to the homomorphism defined by
$$\overline{\Phi\circ\Psi}(\tilde x)=\sum_{\mu(\tilde y)=\mu(\tilde x)}
(\sum^s_{i=1}{\rm EV}^{\nu^-}_-(\tilde x\sharp(-A_i))\cdot{\rm
EV}^{\nu^+}_+(\tilde y\sharp(-A_i)))\cdot \tilde y.\leqno(4.36)$$}

Now we also need to prove that $\overline{\Phi\circ\Psi}$ is a
chain homomorphism from $C_\ast(H, J;\mbox{\Bb Q})$ to itself yet.
The ideas are same as those of Proposition 4.5. Let us outline it
as follows. Firstly, (2.6) implies that if the sum in (4.35) is
not zero then
$${\mathcal  F}_H(\tilde y)\le \max|H|-\min_{1\le i\le s}\omega(A_i)\;\;{\rm and}\;\;
{\mathcal  F}_H(\tilde x)\ge -\max|H|-\max_{1\le i\le
s}\omega(A_i).$$
From these it easily follows that
$\overline{\Phi\circ\Psi}$ maps $C_\ast(H, J;\mbox{\Bb Q})$ to
itself. Next, by the direct computation
 we easily reduce the proof of
$\overline{\Phi\circ\Psi}\circ\partial^F=\partial^F\circ\overline{\Phi\circ\Psi}$
to proving that for given $\tilde x\in\tilde{\mathcal P}_k(H)$ and
 $\tilde z\in\tilde{\mathcal  P}_{k-1}(H)$ the following holds.
\begin{eqnarray*}
&&\sum^s_{i=1}\sum_{\mu(\tilde y)=k}({\rm EV}^{\nu^-}_-(\tilde
x\sharp(-A_i))\cdot {\rm EV}^{\nu^+}_+(\tilde y\sharp(-A_i)))\cdot
n_i(\tilde y, \tilde
z)\\
&&=\sum^s_{i=1}\sum_{\mu(\tilde z^\prime)=k-1} n_i(\tilde x,
\tilde z^\prime)\cdot({\rm EV}^{\nu^-}_-(\tilde
z^\prime\sharp(-A_i))\cdot {\rm EV}^{\nu^+}_+(\tilde
z\sharp(-A_i))).
\end{eqnarray*}
Here and in the following unions  $n_i(\tilde x, \tilde
z^\prime)=\sharp(C(\overline{\mathcal M}^\nu(\tilde x\sharp(-A_i),
 \tilde z^\prime\sharp(-A_i))))$ and
$n_i(\tilde y, \tilde z)=\sharp(C(\overline{\mathcal M}^\nu(\tilde
y \sharp(-A_i), \tilde z\sharp(-A_i))))$. In fact, by (4.34) the
fibre product
$$\overline{\mathcal  M}^{\nu^-}_-(\tilde x\sharp(-A_i))
\times_{{\rm EV}^{\nu^-}_-={\rm EV}^{\nu^+}_+} \overline{\mathcal
M}^{\nu^+}_+(\tilde z\sharp(-A_i))$$ is a collection of compatible
local cornered smooth manifolds of dimension $1$ and with the
natural orientations. Its boundary is given by
\begin{eqnarray*}
&&(\cup_{\mu(\tilde z^\prime)=k-1}n_i(\tilde x, \tilde
z^\prime)\cdot (\overline{\mathcal  M}^{\nu^-}_-(\tilde z^\prime,
\sharp(-A_i))\times_{ {\rm EV}^{\nu_-}_-={\rm EV}^{\nu^+}_+}
\overline{\mathcal  M}^{\nu^+}_+(\tilde
z\sharp(-A_i))))\\
&&\cup(-\cup_{\mu(\tilde y)=k} (\overline{\mathcal
M}^{\nu^-}_-(\tilde x\sharp(-A_i))\times_{ {\rm EV}^{\nu^-}_-={\rm
EV}^{\nu^+}_+} \overline{\mathcal M}^{\nu^+}_+(\tilde
y\sharp(-A_i)))\cdot n_i(\tilde y, \tilde z)).
\end{eqnarray*}
 Notice that the sum in (4.35) is exactly
$$\sharp(\cup^s_{i=1}\overline{\mathcal  M}^{\nu^-}_-(\tilde x\sharp(-A_i))
\times_{{\rm EV}^{\nu^-}_-={\rm EV}^{\nu^+}_+} \overline{\mathcal
M}^{\nu^+}_+(\tilde y\sharp(-A_i))).$$
 They together lead to the conclusions. \vspace{2mm}

\noindent{\it Step 3.}\hspace{2mm}To understand the number in
(4.35) we introduce:\vspace{2mm}

\noindent{\bf Definition 4.11.}\hspace{2mm} Given  $\tilde x,
\tilde y\in\tilde{\mathcal  P}( H)$ and a semistable
 ${\mathcal  F}$-curve $(\Sigma,\underline l)$ with at least two
 principal components (cf. Def.2.1), a continuous map
$$f:\Sigma\setminus\{z_1,\cdots, z_{i-1}, z_{i+1},
\cdots, z_{N_p+1}\}\to M$$ is called a {\bf stable $(J, H)$-broken
trajectory} if we divide $(\Sigma,\underline l)$ into two
semistable ${\mathcal  F}$-curve $(\Sigma_-,\underline l^-)$ and
$(\Sigma_+,\underline l^+)$ at some double point $z_i$ between two
principal components, $1<i<N_p+1$, then $f|_{\Sigma_-}$ and
$f|_{\Sigma_+}$ are the stable $( J, H)_-$-disk with cap $\tilde
x$ and stable $(J, H)_+$-disk with cap $\tilde y$
respectively.\vspace{2mm}

We can also define its equivalence class in an obvious way.
 Denote the space of all such
equivalence classes by
 $\overline{\mathcal M}_{-+}(\tilde x, \tilde y; H, J)$.
 Let ${\mathcal M}_{-+}(\tilde x, \tilde y; H, J)$
 be its subspace consisting of those elements whose domains
have only two principal components and have no any bubble
components. Then the former is the natural compactification of the
latter, and has the virtual dimension  $\mu(\tilde x)-\mu(\tilde
y)$. We can, as before, construct an associated virtual module
cycle $\overline{\mathcal M}^\nu_{-+}(\tilde x, \tilde y; H, J)$
of dimension $\mu(\tilde x)-\mu(\tilde y)$. Specially, if
$\mu(\tilde x)=\mu(\tilde y)$ this virtual moduli cycle determines
a well-defined  rational number $n^\nu_{-+}(\tilde x, \tilde y; H,
J)$ which is independent of a generic choice of $\nu$ in the
obvious way. As in the proof of Theorem 4.1, by carefully checking
the construction of the virtual module cycle we can prove that
$$n^\nu_{-+}(\tilde x, \tilde y; H, J)={\rm EV}^{\nu^-}_-(\tilde x)\cdot
{\rm EV}^{\nu^+}_+(\tilde y).$$ Moreover, as showed in Step 1
there exist only finitely many $A_i\in\Gamma$ such that
$\overline{\mathcal  M}_{-+}(\tilde x\sharp(-A_i), \tilde
y\sharp(-A_i); H, J)\ne\emptyset$, $i=1,\cdots, s$. So the virtual
module cycle  associated with $\cup^s_{i=1}\overline{\mathcal
M}_{-+}(\tilde x\sharp(-A_i), \tilde y\sharp(-A_i); H, J)$ can be
taken as
$$\cup^s_{i=1}\overline{\mathcal  M}^{\nu_i}_{-+}(\tilde x\sharp(-A_i),
 \tilde y\sharp(-A_i); H, J).$$
It follows that
$$\sharp(\cup^s_{i=1}
\overline{\mathcal  M}_{-+}(\tilde x\sharp(-A_i), \tilde
y\sharp(-A_i); H, J))^\nu =\sum^s_{i=1}n^\nu_{-+}(\tilde
x\sharp(-A_i), \tilde y\sharp(-A_i); H, J).$$
 Thus (4.36) becomes
$$\overline{\Phi\circ\Psi}(\tilde x)=\sum_{\mu(\tilde y)=\mu(\tilde x)}
\sharp(\cup^s_{i=1} \overline{\mathcal  M}_{-+}(\tilde
x\sharp(-A_i), \tilde y\sharp(-A_i);
 H, J))^\nu \cdot \tilde y.\leqno(4.37)$$

\noindent{\it Step 4.}\hspace{2mm}Furthermore, for $\rho\ge 0$ and
$f:\mbox{\Bb R}\times S^1\to M$ we define $\bar\partial_{ J,
H^\rho}f$ by
$$\bar\partial_{ J, H^\rho}f(s,t)=
\partial_s f(s,t)+ J(f) (\partial_t f-(\beta_+(s-\rho)+ \beta_+(-s-\rho))X_{ H}(t,
f))
 =0,$$
where $\beta_+$ is as in (4.22). For such a map $f$ we still
define the energy of it by (4.23).\vspace{2mm}

\noindent{\bf Definition 4.12.}\hspace{2mm} Given $[x, v], [y,
w]\in\tilde{\mathcal  P}(H)$ and a  semistable ${\mathcal
F}$-curve $(\Sigma,\underline l)$ as in Definition 2.1, a
continuous map $f:\Sigma\setminus\{z_1,\cdots,z_{N_p+1}\}\to M$ is
called a {\bf stable $( J, H)^\rho$-trajectory}  if there exist
$[x, v]=[{ x}_1,{ u}_1],\cdots, [{ x}_{N_p+1},{ u}_{N_p+1}]=[ y,
w]$ such that  (1), (2), (3) in Def.2.1 are satisfied unless (i)
in Def.2.1(1) is replaced by $\bar\partial_{J, H^\rho}f^P_k=0$ and
$E_\rho(f_k^P)<+\infty$ in some principal component
$P_k$.\vspace{2mm}

Let  $\overline{\mathcal  M}(\tilde x, \tilde y; J, H^\rho)$
denote the space of  equivalence classes of all stable
 $(J, H)^\rho$-trajectories from $\tilde x$ to $\tilde y$. As in \S2 we can
prove that this space is compact according to the weak
$C^\infty$-topology  and construct the corresponding virtual
moduli cycle $\overline{\mathcal  M}^\nu(\tilde x, \tilde y; J,
H^\rho)$
 of dimension $\mu(\tilde x)-\mu(\tilde y)$.
 Specially, if $\mu(\tilde x)=\mu(\tilde y)$ we can associate a rational
number to it, denoted by $m^\nu(\tilde x, \tilde y;  J, H^\rho)$.
Consider the space $\cup_{\rho\ge 0}\overline{\mathcal M}(\tilde
x, \tilde y;  J, H^\rho)$. As before one easily shows that the
natural weak compactification of it is given by
$$(\cup_{\rho\ge 0}\overline{\mathcal  M}(\tilde x, \tilde y;
 J, H^\rho))\cup(\cup^s_{i=1}\overline{\mathcal  M}_{-+}(\tilde x\sharp(-A_i),
  \tilde y\sharp(-A_i); H, J)).$$
Using it we can construct a virtual module cycle of dimension $1$
and prove:\vspace{2mm}

\noindent{\bf Proposition 4.13.}\hspace{2mm}{\it
$\overline{\Phi\circ\Psi}$ in (4.37) is chain homotopy equivalent
to the homomorphism  given by
$$({\Phi\circ\Psi})^0(\tilde x)=\sum_{\mu(\tilde y)=\mu(\tilde x)}
m^\nu(\tilde x, \tilde y;  J, H^0)\cdot \tilde y.$$}

Here we have assumed that $({\Phi\circ\Psi})^0$ is a chain
homomorphism from $C_\ast(H, J;\mbox{\Bb Q})$ to itself. It can be
proved as above. We omit it. As in \cite{F, SZ} and \cite{LiuT1},
by taking a regular homotopy from  $( J, H)$ to $(J,
(\beta_+(\cdot)+\beta_+(-\cdot)) H)$ we can prove that
$({\Phi\circ\Psi})^0$ is chain homotopy equivalent to the
homomorphism defined by
$$\overline{({\Phi\circ\Psi})}^0(\tilde x)=\sum_{\mu(\tilde y)=\mu(\tilde x)}
\sharp(C(\overline{\mathcal  M}^\nu(\tilde x, \tilde y;  J,
H))))\cdot \tilde y. $$ Since $\mu(\tilde y)=\mu(\tilde x)$ it is
easily checked that $\sharp(C(\overline{\mathcal  M}^\nu(\tilde x,
\tilde y; J, H)))=1$ as $\tilde x=\tilde y$, and
$\sharp(C(\overline{\mathcal M}^\nu(\tilde x, \tilde y; J, H)))=0$
otherwise. That is, $\overline{({\Phi\circ\Psi})}^0=id$. This
fact, Proposition 4.13, (4.37) and Proposition 4.10 together prove
Theorem 4.9.\hfill$\Box$

\end{document}